\documentclass[12pt]{article}
\usepackage[utf8]{inputenc}
\usepackage{amsthm,amsmath,amssymb}
\usepackage{algpseudocode}
\usepackage{algorithm}
\usepackage{bbm}
\usepackage{setspace}

\usepackage{xcolor}
\usepackage{natbib}
\usepackage[colorlinks=true, citecolor={blue}]{hyperref}
\usepackage[margin=1in]{geometry}

\usepackage{graphicx}
\usepackage{wrapfig}

\usepackage{subcaption}

\renewcommand{\tilde}{\widetilde}
\renewcommand{\hat}{\widehat}

\DeclareMathOperator*{\argmin}{arg\,min}
\DeclareMathOperator*{\argmax}{arg\,max}

\theoremstyle{plain}  
\newtheorem{theorem}{Theorem}[section]
\newtheorem{lemma}{Lemma}
\newtheorem{proposition}{Proposition}
\newtheorem{corollary}{Corollary}

\theoremstyle{definition}

\newtheorem{example}{Example}
\newtheorem{remark}{Remark}[section]

\title{Choosing the $p$ in $L_p$ loss: rate adaptivity on the symmetric location problem}
\author{Yu-Chun Kao, Min Xu\thanks{Corresponding author: Min Xu, Department of Statistics, Rutgers University, New Brunswick, NJ, USA (email: \texttt{mx76@stat.rutgers.edu})}, and Cun-Hui Zhang \\ 
Department of Statistics \\ Rutgers University, New Brunswick, NJ, USA
}
\date{\today}

\begin{document}

\begin{spacing}{1.1}
\maketitle

\abstract{
Given univariate random variables $Y_1, \ldots, Y_n$ distributed uniformly on $[\theta_0 - 1, \theta_0 + 1]$, the sample midrange $\frac{Y_{(n)}+Y_{(1)}}{2}$ is the MLE for the location parameter $\theta_0$ and has estimation error of order $1/n$, which is much smaller compared with the $1/\sqrt{n}$ error rate of the usual sample mean estimator. However, the sample midrange performs poorly when the data has say the Gaussian $N(\theta_0, 1)$ distribution, with an error rate of $1/\sqrt{\log n}$. In this paper, we propose an estimator of the location $\theta_0$ with a rate of convergence that can, in many settings, adapt to the underlying distribution which we assume to be symmetric around $\theta_0$ but is otherwise unknown. When the underlying distribution is compactly supported, we show that our estimator attains a rate of convergence of $n^{-\frac{1}{\alpha}}$ up to polylog factors, where the rate parameter $\alpha$ can take on any value in $(0, 2]$ and depends on the moments of the underlying distribution. Our estimator is formed by minimizing the $L_\gamma$-loss with respect to the data, for a power $\gamma \geq 2$ chosen in a data-driven way -- by minimizing a criterion motivated by the asymptotic variance. Our approach can be directly applied to the regression setting where $\theta_0$ is a function of observed features and motivates the use of $L_\gamma$ loss function with a data-driven $\gamma$ in certain settings. 
}

\end{spacing}

\section{Introduction}
\label{sec:intro}

Given random variables $Y_1, \ldots, Y_n \stackrel{d}{\sim} \text{Uniform}(\theta_0 - 1, \theta_0 + 1)$, the optimal estimator for the center $\theta_0$ is not the usual sample mean $\bar{Y}$ but rather the sample midrange $Y_{\text{mid}}=\frac{Y_{(n)} + Y_{(1)}}{2}$, which is also the MLE. Indeed,
we have that
\begin{align*}
\mathbb{E} \biggl| \frac{Y_{(n)}+Y_{(1)}}{2}  - \theta_0 \biggr| &\leq \mathbb{E} \biggl| \frac{Y_{(n)} - \theta_0 - 1}{2} \biggr| + \mathbb{E} \biggl| \frac{Y_{(1)} - \theta_0 + 1}{2} \biggr| \\
&= 
1 - \mathbb{E} \frac{Y_{(n)} - \theta_0}{2} + \mathbb{E} \frac{Y_{(1)} - \theta_0}{2} = \frac{2}{n+1},
\end{align*}
which is far smaller than the $1/\sqrt{n}$ error of the sample mean; a two points argument in \citet{lecam1973} shows that the $1/n$ rate is optimal in this case. One may also show via the Lehman--Scheffe Theorem that sample midrange is the uniformly minimum variance unbiased (UMVU) estimator. However, sample midrange is a poor choice when $Y_1, \ldots, Y_n \stackrel{d}{\sim} N(\theta_0, 1)$, where we have that $\mathbb{E} | Y_{\text{mid}} - \theta_0|$ is of order $1/\sqrt{\log n}$. These observations naturally motivate the following question: let $p$ be a univariate density symmetric around $0$ and suppose $Y_1, \ldots, Y_n$ has the distribution $p(\,\cdot - \theta_0)$ which is the location shift of $p$, can we construct an estimator of the location $\theta_0$ whose rate of convergence adapts to the unknown underlying distribution $p$? 

This question has not yet been addressed by the wealth of existing knowledge, dating back to at least \cite{stein1956efficient}, on symmetric location estimation, which focuses on semiparametric efficiency for asymptotically Normal estimators. The classical theory states that when the underlying density $p$ is regular in the sense of being differentiable in quadratic mean (DQM), there exists $\sqrt{n}$-consistent estimator which has the same asymptotic variance as the best estimator when one \emph{does know} the underlying density $p$; in other words, under the regular regime and in terms of asymptotic efficiency, one can perfectly adapt to the unknown distribution. The adaptive estimators rely on being able to consistently estimate the unknown density at an appropriate rate. 

In contrast, the setting where $\theta_0$ can be estimated at a rate faster than $\sqrt{n}$ is \emph{irregular} in that the Fisher information is infinity and any $\sqrt{n}$-asymptotically Normal estimator is suboptimal; the underlying distribution is not DQM and is difficult to estimate. Even the problem of choosing between only the sample mean $\bar{Y}$ and the sample midrange $Y_{\text{mid}}$ is nontrivial, as we show in this paper that tried-and-true method of cross-validation fails in this setting (see Remark~\ref{rem:cv_fail} for the detailed discussion). 

If the underlying density $p$ is known, the optimal rate in estimating the location $\theta_0$ is governed by how quickly the function $\Delta \mapsto H\bigl( p(\cdot), p(\cdot - \Delta) \bigr)$ decreases as $\Delta$ goes to zero, where $H(p, q) := \bigl\{ \int (\sqrt{p(x)} - \sqrt{q(x)})^2 dx \bigr\}^{1/2}$ is the Hellinger distance. To be precise, for any estimator $\hat{\theta}$, we have 
\[
\liminf_{n \rightarrow \infty} \sup_{\theta_0} \mathbb{E}_{\theta_0} \bigl\{ \sqrt{n} H\bigl( p(\cdot - \theta_0), p(\, \cdot - \hat{\theta}) \bigr)\bigr\} > 0,
\]
where the supremum can be taken in a local ball of shrinking radius around any point in $\mathbb{R}$; see for example Theorem 6.1 of Chapter I of \cite{ibragimov2013statistical} for an exact statement. \cite{lecam1973} also showed that the MLE attains this convergence rate under mild conditions. Therefore, if $H^2\bigl( p(\cdot), p(\cdot - \Delta) \bigr)$ is of order $|\Delta|^{\alpha}$ for some $\alpha > 0$, then the optimal rate of the error $\mathbb{E}_{\theta_0} | \hat{\theta} - \theta_0 |$ is $n^{-\frac{1}{\alpha}}$. If the underlying density $p$ is DQM, then we have that $\alpha = 2$ which yields the usual rate of $n^{-\frac{1}{2}}$. But, if $p$ is the uniform density on $[-1, 1]$, we have $\alpha = 1$ which gives an optimal rate of $n^{-1}$.

The behavior of the function $\Delta \mapsto H\bigl( p(\cdot), p(\cdot - \Delta) \bigr)$ depends on the smoothness of the underlying density $p$. In the extreme case where $p$ has a Dirac delta point mass at 0 for instance, $H\bigl( p(\cdot), p(\cdot - \Delta) \bigr)$ is bounded away from 0 for any $\Delta > 0$. This is expected since, in this case, we can estimate $\theta_0$ perfectly by localizing the discrete point mass. More generally, discontinuities in the density function or singularities in its first derivative anywhere can increase $H\bigl( p(\cdot), p(\cdot - \Delta) \bigr)$ and thus lead to a faster rate in estimating the location $\theta_0$. Interested readers can find a detailed discussion and a large class of examples in Chapter VI of \cite{ibragimov2013statistical}. 

When the underlying density $p$ is unknown, it becomes unclear how to design a rate adaptive location estimator. One possible approach is to nonparametrically estimate $p$, but we would need our density estimator to be able to accurately recover the points of discontinuities in $p$ or singularities in $p'$ -- this goes beyond the scope of existing theory on nonparametric density estimation which largely deals with estimating a smooth density $p$. Because of the clear difficulty in analyzing rate adaptive location estimation problem in its fullest generality, we focus on rate adaptivity among compactly supported densities which exhibit discontinuity or singularity at the boundary points of the support; the uniform density on $[-1, 1]$ for instance has discontinuity at the boundary points $-1$ and $1$. 

With the more precise goal in mind, we study a simple class of estimators of the form $\hat{\theta}_{\gamma} = \argmin_{\theta} \sum_{i=1}^n |Y_i - \theta|^\gamma$ where the power $\gamma \geq 2$ is selected in a data-driven way. Estimators of this form encompass both the sample mean $\bar{Y}$, with $\gamma = 2$, and the sample midrange, with $\gamma \rightarrow \infty$. These estimators are easy to interpret, easy to compute, and can be extended in a straightforward way to the regression setting where $\theta_0$ is a linear function of some observed covariates. 

The key step is selecting the optimal power $\gamma$ from the data; in particular, $\gamma$ must be allowed to diverge with $n$ in order for the resulting estimator to have an adaptive rate. Since $\hat{\theta}_{\gamma}$ is unbiased for any $\gamma \geq 2$, the ideal selection criterion is to minimize the variance. In this work, we approximate the variance of $\hat{\theta}_{\gamma}$ by its \emph{asymptotic variance}, which has a finite sample empirical analog that can be computed from the empirical central moments of the data. We then select $\gamma$ by minimizing the empirical asymptotic variance, using Lepski's method to ensure that we consider only those $\gamma$'s for which the empirical asymptotic variance is a good estimate of the population version. For any distribution with a finite second moment, the resulting estimator has rate of convergence at least as fast as $\tilde{O}( n^{-1/2} )$, where we use the $\tilde{O}(\cdot)$ notation to suppress log-factors. Moreover, for any compacted supported density $p$ that satisfies a moment condition of the form $\int |z|^\gamma p(z) \, dz \asymp \gamma^{-\alpha}$ for some $\alpha \in (0, 2]$, our estimator attains an adaptive rate of $\tilde{O}(n^{ - \frac{1}{\alpha}})$. 

Our estimation procedure can be easily adapted to the linear regression setting where we have $Y_i = X_i^\top \beta_0 + Z_i$ where $Z_i$ has a distribution symmetric around 0. It is computationally fast using second order methods and can be directly applied on real data. Importantly, it is robust to violation of the symmetry assumption. More precisely, if $Y_i = \theta_0 + Z_i$ and the noise $Z_i$ has a distribution that is \emph{asymmetric} around 0 but still has mean zero, then our estimator will converge to $\mathbb{E} Y_i = \theta_0$ nevertheless. 

The rest of our paper is organized as follows: we finish Section~\ref{sec:intro} reviewing existing work and defining commonly used notation. In Section~\ref{sec:method}, we formally define the problem and our proposed method; we also show that our proposed method has a rate of convergence that is at least $\tilde{O}(n^{-\frac{1}{2}})$ (Theorem~\ref{thm:basic}). In Section~\ref{sec:theory}, we prove that our proposed estimator has an adaptive rate of convergence $\tilde{O}(n^{-\frac{1}{\alpha}})$ where $\alpha \in (0, 2]$ is determined by a moment condition on the noise distribution. We perform empirical studies in Section~\ref{sec:empirical} and conclude with a discussion of open problems in Section~\ref{sec:discussion}.


\subsection{Literature review}

Starting from the seminal paper by \cite{stein1956efficient}, a long series of work, for example~\cite{stone1975adaptive},~\cite{beran1978efficient}, and many others~\citep{van1970efficiency, bickel1982adaptive,
schick1986asymptotically, mammen1997optimal, dalalyan2006penalized} showed, under the regular DQM setting, we can attain an
asymptotically efficient estimator $\hat{\theta}$ by taking a pilot
estimator $\hat{\theta}_{\text{init}}$, applying a density estimation method on the
residues $\tilde{Z}_i = Y_i - \hat{\theta}_{\text{init}}$ to obtain
a density estimate $\hat{p}$, and then construct $\hat{\theta}$ either by maximizing the estimated log-likelihood, by taking one Newton step using an estimate of the Fisher information, or by various other related schemes; see \cite{bickel1993efficient} for more discussion on adaptive efficiency. Interestingly, \cite{laha2021adaptive} recently showed that the smoothness assumption can be substituted by a log-concavity condition instead.

Also motivated in part by the contrast between sample midrange and sample mean, \cite{baraud2017new} and \cite{baraud2018rho} propose the $\rho$-estimator. When the underlying density $p$ is known, the $\rho$-estimator has optimal rate in estimating the location. When $p$ is unknown, the $\rho$-estimator would need to estimate $p$ nonparametrically; it is not clear under what conditions it would attain adaptive rate. Moreover, computing the $\rho$-estimator in practice is often difficult. 

Our estimator is related to methods in robust statistics \citep{huber2011robust}, although our aim is different. Our asymptotic variance based selector can be seen as a generalization of a procedure proposed by \cite{lai1983adaptive}, which uses the asymptotic variance to select between the sample mean and the median. Another somewhat related line of work is that of \cite{chierichetti2014learning} and \cite{pensia2019estimating}, which study location estimation when $Z_1, \ldots, Z_n$ are allowed to have different distributions, all of which are still symmetric around 0, and construct robust estimators that interestingly adapt to the heterogeneity of the distributions of the $Z_i$'s. 

\subsection{Notation}
We write $[n] := \{1,2,\ldots, n\}$. We write $a\wedge b:=\min(a,b)$, $a\vee b:=\max(a,b)$, $(a)_+:=a\vee 0$ and $(a)_-:=-(a\wedge 0)$. For two functions $f, g$, we write $f \gtrsim g$ if there exists a universal constant $C > 0$ such that $f \geq C g$; we write $f \asymp g$ or $f \propto g$ if $f \gtrsim g$ and $g \gtrsim f$. We use $C$ to denote a positive universal constants whose value may be different from instance to instance. We use the $\tilde{O}(\cdot)$ notation to represent rate of convergence ignoring poly-log factors.




\section{Method}
\label{sec:method}
We observe random variables $Y_1, \ldots, Y_n$ such that
\[
Y_i = \theta_0 + Z_i \text{ for $i \in [n]$}
\]
where $\theta_0 \in \mathbb{R}$ is the unknown location and $Z_1, \ldots, Z_n \stackrel{d}{\sim} P$ where $P$ is an unknown distribution with density $p(\cdot)$ symmetric around zero. Our goal is to estimate $\theta_0$ from the observations $ Y_1, \ldots, Y_n$.



\subsection{A simple class of estimators}

Our approach is motivated by the fact that both the sample mean and the sample midrange minimize the $\ell^\gamma$ norm of the residual for different values of $\gamma$. More precisely,  
\begin{align*}
\bar{Y} &:= \frac{1}{n} \sum_{i=1}^n Y_i = \argmin_{\theta \in \mathbb{R}} \sum_{i=1}^n |Y_i - \theta|^2, \quad \text{ and }\\
Y_{\text{mid}} &:= \frac{Y_{(n)} + Y_{(1)}}{2} = \argmin_{\theta \in \mathbb{R}} \max_{i\in[n]}|Y_i-\theta| = \lim_{\gamma \rightarrow \infty} \argmin_{\theta \in \mathbb{R}}  \sum_{i=1}^n |Y_i - \theta|^\gamma.
\end{align*}

This suggests an estimation scheme where we first select the power $\gamma \geq 2$ in a data-driven way and then output the empirical center with respect to the $\ell^\gamma$ norm:
\[
\hat{\theta}_\gamma := \argmin_{\theta \in \mathbb{R}} \sum_{i=1}^n |Y_i - \theta|^\gamma.
\]
It is clear that $\bar{Y} = \hat{\theta}_2$ and that $\hat{\theta}_{\gamma}$ approaches $Y_{\text{mid}}$ as $\gamma$ increases, that is, $Y_{\text{mid}} \equiv \hat{\theta}_{\infty} := \lim_{\gamma \rightarrow \infty} \hat{\theta}_\gamma$.
We in fact have a deterministic bound of $|\hat{\theta}_{\gamma} - Y_{\text{mid}}|$ in the following lemma:

\begin{lemma}
\label{lem:uniformcontrol}
Let $Y_1, \ldots, Y_n$ be $n$ arbitrary points on $\mathbb{R}$, then
\[
|\hat{\theta}_{\gamma} - Y_{\text{mid}}| \leq 2(Y_{(n)} - Y_{(1)}) \frac{\log n}{\gamma}.
\]
\end{lemma}

We prove Lemma~\ref{lem:uniformcontrol} in Section~\ref{sec:method_appendix} of the appendix. It is important to note that, by Lemma~\ref{lem:uniformcontrol}, we need to consider $\gamma$ as large as $n$ to approximate $Y_{\text{mid}}$ with error that is of order $\frac{\log n}{n}$. Therefore, in settings where $Y_{\text{mid}}$ is optimal, we need $\gamma$ to be able to diverge with $n$. 

Estimators of form $\hat{\theta}_{\gamma}$ is simple, easy to compute via Newton's method (see Section~\ref{sec:optimization_appendix} of the appendix), and interpretable even for asymmetric distributions. The key question is of course, how do we select the power $\gamma$? It is necessary to allow $\gamma$ to increase with $n$ to attain adaptive rate but selecting a power $\gamma$ that is too large can introduce tremendous excess variance. As is often said, \emph{"with great power comes great responsibility"}. 

Before describing our approach in the next subsection, we give some remarks on two approaches that seem reasonable but in fact have significant limitations. 

\begin{remark}
\label{rem:cv_fail}
(Suboptimality of Cross-validation) \\
Cross-validation is a natural method for choosing the best estimator among some family, but this fails in our problem. To illustrate why, we consider the simpler problem where we choose between only the sample mean $\hat{\theta}_2$ and the sample midrange $\hat{\theta}_{\infty}$. We consider held-out validation where we divide our data into training data $D^{\text{train}}$ and test data $D^{\text{test}}$ each with $n$ data points. We compute $\hat{\theta}_2^{\text{train}}, \hat{\theta}_\infty^{\text{train}}$ on training data, evaluate test data MSE
\begin{align}
\hat{R}(\hat{\theta}_2^\text{train}) := \frac{1}{n} \sum_{i=1}^n (Y^{\text{test}}_i - \hat{\theta}_2^{\text{train}} )^2 = \frac{1}{n} \sum_{i=1}^n (Y^{\text{test}}_i - \bar{Y}^{\text{test}} )^2 + (\bar{Y}^{\text{test}} -  \hat{\theta}_2^{\text{train}} )^2, \label{eq:cv_mse}
\end{align}
for $\hat{\theta}^{\text{train}}_2$ and also $\hat{R}(\hat{\theta}_\infty^{\text{train}})$ for the $\hat{\theta}_\infty^{\text{train}}$ midrange estimator. Since the first term on the right hand side of~\eqref{eq:cv_mse} is constant, we select $\gamma = 2$ if $(\bar{Y}^{\text{test}} -  \hat{\theta}_2^{\text{train}} )^2 < (\bar{Y}^{\text{test}} -  \hat{\theta}_\infty^{\text{train}} )^2$. 

Now assume that the data follows the uniform distribution on $[\theta_0 - 1, \theta_0 + 1]$, so that the optimal estimator is the sample midrange $\hat{\theta}_{\infty}$. We observe that $\sqrt{n} (\hat{\theta}_2^{\text{train}} - \theta_0) \stackrel{d}{\rightarrow} N(0, 1/3)$ and $\sqrt{n} (\bar{Y}^{\text{test}} - \theta_0) \stackrel{d}{\rightarrow} N(0, 1/3)$ whereas $\sqrt{n} (\hat{\theta}_{\infty}^{\text{train}} - \theta_0) \rightarrow 0$ in probability. Hence, by the Portmanteau Theorem,
\begin{align*}
\liminf_{n\rightarrow \infty} \mathbb{P}(\text{selecting $\hat{\theta}_2$}) &= \liminf_{n\rightarrow \infty} \mathbb{P}( |\bar{Y}^{\text{test}} -  \hat{\theta}_2^{\text{train}} | < |\bar{Y}^{\text{test}} -  \hat{\theta}_\infty^{\text{train}} | ) \\
&= \liminf_{n\rightarrow \infty} \mathbb{P}( |\sqrt{n}(\bar{Y}^{\text{test}} - \theta_0) -  \sqrt{n}( \hat{\theta}_2^{\text{train}} - \theta_0) | \\
&\hspace{.7in} < | \sqrt{n}(\bar{Y}^{\text{test}} - \theta_0) -  \sqrt{n}( \hat{\theta}_\infty^{\text{train}} - \theta_0) | ) \\
&\geq \mathbb{P}( |W_1 - W_2| < |W_2|) > 0,
\end{align*}
where $W_1$ and $W_2$ are independent $N(0, 1/3)$ random variables. In other words, held-out validation has a non-vanishing probability of incorrectly selecting $\hat{\theta}_2$ over $\hat{\theta}_{\infty}$ even as $n \rightarrow \infty$ and thus has an error of order $1/\sqrt{n}$, which is far larger than the optimal $1/n$ rate. It is straightforward to extend the argument to the setting of $K$-fold cross-validation for any fixed $K$.
\end{remark}

\begin{remark} 
\label{rem:gg}
(Suboptimality of MLE with respect to the generalized Gaussian family) \\
We observe that $\hat{\theta}_\gamma$ is the maximum likelihood estimator for the center when the data follow the Generalized Normal GN$(\theta, \sigma, \gamma)$ distribution, which is also known as the Subbotin distribution~\citep{subbotin1923law}, whose density is of the form 
\begin{align*}
p(x \,;\, \theta, \sigma, \gamma) = \frac{1}{2 \sigma \Gamma(1+1/\gamma)} \exp \biggl( - \biggl| \frac{x- \theta}{\sigma} \biggr|^\gamma \biggr),
\end{align*}
where $\Gamma(t) := \int_0^\infty x^{t-1} e^{-x} dx$ denotes the Gamma function. This suggests a potential approach where we determine $\gamma$ by fitting the data to the potentially misspecified Generalized Gaussian family via likelihood maximization:
\begin{align*}
\argmin_{\gamma}\min_{\theta, \sigma} \frac{1}{n} \sum_{i=1}^n \biggl| \frac{Y_i - \theta}{\sigma} \biggr|^\gamma + \log \sigma + \log ( 2 \Gamma(1 + 1/\gamma)).
\end{align*}
This approach works well if the underlying density $p$ of the noise $Z_i$ belongs in the Generalized Gaussian family. Otherwise, it may be suboptimal: it may select a $\gamma$ that is too small when the optimal $\gamma$ is large and it may select a $\gamma$ that is too large when the optimal $\gamma$ is small. We give a precise and detailed discussion of the drawbacks of the generalized Gaussian MLE in Section~\ref{sec:compare_mle}.
\end{remark}

\subsection{Asymptotic variance}

Under the assumption that the noise $Z_i$ has a distribution symmetric around $0$, 
it is easy to see by symmetry that $\mathbb{E} \hat{\theta}_\gamma = \theta_0$ for any fixed $\gamma > 0$. We thus propose a selection scheme based on minimizing the variance. The finite sample variance of $\hat{\theta}_{\gamma}$ is intractable to compute, but for any fixed $\gamma > 1$, assuming $\mathbb{E}|Y-\theta_0|^{2 ({\gamma-1})}<\infty$, we have that  $\sqrt{n}(\hat{\theta}_\gamma-\theta_0)\stackrel{d}{\to}N(0,V(\gamma))$ as $n \rightarrow \infty$, where
\begin{align}
    V(\gamma) := \frac{\mathbb{E}|Y-\theta_0|^{2(\gamma-1)}}{\left[(\gamma-1)\mathbb{E}|Y-\theta_0|^{\gamma-2}\right]^2} \label{eq:v_defn}
\end{align}
is the asymptotic variance of $\hat{\theta}_{\gamma}$. Thus, from an asymptotic perspective, $\hat{\theta}_{\gamma}$ is a better estimator of $\theta_0$ if $V(\gamma)$ is small. When $\gamma$ is allowed to depend on $n$, $V(\gamma)$ may not be a good approximation of the finite sample variance of $\hat{\theta}_{\gamma}$, but the next example suggests that $V(\cdot)$ is still a sensible selection criterion. 

\begin{example}
When $Y_1, \ldots, Y_n \stackrel{d}{\sim} \text{Uniform}[\theta_0-1, \theta_0 + 1]$, straightforward calculation yields that $\mathbb{E}|Y - \theta_0|^q = \frac{1}{q+1}$ for any $q \in \mathbb{N}$ and thus, we have $V(\gamma) = \frac{1}{2\gamma - 1}$. We see that $V(\gamma)$ is minimized when $\gamma \rightarrow \infty$, in accordance with the fact that the sample midrange $Y_{\text{mid}}$ is the optimal estimator among the class of estimators  $\{\hat{\theta}_{\gamma} \}_{\gamma \geq 2}$. More generally, if $Y_i$ has a density $p(\cdot)$ supported on $[\theta_0 - 1, \theta_0 + 1]$ which is symmetric around $\theta_0$ and satisfies the property that $p(x)$ is bounded away from 0 and $\infty$ for all $x \in [\theta_0 - 1, \theta_0 + 1]$, then one may show that $V(\gamma) \propto \frac{1}{\gamma}$. On the other hand, if $Y_i \sim N(\theta_0, 1)$, then, using the fact that $\mathbb{E}|Y - \theta_0|^{\gamma} \asymp \gamma^{\gamma/2} e^{-\gamma/2}$, we can directly calculate that  that $V(\gamma) \asymp \frac{2^{\gamma}}{\gamma}$, which goes to infinity as $\gamma \rightarrow \infty$ as expected. Using the fact that $\hat{\theta}_2$ is the MLE, we have that $V(\gamma)$ is minimized at $\gamma = 2$ in the Gaussian case.
\end{example}

\subsection{Proposed procedure}
\label{sec:proposed_method}

We thus propose to select $\gamma$ by minimizing an estimate of the asymptotic variance $V(\gamma)$. For simplicity, we restrict our attention to $\gamma \geq 2$ in the main paper and discuss how to select $\gamma \in [1,2)$ in Remark~\ref{rem:robustness}. A natural estimator of $V(\gamma)$ is
\begin{align}
\hat{V}(\gamma) := \frac{\min_{\theta}\frac{1}{n}\sum_{i=1}^n |Y_i-\theta|^{2(\gamma-1)}}{\left[(\gamma-1)\min_{\theta }\frac{1}{n}\sum_{i=1}^n |Y_i-\theta|^{\gamma-2}\right]^2}. \label{eq:vhat_defn}
\end{align}

\begin{figure}
\begin{subfigure}{0.24\textwidth}
\includegraphics[scale=.23]{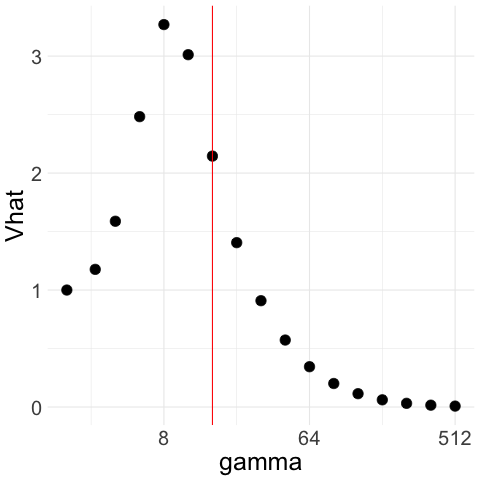}
\caption{$\hat{V}(\gamma)$ for Gaussian}
\label{fig:gauss_vhat}
\end{subfigure}
\begin{subfigure}{0.24\textwidth}
\includegraphics[scale=.23]{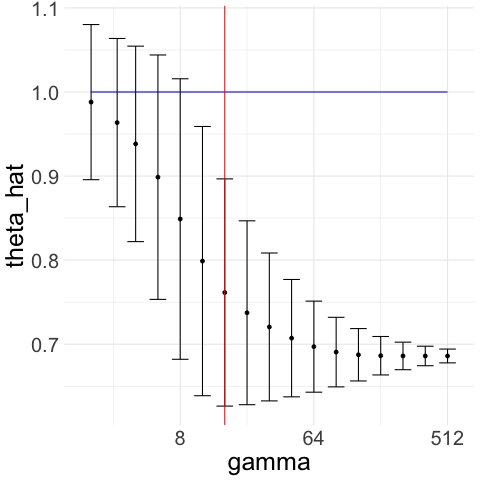}
\caption{$\hat{\theta}_{\gamma} \pm 2\sqrt{ \frac{\hat{V}(\gamma)}{n}}$}
\label{fig:gauss_lepski}
\end{subfigure}
\hspace{0.1in}
\begin{subfigure}{0.24\textwidth}
\includegraphics[scale=.23]{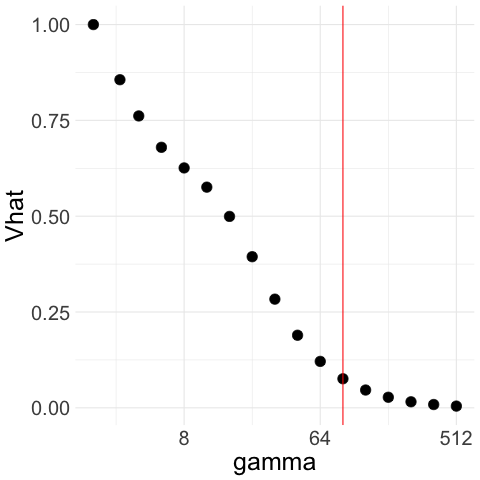}
\caption{$\hat{V}(\gamma)$ for trunc Gaussian}
\label{fig:tgauss_vhat}
\end{subfigure}
\begin{subfigure}{0.24\textwidth}
\includegraphics[scale=.23]{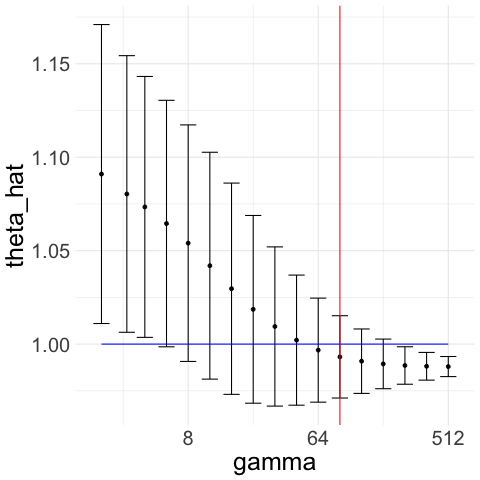}
\caption{$\hat{\theta}_{\gamma} \pm 2\sqrt{ \frac{\hat{V}(\gamma)}{n}}$}
\label{fig:tgauss_lepski}
\end{subfigure}
\caption{Red vertical line gives $\gamma_{\text{max}}$; blue horizontal line is the true $\theta_0$. We use $n=500$. We select $\gamma_{\max}$ as the largest $\gamma$ such that all confidence intervals to the left have a nonempty intersection.}
\label{fig:vhat_plots1}
\end{figure}

Although $\hat{V}(\gamma)$ has pointwise consistency in that it is a consistent estimator of $V(\gamma)$ for any fixed $\gamma$ (see Lemma~\ref{lem:as_min_converge} in Section~\ref{sec:as_min_converge} of the appendix), we require uniform consistency since our goal is to minimize $\hat{V}(\gamma)$ as a surrogate of $V(\gamma)$. This unfortunately does not hold; if we allow $\gamma$ to diverge with $n$, the error $|\hat{V}(\gamma) - V(\gamma)|$ can be arbitrarily large. This occurs because, if we fix $n$ and increase $\gamma$, the finite average $\frac{1}{n}\sum_{i=1}^n |Y_i - \theta|^{\gamma}$ does not approximate the population mean and behaves closer to $\frac{1}{n} \max_i |Y_i - \theta|^{\gamma}$ instead. Indeed, for any fixed $n$ and any deterministic set of points $Y_1, \ldots, Y_n$, we have
\begin{align}\label{eq:Vhat right tail}
\hat{V}(\infty) := \lim_{\gamma \rightarrow \infty} \hat{V}(\gamma) 
&= \lim_{\gamma \rightarrow \infty} \frac{n}{(\gamma-1)^2} \frac{ |Y_{(n)} - Y_{\text{mid}} |^{2(\gamma-1)} + |Y_{(1)} - Y_{\text{mid}}|^{2(\gamma-1)} }
{
\{ |Y_{(n)} - Y_{\text{mid}}|^{\gamma-2} + |Y_{(1)} - Y_{\text{mid}}|^{\gamma-2} \}^2 } \nonumber\\
&= \lim_{\gamma \rightarrow \infty} \frac{n}{2 (\gamma-1)^2} \biggl| \frac{Y_{(n)} - Y_{(1)}}{2} \biggr|^2 = 0.
\end{align}
Therefore, unconstrained minimization of $\hat{V}(\gamma)$ over all $\gamma \geq 1$ would select $\gamma = \infty$. See for example Figure~\ref{fig:gauss_vhat}, where we generate Gaussian noise $Z_i \sim N(0, 1)$ and plot $\hat{V}(\gamma)$ for a range of $\gamma$'s; although the population $V(\gamma)$ tends to infinity when  $\gamma$ is large, the empirical $\hat{V}(\gamma)$ increases for moderately large $\gamma$ but then, as $\gamma$ further increases, $\hat{V}(\gamma)$ decreases and tends to $0$.

Luckily, we can overcome this issue by restricting our attention to $\gamma$'s that are not too large. To be precise, we add an upper bound $\gamma_{\max} \geq 2$ and minimize $\hat{V}(\gamma)$ only among $\gamma \in [2, \gamma_{\max}]$. We select $\gamma_{\max}$ using Lepski's method, which is typically used to select smoothing parameters in nonparametric estimation problems \citep{lepskii1990problem, lepskii1991asymptotically} but can be readily adapted to our setting. The idea is to construct confidence intervals $\hat{\theta}_{\gamma} \pm \tau \sqrt{ \hat{V}(\gamma) / n}$ for a set of $\gamma$'s, starting with $\gamma = 2$, and take $\gamma_{\max}$ to be the largest $\gamma$ such that the confidence intervals all intersect. We would thus exclude $\gamma$ for which $\hat{\theta}_{\gamma}$ is far from $\theta_0$ and $\hat{V}(\gamma)$ is too small. 

This leads to our full estimation procedure below, which we refer to as \textbf{CAVS} (Constrained Asymptotic Variance Selector): 

\begin{algorithm}
\caption{Constrained Asymptotic Variance Selection (CAVS) algorithm}
\label{alg:cavs}
\noindent Let $\tau > 0$ be a tuning parameter and let $\mathcal{N}_n \subseteq [2, \infty]$ be the set of candidate $\gamma$'s. Define $\hat{V}(\gamma)$ as~\eqref{eq:vhat_defn} for $\gamma \in [2, \infty)$ and define $\hat{V}(\infty) := 0$.
\begin{enumerate}
\item Define $\gamma_{\max}$ as the largest $\gamma \in \mathcal{N}_n$ such that 
\begin{align*}
    \bigcap_{ \gamma \in \mathcal{N}_n,\, \gamma \leq \gamma_{\max}} \biggl[\hat{\theta}_{\gamma} - \tau \sqrt{\frac{\hat{V}(\gamma)}{n} },\,  \hat{\theta}_{\gamma} + \tau \sqrt{ \frac{\hat{V}(\gamma)}{n} } \biggr] \neq \emptyset.
\end{align*}

\item Select 
\[
\hat{\gamma} := \argmin_{\gamma \in \mathcal{N}_n,\, \gamma \leq \gamma_{\max}} \hat{V}(\gamma).
\]
\item Output 
\begin{align}
\hat{\theta} \equiv \hat{\theta}_{\hat{\gamma}}=\argmin_{\theta \in \mathbb{R}} \sum_{i=1}^n|Y_i-\theta|^{\hat{\gamma}}.
\label{eq:estimator}
\end{align}
\end{enumerate}

\end{algorithm}

The candidate set $\mathcal{N}_n$ can be the entire half-line $[2,\infty]$. In practice, we take $\mathcal{N}_n$ to be a finite set so that we are able to compute the minimizer of $\hat{V}(\gamma)$. A convenient and computationally efficient choice is $\mathcal{N}_n = \{2, 4, 8, \ldots, n, \infty\}$. 

We illustrate how the CAVS procedure works with two examples in Figure~\ref{fig:vhat_plots1}. In Figure~\ref{fig:gauss_vhat}, we generate Gaussian noise $Z_i \sim N(0, 1)$; we plot $\hat{V}(\gamma)$ for a exponentially increasing sequence of $\gamma$'s ranging from 2 to 512. The constraint upper bound $\gamma_{\max}$ is given by the red line in the figure. Unconstrained minimization of $\hat{V}(\gamma)$ leads to $\hat{\gamma} = 512$.
Figure~\ref{fig:gauss_lepski} illustrates the Lepski method that we use to choose upper bound $\gamma_{\max}$: we compute confidence intervals of width $\tau \sqrt{ \frac{\hat{V}(\gamma)}{n} }$ around $\hat{\theta}_{\gamma}$ for the whole range of $\gamma$'s. To get $\gamma_{\max}$, we pick the largest $\gamma$ such that the intersection of all the confidence intervals to the left of $\gamma_{\max}$ is non-empty. 

This allows us to avoid the region where $\hat{V}(\gamma)$ is very small but the actual asymptotic variance $V(\gamma)$ is very large. Indeed, if $V(\gamma)$ is much larger than the variance $\text{Var}(Z)$, then $\hat{\theta}_{\gamma}$ likely to be quite far from the sample mean $\hat{\theta}_2$ and thus, if $\hat{V}(\gamma)$ is also small, then $\hat{\theta}_{\gamma} \pm \tau \sqrt{ \frac{\hat{V}(\gamma)}{n}}$ is unlikely to overlap with the confidence interval around the sample mean. 

Therefore, with Gaussian noise, CAVS selects $\hat{\gamma} = 2$ by minimizing $\hat{V}(\gamma)$ only to the left of $\gamma_{\max}$ (red line) in Figure~\ref{fig:gauss_vhat}. In contrast, if $V(\gamma)$ decreases as $\gamma$ increases, then $\hat{\theta}_{\gamma}$ remains close to the sample mean $\hat{\theta}_2$ and the confidence interval $\hat{\theta}_{\gamma} \pm \tau \sqrt{ \frac{\hat{V}(\gamma)}{n}}$ overlaps with that of the sample mean even when $\gamma$ is large, which means we would select a large $\gamma_{\max}$ as desired. We illustrate this in Figure~\ref{fig:tgauss_vhat} and~\ref{fig:tgauss_lepski}, where we generate truncated Gaussian noise $Z$ by truncating at $|Z| \leq 2$; that is, we generate Gaussian samples and keep only those that lie in the interval $[-2, 2]$. In this case, the optimal $\gamma$ is $\gamma = \infty$ and the optimal rate is $1/n$. From Figure~\ref{fig:tgauss_vhat}, we see that our procedure picks a large $\hat{\gamma} = 128$. 

\begin{remark} (Selecting $\tau$ parameter) 

Our proposed CAVS procedure has a tuning parameter $\tau$ which governs the strictness of the $\gamma_{\max}$ constraint. Smaller $\tau$ will in general result in a smaller $\gamma_{\max}$ and hence a stronger constraint. For our theoretical results, namely Theorem~\ref{thm:adaptive_rate}, it suffices to choose $\tau$ to be very slowly growing so that $\frac{\tau}{\sqrt{\log \log n}} \rightarrow \infty$. For practical data analysis applications, we recommend $\tau = 1$ as a conservative choice based on simulation studies in Section~\ref{sec:simulation}
\end{remark}

\begin{figure}[htp]
\centering
\begin{subfigure}{0.4\textwidth}
\centering
\includegraphics[scale=.25]{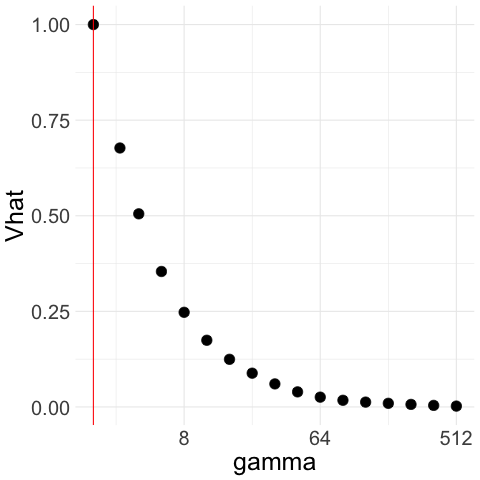}
\caption{$\hat{V}$ for asymmetric density}
\label{fig:asymm1a}
\end{subfigure}
\hspace{.05in}
\begin{subfigure}{0.4\textwidth}
\centering
\includegraphics[scale=.25]{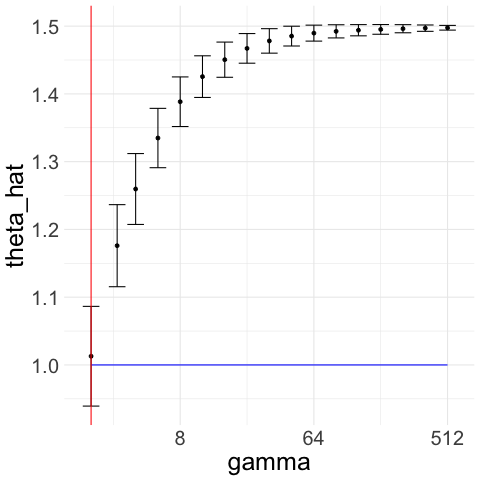}
\caption{$\hat{\theta}_{\gamma} \pm 2 \sqrt{ \frac{\hat{V}(\gamma)}{n} }$}
\label{fig:asymm1b}
\end{subfigure}
\caption{We generate mean zero $Z_i$ from an asymmetric mixture distribution $\frac{2}{3} \text{Unif}[-1, 0] + \frac{1}{3} \text{Unif}[0, 2]$. Note that $\hat{\theta}_3 \pm 2 \sqrt{ \frac{\hat{V}(3)}{n} }$ does not overlap with $\hat{\theta}_2 \pm 2 \sqrt{ \frac{\hat{V}(2)}{n} }$ because $\mathbb{E} \hat{\theta}_2 \neq \mathbb{E} \hat{\theta}_3$ due to the asymmetry. Red vertical line gives $\gamma_{\text{max}}$; blue horizontal line is the true $\theta_0$.}
\label{fig:asymm}
\end{figure}

\begin{remark} (Robustness to asymmetry) 

One important aspect of CAVS is that it is robust to violations of the symmetry assumption. If the density $p$ of the noise $Z_i$ has mean zero but is asymmetric (so that $\theta_0$ is the mean of $Y_i$), then, for various $\gamma$'s that are greater than 2, the $\gamma$-th center of $Y_i = \theta_0 + Z_i$ may be different from $\theta_0$; that is $\theta^*_\gamma := \argmin_{\theta} \mathbb{E}|Y - \theta|^\gamma \neq \theta_0$ so that $\hat{\theta}_{\gamma}$ is a biased estimator of $\theta_0$. In such cases however, the confidence interval $\hat{\theta}_{\gamma} \pm \tau \sqrt{ \frac{\hat{V}(\gamma)}{n} }$ will, for large enough $n$, be concentrated around $\mathbb{E} \hat{\theta}_{\gamma} = \theta^*_\gamma$ and thus not overlap with the confidence interval about the sample mean $\hat{\theta}_2 \pm \tau \sqrt{ \frac{\hat{V}(2)}{n} }$, which will concentrated around $\mathbb{E} \hat{\theta}_2 = \theta_0$. Therefore, we would have $\gamma_{\max} < \gamma$ and the constraint would thus exclude any biased $\hat{\theta}_{\gamma}$. We illustrate an example in Figure~\ref{fig:asymm} where because $\hat{\theta}_3$ is biased, we have that $\gamma_{\max} = 2$ and thus, we select $\hat{\gamma} = 2$ and the resulting estimator $\hat{\theta}_{\hat{\gamma}}$ still converges to $\theta_0$. Indeed, our basic convergence guarantee--formalized in Theorem~\ref{thm:basic}--does not require the noise distribution to be symmetric around $0$, it only requires the noise to have mean zero. 

\end{remark}

\begin{remark} 
\label{rem:regression} (Extension to the regression setting)

We can directly extend our estimation procedure to the linear regression setting. Suppose we observe $(Y_i, X_i)$ for $i = 1,2,\ldots, n$ where $X_i$ is a random vector on $\mathbb{R}^d$, $Y_i = X_i^\top \beta_0 + Z_i$, and $Z_i$ is an independent noise with a distribution symmetric around $0$. 

Then, we would compute, for each $\gamma$ in a set $\mathcal{N}_n \subset [2, \infty]$, 
\begin{align*}
\hat{\beta}_{\gamma} &= \argmin_{\beta \in \mathbb{R}^p} \sum_{i=1}^n |Y_i - X_i^\top \beta|^\gamma  \quad \text{ and } \\
\hat{V}(\gamma) &= \frac{\min_{\beta \in \mathbb{R}^d} \frac{1}{n} \sum_{i=1}^n |Y_i - X_i^\top \beta |^{2(\gamma-1)} }{(\gamma-1)^2 \{ \min_{\beta \in \mathbb{R}^d} \frac{1}{n} \sum_{i=1}^n |Y_i - X_i^\top \beta |^{\gamma-2} \}^2 }.
\end{align*}

We define $\hat{\Sigma}_X := \frac{1}{n} \sum_{i=1}^n X_i X_i^\top$. Using Taylor expansion, it is straightforward to show that $\sqrt{n} \hat{\Sigma}^{1/2}_X (\hat{\beta}_{\gamma} - \beta_0) \stackrel{d}{\rightarrow} N(0, V(\gamma) I_d)$. Thus, for a given $\tau > 0$, our estimation procedure first computes $\gamma_{\max}$ as the largest $\gamma \in \mathcal{N}_n$ such that
\begin{align*}
\bigcap_{\gamma \in \mathcal{N}_n, \, \gamma \leq \gamma_{\max}} \bigotimes_{j=1}^p \biggl[ \bigl(\hat{\Sigma}_X^{1/2} \hat{\beta}_{\gamma}\bigr)_j - \tau \sqrt{ \frac{\hat{V}(\gamma)}{n} }, \,
\bigl( \hat{\Sigma}_X^{1/2} \hat{\beta}_{\gamma} \bigr)_j + \tau \sqrt{ \frac{\hat{V}(\gamma)}{n} } \biggr] \neq \emptyset,
\end{align*}
where we use the $\otimes$ notation to denote the Cartesian product. Then, we select the minimizer $\hat{\gamma} = \argmin_{\gamma \in \mathcal{N}_n, \, \gamma \leq \gamma_{\max}} \hat{V}(\gamma)$ and output $\hat{\beta}_{\hat{\gamma}}$.

\end{remark}

\begin{remark} \label{rem:robustness} (Selecting $\gamma \in [1, 2)$)

When the noise $Z_i$ is heavy-tailed, it is desirable to allow consideration of $\gamma \in [1,2)$; note that $\gamma =1$ corresponds to the sample median $\hat{\theta}_1 = \argmin_{\theta} \sum_{i=1}^n |Y_i - \theta|$. For $\gamma \in [1,2)$, the estimator $\hat{V}(\gamma)$ given in~\eqref{eq:vhat_defn} is not appropriate. In particular, if $Z_i$ has a density $p$ and population median 0 and that $p(0) > 0$, then the asymptotic variance of sample median is $V(1) = \frac{1}{4 p(0)^2}$ instead of~\eqref{eq:v_defn}. For $\gamma \in (1,2)$, expression~\eqref{eq:v_defn} holds but the estimator $\hat{V}(\gamma)$ may behave poorly because of the negative power in the denominator. We do not have a general way of estimating $V(\gamma)$ for $\gamma < 2$. In the specific case of the sample median $(\gamma = 1)$, there are various good estimators of the variance. For instance, \cite{bloch1968simple} proposed an approach based on density estimation and \cite{lai1983adaptive} proposed an approach based on the bootstrap. The general idea of selecting an estimator using asymptotic variance is not specific to the $L_\gamma$-centers; one can also add say Huber loss minimizers into the set of candidate estimators provided that there is a good way to estimate the asymptotic variance.
\end{remark}

\subsection{Basic properties of the estimator}

Using the definition of $\hat{\gamma}$, we can directly show that $\hat{\theta}_{\hat{\gamma}}$ must be close to the sample mean $\bar{Y}$ and that the error of $\hat{\theta}_{\hat{\gamma}}$ is at most $O( \tau \sqrt{\sigma^2 / n} )$ where $\sigma^2 := \text{Var}(Z)$.
\begin{theorem}
\label{thm:basic}
Let $\hat{\sigma}^2$ be the empirical variance of $Y_1, \ldots, Y_n$. For any $n$, it holds surely that
\[
| \hat{\theta}_{\hat{\gamma}} - \bar{Y} | \leq 2 \tau \sqrt{ \frac{ \hat{\sigma}^2}{n} }.
\]
Therefore, if we additionally have that $\sigma^2 := \mathbb{E}|Z|^2 < \infty$, then, writing $\theta_0 = \mathbb{E} Y_1$, 
\[
\mathbb{E} | \hat{\theta}_{\hat{\gamma}} - \theta_0 | \lesssim \tau \sqrt{ \frac{\sigma^2}{n}}.
\]
\end{theorem}

\begin{proof}
Since $\hat{\gamma} \leq \gamma_{\max}$, we have by the definition of $\gamma_{\max}$ that
\begin{align*}
\hat{\theta}_{\hat{\gamma}} + \tau \sqrt{ \frac{\hat{V}(\hat{\gamma})}{n} } & \geq \hat{\theta}_2 - \tau \sqrt{ \frac{\hat{V}(2)}{n} } \, \text{ and } \\
 \hat{\theta}_{\hat{\gamma}} - \tau \sqrt{ \frac{\hat{V}(\hat{\gamma})}{n} } & \leq \hat{\theta}_2 - \tau \sqrt{ \frac{\hat{V}(2)}{n} }.
\end{align*}
Since $\hat{V}(\hat{\gamma}) \leq \hat{V}(2)$ by definition of $\hat{\gamma}$ and since $\hat{\theta}_2 = \bar{Y}$ and $\hat{V}(2) = \hat{\sigma}^2$, the first claim immediately follows. The second claim directly follows from the first claim. 
\end{proof}

It is important to note that Theorem~\ref{thm:basic} does not require symmetry of the noise distribution $P$. If $Y_i$ has a distribution \emph{asymmetric} around $\theta_0$ but $\mathbb{E} Y = \theta_0$, then Theorem~\ref{thm:basic} implies that $\hat{\theta}_{\hat{\gamma}}$ converges to $\theta_0$ as might be desired. 

\begin{remark}
An important property of $\hat{\gamma}$ is that it is shift and scale invariant in the following sense: if we scale our data with the transformation $\tilde{Y}_i = b Y_i + a$ where $b > 0$ and $a \in \mathbb{R}$ and then compute $\tilde{\gamma}$ on $\{ \tilde{Y}_1, \ldots, \tilde{Y}_n \}$, then $\tilde{\gamma} = \hat{\gamma}$. This follows from the fact that $\hat{V}(\gamma)/\hat{V}(2)$ is shift and scale invariant. Likewise, we see that $\hat{\theta}_{\hat{\gamma}}$ is shift and scale equivariant in that if we compute $\tilde{\theta}_{\tilde{\gamma}}$ on $\{ \tilde{Y}_1, \ldots, \tilde{Y}_n \}$, then $\tilde{\theta}_{\tilde{\gamma}} = b \hat{\theta}_{\hat{\gamma}} + a$. 
\end{remark}

\section{Adaptive rate of convergence}
\label{sec:theory}


Theorem~\ref{thm:basic} shows that, so long as $\tau$ is chosen to be not too large and the noise $Z_i$ has finite variance, then our proposed estimator has an error $\mathbb{E} | \hat{\theta}_{\hat{\gamma}} - \theta_0 |$ that is at most $\tilde{O}(n^{-1/2})$. In this section, we show that if the noise $Z_i$ has a density $p(\cdot)$ that is in a class of compactly supported densities, then our estimator can attain an adaptive rate of convergence of $\tilde{O}(n^{-\frac{1}{\alpha}})$ for any $\alpha \in (0, 2]$, depending on a moment property of the noise distribution.

\begin{theorem}
\label{thm:adaptive_rate}
Suppose $Z_1, Z_2, \ldots, Z_n$ are independent and identically distributed with a distribution $P$ symmetric around $0$. Suppose there exists $\alpha \in (0, 2]$, $a_1 \in (0, 1]$ and $a_2 \geq 1$ such that $\frac{a_1}{\gamma^\alpha} \leq \mathbb{E}|Z|^\gamma \leq \frac{a_2}{\gamma^\alpha}$ for all $\gamma \geq 1$. Let $\mathcal{N}_n$ be a subset of $[2, \infty]$ with $M_n := \sup \mathcal{N}_n$ and suppose $\mathcal{N}_n$ contains $2^k$ for all integer $k \leq n \wedge \log_2 M_n$.

Let $C_{a_1, a_2, \alpha} > 0$ be a constant that depends only on $a_1, a_2, \alpha$; let $\hat{\theta}_{\hat{\gamma}}$ be defined as~\eqref{eq:estimator}. The following then hold:
\begin{enumerate}
\item If $\frac{\tau}{\sqrt{\log \log n}} \rightarrow \infty$, then
\[
|\hat{\theta}_{\hat{\gamma}} - \theta_0| \leq O_p \biggl( C_{a_1, a_2, \alpha} \biggl\{ \left(\frac{\log^{\alpha +1} n}{ n}\right)^{\frac{1}{\alpha}} \vee \frac{\log n}{M_n} \biggr\} \biggr).
\]
\item If $\tau \geq \sqrt{\log n}$, then
\[
\mathbb{E} |\hat{\theta}_{\hat{\gamma}} - \theta_0| \leq C_{a_1, a_2, \alpha} \biggl\{ \left(\frac{\log^{\alpha +1} n}{ n}\right)^{\frac{1}{\alpha}} \vee \frac{\log n}{M_n} \biggr\}.
\]
\end{enumerate}
\end{theorem}

Therefore, we can choose $M_n \geq 2^n$ and $\tau = \sqrt{\log n}$, without any knowledge of $\alpha$, so that our estimator has an adaptive rate of convergence 
\[
\mathbb{E} |\hat{\theta}_{\hat{\gamma}} - \theta_0| \lesssim_{a_1, a_2, \alpha} \biggl( \frac{\log^{\alpha+1} n}{n} \biggr)^{\frac{1}{\alpha}}
\]
where $\alpha$ can take on any value in $(0, 2]$ depending on the underlying noise distribution. The adaptive rate $( \frac{\log^{\alpha+1}}{n} )^{1/\alpha}$ is, up to log-factors, minimax optimal for the class of densities satisfying $\mathbb{E}|Z|^{\gamma} \propto \gamma^{-\alpha}$; see Remark~\ref{rem:minimax_opt} for more details. 

We relegate the proof of Theorem~\ref{thm:adaptive_rate} to Section~\ref{sec:adaptive_rate_proof} of the appendix, but give a sketch of the proof ideas here. First, by using the moment condition $\frac{a_1}{\gamma^{\alpha}}\leq \mathbb{E}|Z|^{\gamma} \leq \frac{a_2}{\gamma^{\alpha}}$ as well as Talagrand's inequality, we give the following uniform bound to $\hat{V}(\gamma)$: that $\hat{V}(\gamma) \asymp_{a_1, a_2, \alpha} \gamma^{\alpha-2}$ for all $2 \leq \gamma \leq \bigl( \frac{n}{\log n} \bigr)^{\frac{1}{\alpha}}$. Using this bound in conjunction with another uniform bound on $| \hat{\theta}_{\gamma} - \theta_0|$, we then can guarantee that $\gamma_{\max}$ is large enough in that $\gamma_{\max} \gtrsim_{a_1, a_2} \bigl( \frac{n}{\log n} \bigr)^{\frac{1}{\alpha}} \vee M_n$. These results in turn yields the key fact that $\hat{\gamma}$ is also sufficiently large in that $\hat{\gamma} \gtrsim_{a_1, a_2} \bigl( \frac{n}{\log n} \bigr)^{\frac{1}{\alpha}} \vee M_n $. We then bound the error of  $\hat{\theta}_{\hat{\gamma}}$ by the inequality
\[
|\hat{\theta}_{\hat{\gamma}} - \theta_0| \leq |\hat{\theta}_{\hat{\gamma}} - Y_{\text{mid}}| + | \theta_0 - Y_{\text{mid}} |,
\]
where $Y_{\text{mid}}$ is the sample midrange. We control the first term $|\hat{\theta}_{\hat{\gamma}} - Y_{\text{mid}}|$ through Lemma~\ref{lem:uniformcontrol} and the second term $| \theta_0 - Y_{\text{mid}} |$ using the moment condition. The resulting bound gives the desired conclusion of Theorem~\ref{thm:adaptive_rate}. 


\begin{remark}
The condition that $\mathbb{E}|Z|^\gamma \propto \gamma^{-\alpha}$ for all $\gamma \geq 1$ implies that $Z$ is supported on $[-1, 1]$. This is not as restrictive as it appears: using the fact that $\hat{\theta}_{\hat{\gamma}}$ is scale equivariant, it is straightforward to show that if $Z_i$ takes value on $[-b, b]$ for any $b > 0$ and satisfies $\mathbb{E}|Z|^\gamma \propto b^{\alpha} \gamma^{-\alpha}$, then we have that $\mathbb{E} | \hat{\theta}_{\hat{\gamma}} - \theta_0 | \leq b \cdot  C_{a_1, a_2, \alpha} \bigl\{ \bigl( \frac{\log^{\alpha+1} n}{n} \bigr)^{\frac{1}{\alpha}} \vee \frac{\log n}{M_n} \bigr\}$.
\end{remark}


\subsection{On the moment condition in Theorem~\ref{thm:adaptive_rate}}

The moment condition $\frac{a_1}{\gamma^{\alpha}} \leq \mathbb{E}|Z|^\gamma \leq \frac{a_2}{\gamma^{\alpha}}$ constrains the behavior of the density $p(\cdot)$ around the boundary of the support $[-1, 1]$. The following Proposition formalizes this intuition.

\begin{proposition}
\label{prop:moment_example}
Let $\alpha \in (0, 2)$ and suppose $X$ is a random variable with density $p(\cdot)$ satisfying
\[
C_{\alpha,1} (1 - |x|)_+^{\alpha-1} \leq p(x) \leq C_{\alpha,2} (1 - |x|)_+^{\alpha-1},\quad \forall x \in [-1, 1],
\]
for $C_{\alpha,1}, C_{\alpha,2} > 0$ dependent only on $\alpha$. Then, there exists $C'_{\alpha,1}, C'_{\alpha,2} > 0$, dependent only on $\alpha$, such that, for all $\gamma \geq 1$,
\[
\frac{C'_{\alpha,1}}{\gamma^\alpha} \leq \mathbb{E}|X|^\gamma \leq \frac{C'_{\alpha,2}}{\gamma^\alpha}.
\]
\end{proposition}
We prove Proposition~\ref{prop:moment_example} in Section~\ref{sec:proof_of_examples} of the Appendix. 

\begin{figure}[htp]
\centering
\begin{subfigure}{0.45\textwidth}
\centering
\includegraphics[scale=.2]{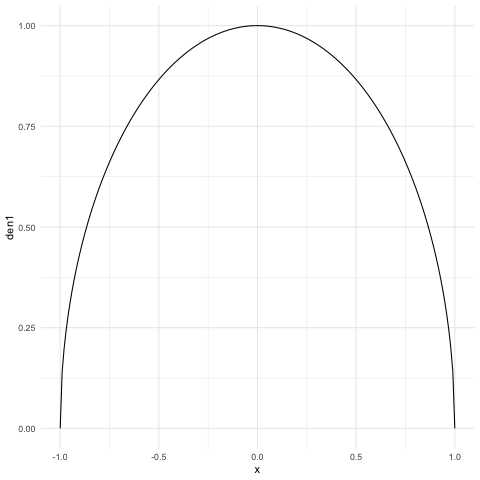}
\caption{Semicircle density}
\label{fig:semicircle}
\end{subfigure}
\begin{subfigure}{0.45\textwidth}
\centering
\includegraphics[scale=.2]{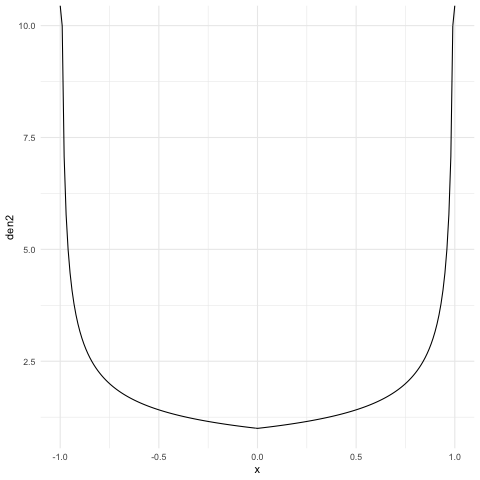}
\caption{U-shaped density}
\label{fig:ushape}
\end{subfigure}
\caption{}
\label{fig:densities}
\end{figure}

\begin{example} 
\label{ex:moment_densities}
Using Proposition~\ref{prop:moment_example}, we immediately obtain examples of noise distributions where the rates of convergence of our location estimator $\hat{\theta}_{\hat{\gamma}}$ vary over a wide range. 
\begin{enumerate}
\item When $Z$ has the semicircle density $p(x) \propto (1 - |x|^2)^{1/2}$ (see Figure~\ref{fig:semicircle}), then $\mathbb{E}|Z|^{\gamma} \propto \gamma^{-\frac{3}{2}}$ so that $\hat{\theta}_{\hat{\gamma}}$ has rate $\tilde{O}(n^{- \frac{2}{3}})$, where we use the $\tilde{O}(\cdot)$ notation to indicate that we have ignored polylog terms.
\item When $Z \sim\text{Unif}[-1,1]$, we have that $\mathbb{E}|Z|^\gamma = \frac{1}{\gamma+1}$ so that $\hat{\theta}_{\hat{\gamma}}$ has rate $\tilde{O}(n^{-1})$.
\item More generally, let $q$ be a symmetric continuous density on $\mathbb{R}$ and let $p$ be a density that results from truncating $q$, that is, $p(x) \propto q(x) \mathbbm{1}\{ |x| \leq 1 \}$. If $p(1) = p(-1) > 0$, then $\frac{a_1}{\gamma} \leq \mathbb{E}|Z|^\gamma \leq \frac{a_2}{\gamma}$ where $a_1, a_2$ depend on $q$. In particular, if $Z$ is a truncated Gaussian, then $\hat{\gamma}_{\hat{\gamma}}$ also has $\tilde{O}(n^{-1})$ rate. 
\item Suppose $Z$ has a U-shaped density of the form $p(x) \propto (1 - |x|)^{-\frac{1}{2}}$ (Figure~\ref{fig:ushape}), then $\mathbb{E}|Z|^{\gamma} \propto \gamma^{-\frac{1}{2}}$ so that $\hat{\theta}_{\hat{\gamma}}$ has rate $\tilde{O}(n^{-2})$. 
\end{enumerate}
\end{example}

\begin{remark}
\label{rem:minimax_opt}
By Proposition~\ref{prop:hellinger_bound} and the subsequent Remark~\ref{rem:hellinger_bound} in Section~\ref{sec:proof_of_examples} of the Appendix, we have that if a density $p$ is of the form $p(x) = C_\alpha (1 - |x|)^{\alpha - 1} \mathbbm{1}\{ |x| \leq 1\}$ for $\alpha \in (0, 2)$, then we have that, writing $H^2(\theta_1, \theta_2) := \int \bigl(\sqrt{ p(x - \theta_1) } - \sqrt{p(x - \theta_2)} \bigr)^2 dx$, 
\[
C_{\alpha, 1} |\theta_1 - \theta_2|^{\alpha} \leq H^2(\theta_1, \theta_2) \leq C_{\alpha, 2} |\theta_1 - \theta_2|^{\alpha},
\]
for $C_{\alpha, 1}$ and $C_{\alpha, 2}$ dependent only on $\alpha$. 
From \citet[][Proposition 1]{lecam1973}, any estimator $\hat{\theta}$ has a rate lower bounded by the fact that $H^2(\hat{\theta}, \theta_0) \gtrsim \frac{1}{n}$ so that among the class of densities 
\begin{align}
\mathcal{P}_{a_1, a_2} := \biggl\{ p \,:\, \text{ symmetric},\, \frac{a_1}{\gamma^\alpha} \leq \int |x|^\gamma p(x) dx \leq \frac{a_2}{\gamma^\alpha}, \text{$\forall \gamma \geq 1$, for some $\alpha \in (0, 2]$} \biggr\},
\label{eq:moment_density_class}
\end{align}
our proposed estimator $\hat{\theta}_{\hat{\gamma}}$ has a rate of convergence that is minimax optimal up to poly-log factors. 
\end{remark}

\subsection{Comparison with the MLE}
\label{sec:compare_mle}

Recall from Remark~\ref{rem:gg} that for $\theta \in \mathbb{R}$ and $\sigma, \gamma > 0$, the generalized Gaussian distribution (also known as the Subbotin distribution) has a density of the form $p(x : \theta, \sigma, \gamma) = \frac{1}{2 \sigma \Gamma(1 + \gamma^{-1})} \exp \bigl( - \bigl| \frac{x - \theta}{\sigma} \bigr|^\gamma \bigr)$. We note that the uniform distribution on $[-\sigma, \sigma]$ is a limit point of the generalized Gaussian class where we let $\gamma \rightarrow \infty$. 

Using univariate observations $Y_1, \ldots, Y_n$, we may then compute the MLE of $\gamma$ with respect to the generalized Gaussian family:
\begin{align*}
\hat{\gamma}_{\text{MLE}} = \argmin_{\gamma} \min_{\theta, \sigma} \frac{1}{n} \sum_{i=1}^n \biggl| \frac{Y_i - \theta}{\sigma} \biggr|^\gamma + \log \sigma + \log  \Gamma\biggl(1 + \frac{1}{\gamma}\biggr).
\end{align*}
For any fixed $\gamma$, we may minimize over $\theta$ and $\sigma$ to obtain that 
\[
\hat{\gamma}_{\text{MLE}} = \argmin_{\gamma} L_n(\gamma),
\]where
\[
L_n(\gamma):=\frac{1}{\gamma} \log \biggl( \min_{\theta} \frac{1}{n} \sum_{i=1}^n |Y_i - \theta|^\gamma \biggr) + \frac{1+ \log \gamma}{\gamma} + \log \Gamma\biggl(1 + \frac{1}{\gamma}\biggr).
\]
A natural question then is how good is $\hat{\gamma}_{\text{MLE}}$ as a selection procedure? Would the resulting location estimator $\hat{\theta}_{\hat{\gamma}_{\text{MLE}}}$ have good properties? 
If the density of $Y_i$ belongs to the generalized Gaussian class, then we expect $\hat{\gamma}_{\text{MLE}}$ to perform well. But when there is model misspecification, we show in this section that $\hat{\gamma}_{\text{MLE}}$ performs suboptimally compared to the CAVS estimator that we propose in Section~\ref{sec:proposed_method}.

To start, let us define the population level likelihood function for every $\gamma > 0$
\begin{align*}
L(\gamma) &:= \min_{\theta, \sigma} \mathbb{E} \biggl| \frac{Y - \theta}{\sigma} \biggr|^\gamma + \log (2\sigma) + \log \Gamma\bigl(1 + \frac{1}{\gamma}\bigr) \\
&= \frac{1}{\gamma} \log \bigl( \min_{\theta} \mathbb{E}|Y - \theta|^\gamma \bigr) + \frac{1+\log \gamma}{\gamma} + \log \Gamma\biggl( 1 + \frac{1}{\gamma} \biggr).
\end{align*}
We define $L(\infty) := \lim_{\gamma \rightarrow \infty} L(\gamma)$ and $L_n(\infty) := \lim_{\gamma \rightarrow \infty} L_n(\gamma)= \log \bigl\{ (Y_{(n)}-Y_{(1)})/2 \bigr\}$. We note that if $\mathbb{E}|Y|^\gamma=\infty$, then $L(\gamma)=\infty$ and if $Y$ is supported on the real line, then $L(\infty)=\infty$. Moreover, by Lemma~\ref{lem:as_min_converge} (in Section~\ref{sec:as_min_converge} of the appendix), we have that, for any fixed $\gamma\in \mathbb{R}\cup\{\infty\}$, we have that $L_n(\gamma) \stackrel{a.s.}{\rightarrow} L(\gamma)$.

Define $\gamma_{\text{MLE}}^* = \argmin_{\gamma \geq 2} L(\gamma)$ as the minimizer of $L(\gamma)$. We show in the next Proposition that when the noise $Z_i$ is supported on $[-1,1]$ with a small but positive density value at the boundary, then $\gamma_{\text{MLE}}^* < \infty$ even though the optimal selection of $\gamma$ is to take $\gamma \rightarrow \infty$ since the sample midrange $\hat{\theta}_{\infty}$ would have a rate of convergence that is at least as fast as $\tilde{O}(n^{-1})$.

\begin{proposition}
\label{prop:pop_mle}
Suppose $Y = Z + \theta_0$ where $Z$ has a distribution symmetric around $0$. Define $\gamma^*_{\text{MLE}} = \argmin_{\gamma > 0} L(\gamma)$. 
\begin{enumerate}
\item If $Z$ is supported on all of $\mathbb{R}$, then $\gamma^*_{\text{MLE}} < \infty$.
\item Suppose $Z$ has a density $p$ supported and continuous on $[-1, 1]$. Let $\gamma_{\text{E}} \approx 0.57721$ be the Euler--Mascheroni constant. If the density value at the boundary satisfies $p(1) < \frac{1}{2} e^{\gamma_{\text{E}} - 1}$, then $\gamma^*_{\text{MLE}} < \infty$.
\item Suppose $Z$ has a density $p$ supported and continuous on $[-1, 1]$. If the density value at the boundary satisfies $p(1) > \frac{1}{2} e^{\gamma_{\text{E}} - 1}$, then $\gamma = \infty$ is a local minimum of $L(\gamma)$. 
\end{enumerate}
\end{proposition}
We relegate the proof of Proposition~\ref{prop:pop_mle} to Section~\ref{sec:pop_mle_proof} of the Appendix.

If the noise density $p$ is continuous and has boundary value $p(1) \in (0,\frac{1}{2} e^{\gamma_{\text{E}} - 1})$, then Proposition~\ref{prop:pop_mle} suggests that we would not expect  $\hat{\gamma}_{\text{MLE}} \rightarrow \infty$. More precisely, we have that $L(\gamma^*_{\text{MLE}}) < L(\infty)$ and thus, by Lemma~\ref{lem:as_min_converge}, when $n$ is large enough, we also have $L_n(\gamma^*_{\text{MLE}}) < L_n(\infty)$ almost surely. Therefore, selecting $\gamma$ by minimizing $L_n$ would always favor a finite $\gamma = \gamma^*_{\text{MLE}}$ over $\gamma = \infty$. As a result, selecting $\gamma$ based on MLE yields a suboptimal rate of $n^{-1/2}$.

In contrast, Theorem~\ref{thm:adaptive_rate} shows that under the same setting, our proposed CAVS estimator selects a divergent $\hat{\gamma}$ which can yield an error that is smaller than $\tilde{O}(n^{-1/2})$ for $\hat{\theta}_{\hat{\gamma}}$. In fact, there are settings in which the density at the boundary is equal to zero, that is, $p(1) = 0$, where our proposed estimator can $\hat{\theta}_{\hat{\gamma}}$ have a rate of convergence that is faster than $n^{-1/2}$; for example, we see in that $| \hat{\theta}_{\hat{\gamma}} - \theta_0 |$ is $\tilde{O}( n^{-2/3} )$ when the noise has the semicircle density.



We note that although Proposition~\ref{prop:pop_mle} is stated for $Z$ supported on $[-1,1]$, by scale invariance of $\gamma^*_{\text{MLE}}$, Proposition~\ref{prop:pop_mle} holds for support of the form $[-b, b]$, where the the condition on the density generalizes to $p(b) > \frac{1}{2b} e^{\gamma_0 - 1}$.

\begin{remark}
Another drawback, one that is perhaps more alarming, of selecting $\gamma$ based on the Generalized Gaussian likelihood is that the resulting location estimator may have a standard deviation (and hence error) that is larger than $O(n^{-1/2})$.

Consider the following example: let $p_1$ be the density of $|W|^{\frac{1}{3}}\text{sign}(W)$, where $W$ follows the standard Cauchy distribution, let $p_2(x) \propto \exp\bigl(-|x|^3\bigr)$, and let the noise $Z$ have a mixture density $p = \delta p_1+(1-\delta) p_2$ for some $\delta \in [0, 1]$. We let $Y = Z + \theta_0$ as usual. 

If $\delta = 0$ so that $Z \sim p_2$, then $L(\gamma)$ is minimized at $\gamma = 3$. It also holds, when $\delta$ is sufficiently small (see Lemma~\ref{lem:conti_argmin_L}), the likelihood $L(\gamma)$ is also minimized at $\gamma = 3$ so that the likelihood based selector would likely output $\hat{\theta}_3$. However, for any $\delta > 0$, we have that $V(3)$, the asymptotic variance of $\hat{\theta}_3$, is $\frac{\mathbb{E}|Z|^4}{(2 \mathbb{E}|Z|^2)} = \infty$. In contrast, our proposed procedure would output the sample mean $\bar{Y} = \hat{\theta}_2$, which has finite asymptotic variance. Intuitively, the CAVS procedure behaves better because it takes into account the higher moment $\mathbb{E}|Z|^{2(\gamma -1)}$ whereas the likelihood selector is based only on $\mathbb{E}|Z|^\gamma$.


\end{remark}

\section{Empirical studies}
\label{sec:empirical}

We perform empirical studies on simulated data to verify our theoretical results in Section~\ref{sec:theory}. We also analyze a dataset of NBA player statistics for the 2020-2021 season to show that our proposed CAVS estimator can be directly applied to real data. 

\subsection{Simulations}
\label{sec:simulation}

\begin{figure}[htp]
\centering
\begin{subfigure}{0.45\textwidth}
\centering
\includegraphics[scale=.4]{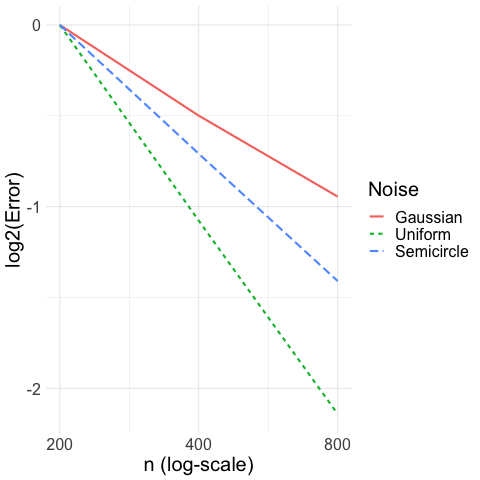}
\caption{Location estimation}
\label{fig:loc_rate}
\end{subfigure}
\begin{subfigure}{0.45\textwidth}
\centering
\includegraphics[scale=.4]{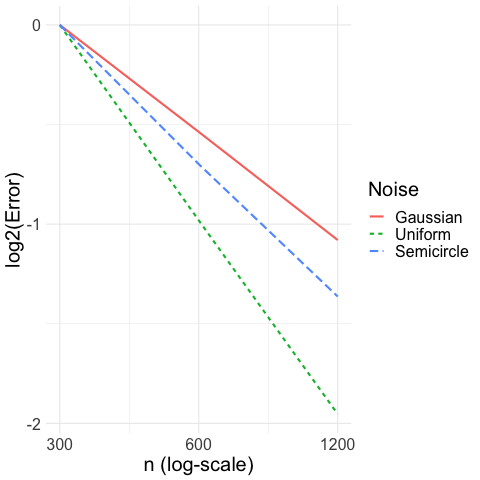}
\caption{Regression}
\label{fig:regr_rate}
\end{subfigure}
\caption{Log-error vs. sample size plots. Sample size $n$ is plotted on a log-scale.}
\label{fig:rate}
\end{figure}

\textbf{Convergence rate for location estimation:} Our first simulation takes the location estimation setting where $Y_i = \theta_0 + Z_i$ for $i = 1, \ldots, n$. We let the distribution of the noise $Z_i$ be either Gaussian $N(0, 1)$, uniform $\text{Unif}[-1, 1]$, or semicircle (see Example~\ref{ex:moment_densities}). We let the sample size $n$ vary between $(200, 400, 800)$. We compute our proposed CAVS estimator  $\hat{\theta}_{\hat{\gamma}}$ (with $\tau = \sqrt{ \log \frac{4n}{200} }$) and plot, in Figure~\ref{fig:loc_rate}, log-error versus the sample size $n$, where $n$ is plotted on a logarithmic scale. Hence, a rate of convergence of $n^{-t}$ would yield an error line of slope $-t$ in Figure~\ref{fig:loc_rate}. We normalize the errors so that all the lines have the same intercept. We see that error under uniform noise has a slope of $-1$, error under semicircle noise has a slope of $-2/3$, and error under Gaussian noise has a slope of $-1/2$ exactly as predicted by Theorem~\ref{thm:basic} and Theorem~\ref{thm:adaptive_rate}.\\

\textbf{Convergence rate for regression:} Then, we study the regression setting where $Y_i = X_i^\top \beta_0 + Z_i$ for $i = 1, 2, \ldots, n$. We let the distribution of the noise $Z_i$ be either Gaussian $N(0, 1)$, uniform $\text{Unif}[-1, 1]$, or the semicircle density given in Example~\ref{ex:moment_densities}. We let the sample size $n$ vary between $(200, 400, 600)$. We apply the regression version of the CAVS estimate $\hat{\beta}_{\hat{\gamma}}$ as described in Remark~\ref{rem:regression} (with $\tau = \sqrt{\log \frac{4n}{200}}$), and plot, in Figure~\ref{fig:loc_rate}, log-error versus the sample size $n$, where $n$ is plotted on a logarithmic scale. We see that CAVS also has adaptive rate of convergence; the uniform noise yields a rate of $n^{-1}$, the semicircle noise yields a rate of $n^{-2/3}$, and the Gaussian noise yields a rate of $n^{-1/2}$ as $n$ increases, as predicted by our theory. \\

\begin{figure}[htp]
\centering
\begin{subfigure}{0.45\textwidth}
\centering
\includegraphics[scale=.4]{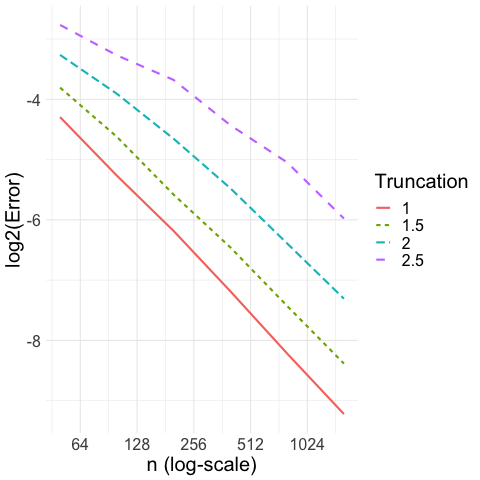}
\caption{Gaussian truncated at different levels}
\label{fig:trunc_rate}
\end{subfigure}
\begin{subfigure}{0.45\textwidth}
\centering
\includegraphics[scale=.4]{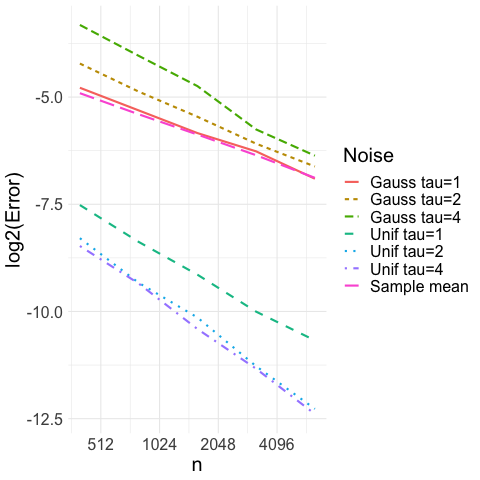}
\caption{Different $\tau$}
\label{fig:tau_rate}
\end{subfigure}
\caption{Log-error vs. sample size plots. Sample size $n$ is plotted on a log-scale.}
\label{fig:rate2}
\end{figure}

\textbf{Convergence rate for truncated Gaussian at different truncation levels:} In Figure~\ref{fig:trunc_rate}, we take the location model $Y_i = \theta_0 + Z_i$ where $Z_i$ has the density $p_t(x) \propto  \exp\{ -\frac{1}{2} \frac{x^2}{\sigma^2_t}\} \mathbbm{1}(|x|\leq t/\sigma_t)$ for some $t > 0$ and where $\sigma_t > 0$ is chosen so that $Z_i$ always has unit variance. In other words, we sample $Z_i$ by first generating $W \sim N(0, 1)$, keep $W$ only if $|W| \leq t$, and then take $Z_i = \sigma_t W$ where $\sigma_t > 0$ is chosen so that $\text{Var}(Z_i) = 1$. We use four different truncation levels $t = 1, 1.5, 2, 2.5$; we let the sample size vary from $n=50$ to $n=1600$ and compute our CAVS estimate $\hat{\theta}_{\hat{\gamma}}$ (with $\tau = \sqrt{ \log \frac{4n}{50} }$). We plot in Figure~\ref{fig:trunc_rate}, the log-error versus the sample size $n$, where $n$ is plotted on a logarithmic scale. We observe that when the truncation level is $t = 1$ or $1.5$ or $2$, the error is of order $n^{-1}$. When the truncation level is $t=2.5$, the error behaves like $n^{-1/2}$ for small $n$ but transitions to $n^{-1}$ when $n$ becomes large. This is not surprising since, when $n$ is small, it is difficult to know whether the $Z_i$'s are drawn from $N(0, 1)$ or drawn from truncated Gaussian with a large truncation level. \\

\textbf{Convergence rate for different $\tau$:} In Figure~\ref{fig:tau_rate}, we take the location model $Y_i = \theta_0 + Z_i$ and take $Z_i$ to be either Gaussian $N(0,1)$ or uniform $\text{Unif}[-M, M]$ where $M > 0$ is chosen so that $Z_i$ has unit variance. We then apply our proposed CAVS procedure for different levels of $\tau$, ranging from $\tau \in \{1, 2, 4\}$. We let the sample size vary from $n=400$ to $n=6400$ and plot the log-error versus the sample size $n$, where $n$ is plotted on a logarithmic scale. For comparison, we also plot the error of the sample mean $\bar{Y}$, which does not depend on the distribution of $Z_i$ since we scale $Z_i$ to have unit variance in both settings. We observe in Figure~\ref{fig:tau_rate} that when $\tau = 1$, the CAVS estimate $\hat{\theta}_{\hat{\gamma}}$ basically coincides with the sample mean if $Z_i \sim N(0,1)$ but has much less error when $Z_i$ is uniform. As we increase $\tau$, CAVS estimator has increased error under the Gaussian setting when $Z_i \sim N(0, 1)$ since we select $\hat{\gamma} > 2$ more often; under the uniform setting, it has less error. Based on these studies, we recommend $\tau = 1$ in practice as a conservative choice.

\subsection{Real data experiments}

Uniform or truncated Gaussian data are not ubiquitous but they do appear in real world datasets. In this section, we use the CAVS location estimation and regression procedure to analyze a dataset of 626 NBA players in the 2020--2021 season. We consider variables AGE, MPG (average minutes played per game), and GP (games played). 

\begin{figure}[htp]
\centering
\begin{subfigure}{0.32\textwidth}
\centering
\includegraphics[scale=.27]{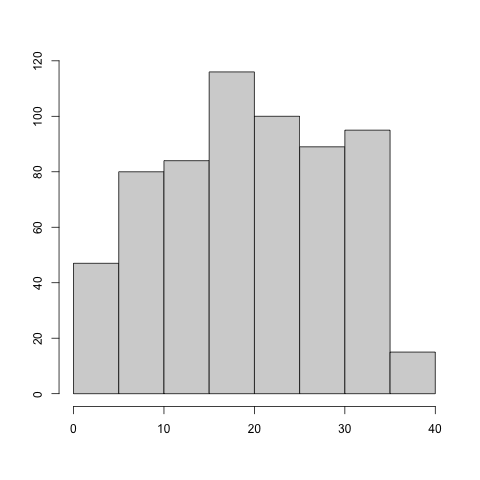}
\caption{Histogram of MPG}
\label{fig:nba1hist}
\end{subfigure}
\begin{subfigure}{0.32\textwidth}
\centering
\includegraphics[scale=.27]{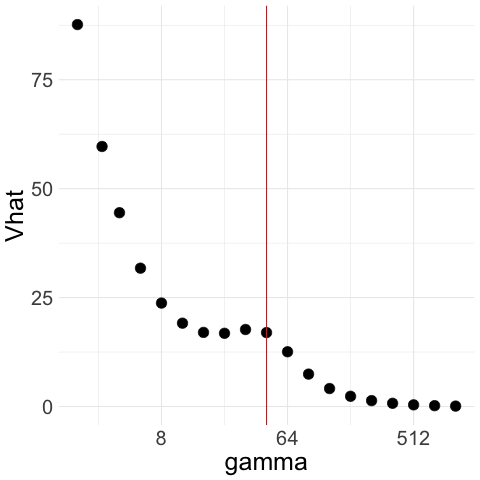}
\caption{$\hat{V}(\gamma)$}
\label{fig:nba1a}
\end{subfigure}
\begin{subfigure}{0.32\textwidth}
\centering
\includegraphics[scale=.27]{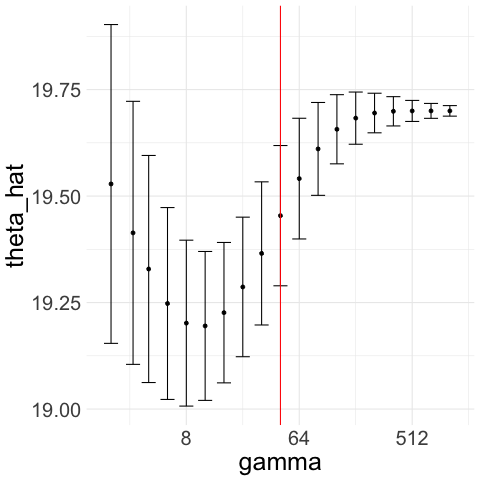}
\caption{$\hat{\theta}_{\gamma} \pm \sqrt{ \frac{\hat{V}(\gamma)}{n}}$}
\label{fig:nba1b}
\end{subfigure}
\caption{Analysis on MPG (average minutes played per game) in the NBA 2021 data}
\label{fig:nba1}
\end{figure}

Both MPG and GP variables are compactly supported. They also do not exhibit clear signs of asymmetry; MPG has an empirical skewness of $-0.064$ and GP has an empirical skewness of $0.013$. We apply the CAVS procedure to both with $\tau = 1$ and we obtain $\hat{\gamma} = 32$ for MPG variable and $\hat{\gamma} = 2048$ for the GP variable. In contrast, the AGE variable has a skewness of $0.56$ and when we apply CAVS procedure (still with $\tau = 1$), we obtain $\hat{\gamma} = 2$. These results suggest that CAVS can be useful for practical data analysis. 

Moreover, we also study the CAVS regression method by considering two regression models:
\[
(\text{MODEL 1}) \,\, \text{MPG} \sim \text{GP} + \text{AGE} + \text{W}, \qquad (\text{MODEL 2}) \,\, \text{MPG} \sim \text{AGE} + \text{W},
\]
where $\text{W}$ is an independent Gaussian feature add so that we can assess how close the estimated coefficient $\hat{\beta}_{\text{W}}$ is to zero to gauge the estimation error. We estimate $\hat{\beta}_{\hat{\gamma}}$ on 100 randomly chosen training data points and report the predictive error on the remaining test data points; we also report the average value of $|\hat{\beta}_{\text{W}}|$, which we would like to be as close to 0 as possible. We perform 1000 trials of this experiment (choosing random training set in each trial) and report the performance of CAVS versus OLS estimator in Table~\ref{tab:errs}. 

\begin{table}[htp]
\centering
\begin{tabular}{|c|c|c|c|c|}
\hline
 & Model 1 Pred. Error & Model 1 $|\hat{\beta}_{\text{W}}|$ & Model 2 Pred. Error & Model 2 $|\hat{\beta}_{\text{W}}|$ \\ 
\hline
CAVS & {\bf 0.686} & {\bf 0.045} & {\bf 0.95} & {\bf 0.082} \\ 
\hline
OLS & 0.689 & 0.140 & 1.04 & 0.205 \\ 
\hline
\end{tabular}
\caption{Comparison of CAVS vs. OLS on two simple regression models.}
\label{tab:errs}
\end{table}

\begin{figure}[htp]
\centering
\begin{subfigure}{0.32\textwidth}
\centering
\includegraphics[scale=.27]{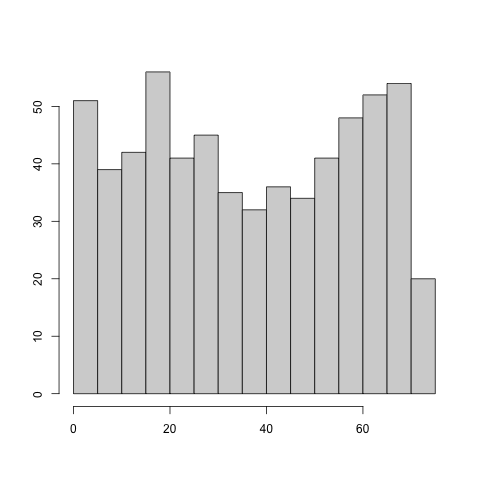}
\caption{Histogram of GP}
\label{fig:nba2hist}
\end{subfigure}
\begin{subfigure}{0.32\textwidth}
\centering
\includegraphics[scale=.27]{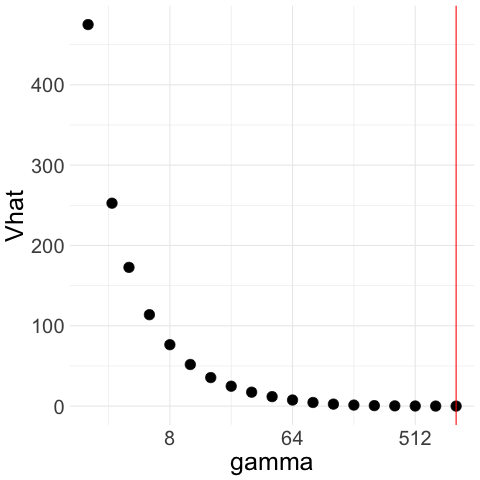}
\caption{$\hat{V}(\gamma)$}
\label{fig:nba2a}
\end{subfigure}
\begin{subfigure}{0.32\textwidth}
\centering
\includegraphics[scale=.27]{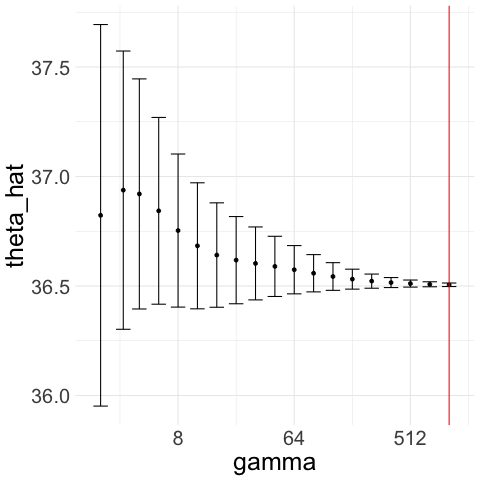}
\caption{$\hat{\theta}_{\gamma} \pm \sqrt{ \frac{\hat{V}(\gamma)}{n}}$}
\label{fig:nba2b}
\end{subfigure}
\caption{Analysis on GP (games played) in the NBA 2021 data}
\label{fig:nba2}
\end{figure}

\section{Discussion}
\label{sec:discussion}

In this paper, we give an estimator of the location of a symmetric distribution whose rate of convergence can be faster than the usual $\sqrt{n}$ and can adapt to the unknown underlying density. There are a number of interesting open questions that remain:

\begin{itemize}
\item It is unclear whether the excess log-factors in our adaptive rate result is an artifact of the analysis or an unavoidable cost of adaptivity.   
\item We emphasize that it is the discontinuity of the noise density $p$ (or the singularities of the derivative $p'$) on the boundary of the support $[-1, 1]$ that allows our estimator to have rates of convergence faster than $1/\sqrt{n}$. Any discontinuities of the noise density (or the singularities of the density derivative) in the \emph{interior} of the support will also lead to an infinite Fisher information (for the location parameter) and open the possibility of faster-than-root-n rate. Our estimator unfortunately cannot adapt to discontinuities in the interior. For example, if the noise density is a mixture of Uniform$[-1,1]$ and Gaussian $N(0,1)$, the tail of the Gaussian component would imply that our estimator cannot have a rate faster than root-n. On the other hand, an oracle with knowledge of the discontinuity points at $\pm 1$ could still estimate $\theta_0$ at $n^{-1}$ rate with the estimator $\argmax_{\theta}\frac{1}{n}\sum_{i=1}^n\mathbbm{1}\{Y_i\in[\theta-1,\theta+1]\}$. One potential approach for adapting to discontinuities in the interior is to first estimate the position of these discontinuity points. However, the position must be estimated with a high degree of accuracy as any error would percolate to the down-stream location estimator. We leave a formal investigation of this line of inquiry to future work. 

\item It would be nontrivial to extend our rate adaptivity result to the multivariate setting, for instance, if $Y_i = \theta_0 + Z_i$ where $\theta_0 \in \mathbb{R}^d$ and $Z_i$ is uniformly distributed on a convex body $K \subset \mathbb{R}^d$ that is balanced in that $K = -K$ so that $Z_i \stackrel{d}{=} -Z_i$. When $K$ is known, it would be natural to study estimators of the form $\hat{\theta}_{\gamma} = \argmin_{\theta \in \mathbb{R}^d} \sum_{i=1}^n \| Y_i - \theta \|_K^\gamma$ where $\| \cdot \|_K$ is the gauge function (Minkowski functional) associated with $K$. In general, it would be necessary to simultaneously estimate $K$ and $\theta_0$. \cite{xu2021high} studies an approach where one first estimates $\theta_0$ via the sample mean and then compute $\hat{K}$ using the convex hull of the directional quantiles of the data. This however cannot achieve a rate faster than root-n.

\item When applied in the linear regression setting, our CAVS procedure performs well empirically on both synthetic and real data. It would thus be interesting to rigorously establish a rate adaptivity result in the linear regression model. More generally, in a nonparametric regression model $Y_i = f(X_i) + Z_i$ where the noise $Z_i$ has a noise distribution symmetric around $0$ and the regression function $f$ lies in some nonparametric function class $\mathcal{F}$, we can still use our procedure to select amongst estimators of the form $\hat{f}_{\gamma} = \argmin_{f \in \mathcal{F}} \sum_{i=1}^n |Y_i - f(X_i)|^{\gamma}$. Understanding the statistical properties of this procedure would motivate the use of $\ell^{\gamma}$ loss functions, for $\gamma > 2$, in general regression problems.

\end{itemize}

\section*{Acknowledgement}

The first and second authors are supported by NSF grant DMS-2113671. The authors are very grateful to Richard Samworth for suggesting the use of Lepski's method. The authors further thank Jason Klusowski for many insightful discussions. 

\bibliographystyle{dcu}
\bibliography{reference}

\clearpage

\setcounter{section}{0}
\setcounter{equation}{0}
\setcounter{theorem}{0}
\def\theequation{S\arabic{section}.\arabic{equation}}
\def\thesection{S\arabic{section}}
\def\thetheorem{S\arabic{section}.\arabic{theorem}}

\begin{center}
\Large{Supplementary material to ``Rate optimal and adaptive estimation of the center of a symmetric distribution''} \\ \vspace{0.2in}
\large{Yu-Chun Kao, Min Xu, and Cun-Hui Zhang}
\end{center}

\begin{spacing}{1.2}

\section{Supplementary material for Section~\ref{sec:method}}
\label{sec:method_appendix}
\subsection{Proof of Lemma~\ref{lem:uniformcontrol}}
\begin{proof} (of Lemma~\ref{lem:uniformcontrol})

First, we observe that if $4 \frac{\log n}{\gamma} \geq 1$, then, by the fact that $\hat{\theta}_{\gamma} \in [X_{(1)}, X_{(n)}]$, we have that 
\[
|\hat{\theta}_{\gamma} - X_{\text{mid}}| \leq \frac{1}{2} (X_{(n)} - X_{(1)}) \leq 2 (X_{(n)} - X_{(1)}) \frac{\log n}{\gamma}.
\]
Therefore, we assume that $4 \frac{\log n}{\gamma} \leq 1$. 

We apply Lemma~\ref{lem:argmin_approx} with $f(\theta) = \bigl\{ \frac{1}{n}\sum_{i=1}^n |X_i - \theta|^{\gamma} \bigr\}^{1/\gamma}$ and $g(\theta) = \max_i | X_i - \theta |$ so that $\theta_g := \argmin g(\theta) = \hat{\theta}_{\text{mid}}$ and $\theta_f := \argmin f(\theta) = \hat{\theta}_{\gamma}$. Fix any $\delta > 0$. We observe that
\begin{align*}
g(\theta_g) &= \max_i | X_i - X_{\text{mid}}| = \frac{X_{(n)}-X_{(1)}}{2} \\
g(\theta_g + \delta) &= \frac{X_{(n)}-X_{(1)}}{2} + \delta, \qquad
\text{ and }
g(\theta_g - \delta) = \frac{X_{(n)}-X_{(1)}}{2} + \delta.
\end{align*}

Therefore, for $\theta \in \{ \theta_g - \delta, \theta_g + \delta\}$, we have that $\frac{1}{2}(g(\theta) - g(\theta_g)) = \frac{\delta}{2}$. On the other hand, by the fact that
\[
\biggl\{ \frac{1}{n}\sum_{i=1}^n |X_i - \theta|^{\gamma} \biggr\}^{\frac{1}{\gamma}} \geq n^{-\frac{1}{\gamma}} \max_{i \in [n]} |X_i - \theta|,
\]
we have that
\begin{align}
g(\theta) \geq f(\theta) \geq n^{- \frac{1}{\gamma}} g(\theta) \quad \forall \theta \in \mathbb{R}. \label{eq:fg_comparison}
\end{align}
Therefore, for $\theta \in \{ \theta_g - \delta, \theta_g, \theta_g + \delta \}$, we have that
\begin{align*}
|f(\theta) - g(\theta)| 
&= g(\theta) - f(\theta) \leq (1 - n^{-\frac{1}{\gamma}}) g(\theta) \\
&\leq \frac{\log n}{\gamma} \biggl( \frac{X_{(n)} - X_{(1)}}{2} + \delta \biggr).
\end{align*}
Using our assumption that $4 \frac{\log n}{\gamma} \leq 1$, we have that for any $\delta \geq 2(X_{(n)} - X_{(1)}) \frac{\log n}{\gamma}$ and any $\theta \in \{\theta_g - \delta, \theta_g, \theta_g + \delta\}$,
\[
|f(\theta) - g(\theta)| \leq \frac{\log n}{\gamma} \biggl( \frac{X_{(n)} - X_{(1)}}{2} + \delta \biggr) \leq \frac{\delta}{2} = \frac{1}{2} (g(\theta) - g(\theta_g)).
\]
The Lemma thus immediately follows from Lemma~\ref{lem:argmin_approx}.
\end{proof}

\subsection{ Lemma~\ref{lem:as_min_converge} on the convergence of $\hat{V}(\gamma)$}
\label{sec:as_min_converge}

The following lemma implies that, for a fixed $\gamma$ such that $V(\gamma)$ is well-defined, our asymptotic variance estimator $\hat{V}(\gamma)$ is consistent. For a random variable $Y$, we define its essential supremum to be
\[
\text{ess-sup}(Y) := \inf \bigl\{ M \in \mathbb{R} \,:\, \mathbb{P}(Y \leq M) = 1 \bigr\},
\]
where the infimum of an empty set is taken to be infinity. Note that $\text{ess-sup}(|Y|) < \infty$ if and only if $Y$ is compacted supported and that $\lim_{\gamma \rightarrow \infty} \{ \mathbb{E}|Y|^{\gamma} \}^{\frac{1}{\gamma}} = \text{ess-sup}(|Y|)$. 

We may define $\text{ess-inf}(Y)$ is the same way. For an infinite sequence $Y_1, Y_2, \ldots$ of independent and identically distributed random variables, it is straightforward to show that $Y_{(n), n} := \max_{i \in [n]} Y_i \stackrel{\text{a.s.}}{\rightarrow} \text{ess-sup}(Y)$ and $Y_{(1), n} := \min_{i \in [n]} Y_i \stackrel{\text{a.s.}}{\rightarrow} \text{ess-inf}(Y)$ regardless of whether the essential supremum and infimum are finite or not.

\begin{lemma}
\label{lem:as_min_converge}
Let $Y_1, Y_2, \ldots $ be a sequence of independent and identically distributed random variables and let $\gamma > 1$. The following hold:
\begin{enumerate}
\item If $\mathbb{E}|Y|^\gamma < \infty$, then $\min_{\theta \in \mathbb{R}} \frac{1}{n} \sum_{i=1}^n |Y_i - \theta|^\gamma \stackrel{a.s.}{\rightarrow} \min_{\theta \in \mathbb{R}} \mathbb{E}|Y - \theta|^\gamma$.
\item If $\mathbb{E}|Y|^\gamma = \infty$, then $\min_{\theta \in \mathbb{R}} \frac{1}{n} \sum_{i=1}^n |Y_i - \theta|^\gamma \stackrel{a.s.}{\rightarrow} \infty$.
\item If $Y$ is compactly supported, then we have that $\min_{\theta \in \mathbb{R}} \max_{i\leq n} |Y_i - \theta|=(Y_{(n), n} - Y_{(1), n})/2 \stackrel{a.s.}{\rightarrow} \min_{\theta} \text{ess-sup}(|Y - \theta|)$.
\item If $\text{ess-sup}(|Y|) = \infty$, then $\min_{\theta \in \mathbb{R}} \max_{i\leq n} |Y_i - \theta|\stackrel{a.s.}{\rightarrow}\infty$.
\end{enumerate}
As a direct consequence, for any $\gamma > 1$ such that $\mathbb{E}|Y|^{\gamma-2} < \infty$, we have $\hat{V}(\gamma) \stackrel{a.s.}{\rightarrow} V(\gamma)$, even when $V(\gamma) = \infty$.
\end{lemma}


\begin{proof} (of Lemma~\ref{lem:as_min_converge}) 

For the first claim, we apply Proposition~\ref{prop:convex_min_converge} with $g(y, \theta) = |y - \theta|^\gamma$ and $\psi(\theta) = \mathbb{E}| Y - \theta |^\gamma$ and immediately obtain the desired conclusion.

We now prove the second claim by a truncation argument. Suppose $\mathbb{E}|Y|^\gamma = \infty$ so that $\min_{\theta} \mathbb{E}|Y - \theta|^\gamma = \infty$. Fix $M > 0$ arbitrarily. We claim there then exists $\tau > 0$ such that 
\[
\min_{\theta \in \mathbb{R}} \mathbb{E}\bigl[|Y - \theta|^\gamma \mathbbm{1}\{ |Y| \leq \tau \} \bigr] > M.
\]
To see this, for any $\tau > 0$, define $\theta_\tau = \argmin_{\theta \in \mathbb{R}} \mathbb{E}\bigl[|Y - \theta|^\gamma \mathbbm{1}\{ |Y| \leq \tau \} \bigr]$. The argmin is well-defined since $\theta \mapsto \mathbb{E}\bigl[|Y - \theta|^\gamma \mathbbm{1}\{ |Y| \leq \tau \} \bigr]$ is strongly convex and goes to infinity as $|\theta| \rightarrow \infty$. If $\{ \theta_{\tau} \}_{\tau=1}^\infty$ is bounded, then the claim follows because $\mathbb{E}|Y - \theta_\tau|^{\gamma} \mathbbm{1}\{|Y|\leq \tau\} \geq \{ (\mathbb{E}|Y|^{\gamma} \mathbbm{1}\{ |Y| \leq \tau \} )^{1/\gamma} - \theta_\tau \}^\gamma$. If $\{\theta_\tau\}_{\tau=1}^\infty$ is unbounded, then there exists a sub-sequence $\tau_m$ such that $\lim_{m \rightarrow \infty} \theta_{\tau_m} \rightarrow \infty$ say. For any $a > 0$ such that $\mathbb{P}(|Y| \leq a) > 0$, we have $\lim_{m \rightarrow \infty} \mathbb{E}|Y - \theta_{\tau_m}|^{\gamma} \mathbbm{1}\{|Y|\leq \tau_m \} \geq \lim_{m \rightarrow \infty} |a - \theta_{\tau_m} |^{\gamma}\mathbb{P}(|Y| \leq a) = \infty$. Therefore, in either cases, our claim holds. 

Using Proposition~\ref{prop:convex_min_converge} again with $g(x, \theta) = |x - \theta|^{\gamma} \mathbbm{1}\{|x| \leq \tau\}$, we have that
\[
\min_{\theta} \frac{1}{n} \sum_{i=1}^n |Y_i - \theta|^{\gamma} \mathbbm{1}\{ |Y_i| \leq \tau\} \stackrel{a.s.}{\rightarrow} \min_{\theta} \mathbb{E} \bigl[|Y - \theta|^\gamma \mathbbm{1}\{ |Y| \leq \tau \} \bigr] > M.
\]
In other words, there exists an event $\tilde{\Omega}_M$ with probability 1 such that, for any $\omega \in \tilde{\Omega}_M$, there exists $n_\omega$ such that for all $n \geq n_{\omega}$, 
\begin{align*}
\min_{\theta} \frac{1}{n} \sum_{i=1}^n |Y_i - \theta|^{\gamma} 
&\geq 
\min_{\theta} \frac{1}{n} \sum_{i=1}^n |Y_i - \theta|^{\gamma} \mathbbm{1}\{ |Y_i| \leq \tau\} \geq M/2.
\end{align*}

Thus, on $\tilde{\Omega}_M$, we have that
\[
\liminf_{n \rightarrow \infty} \min_{\theta} \frac{1}{n} \sum_{i=1}^n |Y_i - \theta|^{\gamma} > M/2.
\]
Thus, on the event $\tilde{\Omega} = \cap_{M=1}^\infty \tilde{\Omega}_M$, we have that $\min_{\theta} \frac{1}{n} \sum_{i=1}^n |Y_i - \theta|^{\gamma} \rightarrow \infty$. Since $\tilde{\Omega}$ has probability 1, the second claim follows.
For the third claim, without loss of generality, we can assume that $\text{ess-sup}(Y) = 1$ and $\text{ess-inf}(Y) = -1$. Define $X_n=(Y_{(n), n}-Y_{(1), n})/2$, then we have $X_n \leq \min_\theta \text{ess-sup}(|Y-\theta| )=1$ and
\begin{align*}
    \mathbb{P}\{X_n<1-\delta\} \leq\mathbb{P}\{Y_{(n),n}<1-\delta\} +\mathbb{P}\{Y_{(1),n}>-1+\delta\},
\end{align*}where, as $n \rightarrow \infty$, the right hand side tends to $0$ for every $\delta>0$. $X_n$ thus converges to $1$ in probability. Since the collection $\{X_n\}_{n=1}^\infty$ is defined on the same infinite sequence $\{ Y_1, Y_2, \ldots \}$ of independent and identically distributed random variables, we have that $1 \geq X_n \geq X_{n-1} \geq 0$ so that $X_n \stackrel{a.s.}{\rightarrow } 1$ by the monotone convergence theorem.

For the forth claim, suppose without loss of generality that $\text{ess-inf}(Y) \leq -1$ and that $\text{ess-sup}(Y) = \infty$. Let $X_n = (Y_{(n),n} - Y_{(1),n})/2$ as with the proof of the third claim. Then, 
\begin{align*}
    \mathbb{P}\{X_n < M\} \leq \mathbb{P}\{Y_{(n),n} < 2M \} + \mathbb{P}\{Y_{(1),n} \geq 0 \}.
\end{align*}
Since the right hand side tends to $0$ for every $M>0$, we have that $X_n$ converges to infinity almost surely. The Lemma follows as desired. 
\end{proof}

\subsection{Bound on $V$}

The following lower bound on $V(\gamma)$ holds regardless of whether $Y$ is symmetric around $\theta_0$ or not. We have
\begin{align*}
V(\gamma) &= \frac{ \mathbb{E}| Y - \theta_0|^{2(\gamma-1)} }{(\gamma-1)^2 \{ \mathbb{E} |Y - \theta_0|^{\gamma-2} \}^2 } \\
&= \frac{ \mathbb{E}| Y - \theta_0|^{2(\gamma-1)} }{ \{ \mathbb{E}| Y - \theta_0|^{\gamma-1} \}^2 } 
\biggl( \frac{ \mathbb{E}|Y - \theta_0|^{\gamma-1} }{\mathbb{E} |Y - \theta_0|^{\gamma-2} } \biggr)^2 \frac{1}{(\gamma-1)^2} \\
&\geq \frac{ \mathbb{E}| Y - \theta_0|^{2(\gamma-1)} }{ \{ \mathbb{E}| Y - \theta_0|^{\gamma-1} \}^2 } 
\{ \mathbb{E} |Y - \theta_0|^{\gamma-2} \}^{\frac{2}{\gamma - 2}} \frac{1}{(\gamma-1)^2} \\
&\geq \frac{ \mathbb{E}| Y - \theta_0|^{2(\gamma-1)} }{ \{ \mathbb{E}| Y - \theta_0|^{\gamma-1} \}^2 } \mathbb{E}|Y-\theta_0|^2 \frac{1}{(\gamma-1)^2},
\end{align*}

where the first inequality follows from the fact that $\mathbb{E}|Y - \theta_0|^{\gamma-1} \geq \bigl( \mathbb{E} |Y-\theta_0|^{\gamma-2} \bigr)^{\frac{\gamma-1}{\gamma-2}}$. In particular, we have that $V(\gamma) \geq \frac{\mathbb{E}|Y - \theta_0|^2}{(\gamma-1)^2}$. Equality is attained when $Y - \theta_0$ is a Rademacher random variable.

\subsection{Optimization algorithm}
\label{sec:optimization_appendix}

We give the Newton's method algorithm for computing $\hat{\theta}_\gamma = \argmin_{\theta \in \mathbb{R}} \frac{1}{n} \sum_{i=1}^n |Y_i - \theta|^\gamma$. It is important to note that to avoid numerical precision issues when $\gamma$ is large, we have to transform the input $Y_1, \ldots, Y_n$ so that they are supported on the unit interval $[-1, 1]$. 

\begin{algorithm}[htp]
\caption{Newton's method for location estimation}
\label{alg:newton}
\begin{flushleft}
\textbf{INPUT}: observations $Y_1, \ldots, Y_n \in \mathbb{R}$ and $\gamma \geq 2$. \\
\textbf{OUTPUT:} $\hat{\theta}_\gamma := \min_{\theta \in \mathbb{R}} \frac{1}{n} \sum_{i=1}^n |Y_i - \theta|^\gamma$.
\begin{algorithmic}[1]
\State Compute $S = (Y_{(n)} - Y_{(1)})$ and $M = (Y_{(n)}+Y_{(1)})/2$ and transform $Y_i \leftarrow 2(Y_i - M)/S$.
\State Initialize $\theta^{(0)} = 0$.
\For{$t = 1, 2, 3, \ldots $}
    \State Compute $f' = - \frac{1}{n} \sum_{i=1}^n |Y_i - \theta^{(t-1)} |^{\gamma-1} \text{sign}(Y_i - \theta^{(t-1)})$
    \State Compute $f'' = \frac{\gamma-1}{n} \sum_{i=1}^n |Y_i - \theta^{(t-1)}|^{\gamma-2}$.
    \State Set $\theta^{(t)} = \theta^{(t-1)} - \frac{f'}{f''}$.
    \State If $|f'| \leq \varepsilon$, break and output $S \theta^{(t)} /2 + M$.
\EndFor
\end{algorithmic}
\end{flushleft}
\end{algorithm}

To compute $\hat{\theta}_{\gamma}$ for a collection of $\gamma_1 < \gamma_2 < \ldots $, we can warm start our optimization of $\hat{\theta}_{\gamma_2}$ by initializing with $\hat{\theta}_{\gamma_1}$. In the regression setting where $\gamma$ is large, we find that it improves numerical stability to to apply a quasi-Newton's method where we add a an identity $\varepsilon I$ to the Hessian for a small $\varepsilon > 0$. 

\subsection{Supporting Lemmas}

\begin{lemma}
\label{lem:argmin_approx}
Let $f, g \,:\, \mathbb{R}^d \rightarrow \mathbb{R}$ and suppose $f$ is convex. Let $x_g \in \argmin g(x)$ and $x_f \in \argmin f(x)$. Suppose there exists $\delta > 0$ such that
\begin{align*}
|f(x) - g(x)| \vee |f(x_g) - g(x_g)| < \frac{1}{2} (g(x) - g(x_g)), \quad \text{ for all $x$ s.t. $\| x - x_g \| = \delta$}.
\end{align*}
Then, we have that
\[
\| x_f - x_g \| \leq \delta. 
\]
\end{lemma}

\begin{proof}

Let $\delta > 0$ and suppose $\delta$ satisfies the condition of the Lemma. Fix $x \in \mathbb{R}^d$ such that $\| x - x_g \| > \delta$. Define $\xi = x_g + \frac{\delta}{\| x - x_g \|} (x - x_g)$ so that $\| \xi - x_g \| = \delta$. Note by convexity of $f$ that $f(\xi) \leq (1 - \frac{\delta}{\| x - x_g \|}) f(x_g) + \frac{\delta}{\|x - x_g\|} f(x)$. 

Therefore, we have that
\begin{align*}
\frac{\delta}{\| x - x_g \|} (f(x) - f(x_g)) &\geq f(\xi) - f(x_g) \\
&= f(\xi) - g(\xi) + g(\xi) - g(x_g) + g(x_g) - f(x_g) > 0
\end{align*}
under the condition of the Theorem. Therefore, we have $f(x) > f(x_g)$ for any $x$ such that $\| x - x_g \| > \delta$. The conclusion of the Theorem follows as desired.
\end{proof}
\subsubsection{LLN for minimum of a convex function}

\begin{proposition}
\label{prop:convex_min_converge}
Suppose $\theta \mapsto g(y, \theta)$ is convex on $\mathbb{R}$ for all $y \in \mathcal{Y}$. Define $\psi(\theta) := \mathbb{E} g(Y, \theta)$ and suppose $\psi$ is finite on an open subset of $\mathbb{R}$ and $\lim_{|\theta| \rightarrow \infty} \psi(\theta) = \infty$.

Then, we have that
\[
\min_{\theta \in \mathbb{R}} \frac{1}{n} \sum_{i=1}^n g(Y_i, \theta) \stackrel{a.s.}{\rightarrow} \min_{\theta \in \mathbb{R}} \mathbb{E} g(Y, \theta),
\]and
\[
\sup_{\theta_1\in\Theta_n}\min_{\theta_2\in\Theta_0}|\theta_1-\theta_2|\stackrel{a.s.}{\rightarrow} 0,
\]where $\Theta_n:=\argmin_{\theta \in \mathbb{R}} \frac{1}{n} \sum_{i=1}^n g(Y_i, \theta)$ and $\Theta_0:=\argmin_{\theta \in \mathbb{R}} \mathbb{E} g(Y, \theta)$
\end{proposition}

\begin{proof}
Define $\hat{\psi}_n(\theta) := \frac{1}{n} \sum_{i=1}^n g(Y_i, \theta)$ and observe that $\hat{\psi}_n$ is a convex function on $\mathbb{R}$. We also observe that $\argmin \psi$ is a closed bounded interval on $\mathbb{R}$ and we define $\theta_0$ to be its midpoint. 

Fix $\epsilon > 0$ arbitrarily. We may then choose $\theta_L \in (-\infty, \theta_0)$ and $\theta_R \in (\theta_0, \infty)$ such that
\begin{enumerate}
\item $\psi(\theta_L) > \psi(\theta_0)$ and $\psi(\theta_R) > \psi(\theta_0)$, 
\item $(\psi(\theta_L) - \psi(\theta_0)) \vee (\psi(\theta_R) - \psi(\theta_0)) \leq \epsilon$,
\item $\theta_0 - \theta_L = \theta_R - \theta_0$,
\item and $\min_{\theta\in\Theta_0}|\theta_R-\theta|\vee \min_{\theta\in\Theta_0}|\theta_L-\theta|<\epsilon$.
\end{enumerate}

Define $\tilde{\epsilon} := (\psi(\theta_L) - \psi(\theta_0)) \wedge (\psi(\theta_R) - \psi(\theta_0))$ and note that $0 < \tilde{\epsilon} < \epsilon$ by our choice of $\theta_L$ and $\theta_R$. By LLN, there exists an event $\tilde{\Omega}_{\epsilon}$ with probability 1 such that, for every $\omega \in \tilde{\Omega}_{\epsilon}$, there exists $n_{\omega} \in \mathbb{N}$ where for all $n \geq n_{\omega}$,
\[
|\hat{\psi}_n(\theta_L) - \psi(\theta_L) | \vee |\hat{\psi}_n(\theta_R) - \psi(\theta_R) | \vee |\hat{\psi}_n(\theta_0) - \psi(\theta_0) | \leq \tilde{\epsilon}/3.
\]

Fix any $\omega \in \tilde{\Omega}_{\epsilon}$ and fix $n \geq n_{\omega}$, we have that $\hat{\psi}_n(\theta_L) \geq \psi(\theta_L) - \tilde{\epsilon}/3 > \psi(\theta_0)$ and likewise for $\hat{\psi}_n(\theta_R)$. Thus, $\hat{\psi}_n$ must attain its minimum in the interval $(\theta_L, \theta_R)$, i.e., $\sup_{\theta_1\in\Theta_n}\min_{\theta_2\in\Theta_0}|\theta_1-\theta_2|<\epsilon$. We then have by Lemma~\ref{lem:three_points1} that
\begin{align*}
\min_{\theta \in \mathbb{R}} \hat{\psi}(\theta)
&= \min_{\theta \in (\theta_L, \theta_R)} \hat{\psi}_n(\theta)\\
&\geq \hat{\psi}_n(\theta_0) - |\hat{\psi}_n(\theta_0) - \hat{\psi}_n(\theta_R) | \vee |\hat{\psi}_n(\theta_0) - \hat{\psi}_n(\theta_L)| \\
&\geq \psi(\theta_0) - \tilde{\epsilon} - \epsilon \geq \psi(\theta_0) - 2 \epsilon.
\end{align*}

On the other hand, 
\[
\min_{\theta \in \mathbb{R}} \hat{\psi}_n(\theta) \leq \hat{\psi}_n(\theta_0) \leq \psi(\theta_0) + \epsilon.
\]

Therefore, for all $\omega \in \tilde{\Omega}_{\epsilon}$, we have that
\[
\limsup_{n \rightarrow \infty} \bigl| \min_{\theta \in \mathbb{R}} \hat{\psi}_n(\theta) - \psi(\theta_0) \bigr| \leq 2 \epsilon,
\]and
\[
\limsup_{n \rightarrow \infty} \sup_{\theta_1\in\Theta_n}\min_{\theta_2\in\Theta_0}|\theta_1-\theta_2|<\epsilon.
\]

We then define $\tilde{\Omega} := \cap_{k=1}^\infty \tilde{\Omega}_{1/k}$ and observe that $\tilde{\Omega}$ has probability 1 and that on $\tilde{\Omega}$, 
\[
\lim_{n \rightarrow \infty} \bigl| \min_{\theta \in \mathbb{R}} \hat{\psi}_n(\theta) - \psi(\theta_0)\bigr| = 0,
\]and
\[
\lim_{n \rightarrow \infty} \sup_{\theta_1\in\Theta_n}\min_{\theta_2\in\Theta_0}|\theta_1-\theta_2|=0.
\]
The Proposition follows as desired.

\end{proof}

\begin{lemma}
\label{lem:three_points1}
Let $f \,:\, \mathbb{R} \rightarrow \mathbb{R}$ be a convex function. For any $x_0 \in \mathbb{R}$, $x_L \in (-\infty, x_0)$ and $x_R \in (x_0, \infty)$, we have 
\begin{align*}
\text{for all $x \in (x_L, x_0)$, } f(x) &\geq f(x_0) + \{f(x_R) - f(x_0)\} \frac{x - x_0}{x_R - x_0} \\
\text{for all $x \in (x_0, x_R)$, } f(x) &\geq f(x_0) + \{f(x_0) - f(x_L)\} \frac{x - x_0}{x_0 - x_L}.
\end{align*}
As a direct consequence, if $x_0 - x_L = x_R - x_0$, then we have that for all $x \in (x_L, x_0)$, 
\[
f(x) \geq f(x_0) - |f(x_0) - f(x_R)|
\]
and that for all $x \in (x_0, x_R)$,
\[
f(x) \geq f(x_0) - |f(x_0) - f(x_L)|.
\]
\end{lemma}

\begin{proof}
Let $x \in (x_L, x_0)$; using the fact that $f'(x_0) \leq \frac{f(x_R) - f(x_0)}{x_R - x_0}$, we have
\begin{align*}
f(x) &\geq f(x_0) + f'(x_0)(x - x_0) \\
&\geq f(x_0) + \{ f(x_R) - f(x_0)\} \frac{x - x_0}{x_R - x_0}.
\end{align*}

Likewise, for $x \in (x_0, x_R)$, we have $f'(x_0) \geq \frac{f(x_0) - f(x_L)}{x_0 - x_L}$ and hence,
\begin{align*}
f(x) &\geq f(x_0) + f'(x_0)(x - x_0) \\
&\geq f(x_0) + \{ f(x_0) - f(x_L)\} \frac{x - x_0}{x_0 - x_L}.
\end{align*}

\end{proof}
\section{Supplementary material for Section~\ref{sec:theory}}

\subsection{Proof of Theorem~\ref{thm:adaptive_rate}}
\label{sec:adaptive_rate_proof}

\textbf{Structure of intermediate results:} The proof is long and uses various intermediate technical results. The key intermediate theorems are (1) Theorem~\ref{thm:thetahat_gamma_bound} which is essentially a corollary of Proposition~\ref{prop:tailbound1} and (2) Theorem~\ref{thm:vhat_bound} which follow from Proposition~\ref{prop:tailbound2} as well as Theorem~\ref{thm:thetahat_gamma_bound}.\\

\textbf{Notation for constants:} For all the proofs in this section, we let $C$ indicate a generic universal constant whose value could change from instance to instance. We let $C_1, C_2, C_3, C_4$ be specific universal constants where $C_1, C_2$ are defined in the proof of Proposition~\ref{prop:tailbound1} and where $C_3, C_4$ are defined in Theorem~\ref{thm:vhat_bound}. \\

\begin{proof} (of Theorem~\ref{thm:adaptive_rate})\\

We first prove the following: assume that $n$ is large enough such that
\begin{align}
\tau \geq \frac{ C_1 \sqrt{C_4} a_2^{3/2}}{a_1^{3/2}} \sqrt{\log \log n}, \quad \text{and that} \nonumber \\
\biggl\{ \frac{1}{C_{1\alpha} \vee C_{2\alpha} \vee C_4} \bigl( \frac{c_0^2 a_1^6}{a^3_2} \alpha^2 \bigr) \frac{n}{\log n} \biggr\}^{\frac{1}{\alpha}} \geq e^{C_1 \frac{a_2}{a_1}} \geq 2, \label{eq:n_large}
\end{align}
where $C_1, C_4$ are universal constants and $C_{1\alpha}, C_{2\alpha}$ are constants depending only on $\alpha$ -- the value of these are specified in Theorem~\ref{thm:thetahat_gamma_bound} and Theorem~\ref{thm:vhat_bound}.

We claim that 
\begin{align}
\mathbb{P}\biggl\{ | \hat{\theta}_{\hat{\gamma}} - \theta_0 | \leq C_{a_1, a_2,\alpha} \biggl( \frac{\log^{1 + \frac{1}{\alpha}} n}{n^\frac{1}{\alpha}} \vee \frac{\log n}{M_n} \biggr) \biggr\} \geq 1 - \frac{4}{n^{\frac{1}{\alpha}}} - \exp( - \frac{1}{\alpha} (\tau \wedge \sqrt{\log n}) \sqrt{\log n}).
\label{eq:err_in_prob}
\end{align}
This immediately proves the first claim of the theorem. To see that the second claim of the theorem also holds, note that if~\eqref{eq:err_in_prob} holds and if $\tau \geq \sqrt{\log n}$, then, by inflating the constant $C_{a_1, a_2, \alpha}$ if necessary, we have that, for all $n \in \mathbb{N}$,
\begin{align*}
\mathbb{E} | \hat{\theta}_{\hat{\gamma}} - \theta_0 | &\leq C_{a_1, a_2,\alpha} \biggl( \frac{\log^{1 + \frac{1}{\alpha}} n}{n^\frac{1}{\alpha}} \vee \frac{\log n}{M_n} \biggr) + \frac{7}{n^{1/\alpha}} \\
&\leq C_{a_1, a_2,\alpha} \biggl( \frac{\log^{1 + \frac{1}{\alpha}} n}{n^\frac{1}{\alpha}} \vee \frac{\log n}{M_n} \biggr),
\end{align*}
where the first inequality uses the fact that $| \hat{\theta}_{\hat{\gamma}} - \theta_0| \leq 1$. The desired conclusion would then immediately follow. \\

We thus prove~\eqref{eq:err_in_prob} under assumption~\eqref{eq:n_large}. To that end, let $c_0 = 2^{-8}$ and define $\gamma_u=\bigl\{ \frac{1}{C_{1\alpha} \vee C_{2\alpha} \vee C_4} \bigl( \frac{c_0^2 a_1^6}{a^3_2} \alpha^2 \bigr) \frac{n}{\log n} \bigr\}^{\frac{1}{\alpha}}$ and note that $\gamma_u \geq e^{C_1 \frac{a_2}{a_1}} \geq 2$ under assumption~\eqref{eq:n_large}. 

Let $C_4$ be a sufficiently large universal constant as defined in Theorem~\ref{thm:vhat_bound} and define the event
\begin{align}
\mathcal{E}_1 := \biggl\{ \frac{1}{C_4} \frac{a_1}{a_2^2} \gamma^{\alpha-2} \leq \hat{V}(\gamma) \leq C_4 \frac{a_2}{a_1^2} \gamma^{\alpha-2},\quad \text{ for all $\gamma \in [2, (\gamma_u+1)/2]$} \biggr\},
\label{eq:event1a}
\end{align}
It holds by Theorem~\ref{thm:vhat_bound} that $\mathbb{P}( \mathcal{E}_1) \geq 1 - 2 n^{ - \frac{1}{\alpha}}$. 

Now define $\tau' = \frac{1}{\sqrt{C_4}} \frac{\sqrt{a_1}}{a_2} \tau$ and note that $\tau' \geq \frac{C_1 \sqrt{a_2}}{a_1} \sqrt{\log \log n}$ under assumption~\eqref{eq:n_large}. Define the event
\begin{align}
\mathcal{E}_2 := \biggl\{ | \hat{\theta}_{\gamma} - \theta_0 | \leq \tau' \sqrt{ \frac{\gamma^{\alpha-2} }{n}},\, \text{ for all $\gamma \in [2, \gamma_u]$ } \biggr\}.
\end{align}

Then we have by Theorem~\ref{thm:thetahat_gamma_bound} that
\begin{align*}
\mathbb{P}(\mathcal{E}^c_2) &\leq \exp \bigl\{ - \frac{a_1^2}{C_2 a_2} (\tau' \wedge \sqrt{\log n}) \sqrt{ \frac{n}{\gamma_u^\alpha}} \bigr\} \\
&\leq \exp \bigl\{ - \frac{1}{\alpha} \frac{\sqrt{C_4} a_2}{\sqrt{a_1}} (\tau' \wedge \sqrt{\log n}) \sqrt{\log n} \bigr\} \\
&\leq \exp \bigl\{ - \frac{1}{\alpha} (\tau \wedge \sqrt{\log n}) \sqrt{\log n} \bigr\}.
\end{align*}

On the event $\mathcal{E}_1 \cap \mathcal{E}_2$, we have that, for all $\gamma \in [2, (\gamma_u+1)/2]$,
\[
| \hat{\theta}_{\gamma} - \theta_0| \leq 
\tau' \sqrt{ \frac{\gamma^{\alpha-2}}{n}} \leq 
\tau' \sqrt{C_4} \frac{a_2}{\sqrt{a_1}} \sqrt{ \frac{\hat{V}(\gamma)}{n}}
\leq \tau \sqrt{ \frac{\hat{V}(\gamma)}{n} }.
\]

Therefore, we have that 
\begin{align*}
\theta_0 \in \bigcap_{\gamma \in \mathcal{N}_n,\, \gamma \leq (\gamma_u + 1)/2} \biggl[ \hat{\theta}_{\gamma} - \tau \sqrt{ \frac{\hat{V}(\gamma)}{n}}, \, \hat{\theta}_{\gamma} + \tau \sqrt{ \frac{\hat{V}(\gamma)}{n}} \biggr].
\end{align*}
Since $\mathcal{N}_n$ contains $\{2^k \,:\, k \leq \log_2 M_n\}$, either $\frac{\gamma_u+1}{2} \geq M_n$ or there exists $\gamma \in \mathcal{N}_n$ such that $\gamma \geq \frac{\gamma_u+1}{4}$. In either case, it holds by the definition of $\gamma_{\max}$ that
$\gamma_{\max} \geq \frac{\gamma_u+1}{4} \wedge M_n$.  Write $\tilde{\gamma} := \frac{\gamma_u+1}{4} \wedge M_n$.
For any $\gamma < \frac{1}{C_4^2} \bigl(\frac{a_1}{a_2}\bigr)^{\frac{3}{2-\alpha}} \tilde{\gamma} $, we have
\begin{align*}
\hat{V}(\gamma) \geq \frac{1}{C_4} \frac{a_1}{a_2^2} \gamma^{\alpha - 2} > C_4 \frac{a_2}{a_1^2} \tilde{\gamma}^{\alpha-2} \geq \hat{V}(\tilde{\gamma}).
\end{align*}
Since $\hat{\gamma} = \argmin_{\gamma \in \mathcal{N}_n ,\, \gamma \leq \gamma_{\max}} \hat{V}(\gamma)$ and since $\frac{1}{C_4^2} \bigl(\frac{a_1}{a_2}\bigr)^{\frac{3}{2-\alpha}} \leq 1$ so that there exists $\gamma \in \mathcal{N}_n$ such that $\gamma_{\max} \geq \gamma \geq \frac{1}{C_4^2} \bigl(\frac{a_1}{a_2}\bigr)^{\frac{3}{2-\alpha}} \tilde{\gamma}$, it must be that
\begin{align*}
\hat{\gamma} \geq \frac{1}{C_4^2} \bigl(\frac{a_1}{a_2}\bigr)^{\frac{3}{2-\alpha}} \tilde{\gamma} \geq
\frac{1}{C_4^2} \bigl(\frac{a_1}{a_2}\bigr)^{\frac{3}{2-\alpha}} \biggl( \frac{\gamma_u+1}{2} \wedge M_n \biggr) \geq
\tilde{C}^{-1}_{a_1, a_2, \alpha} \biggl\{ \biggl( \frac{n}{\log n} \biggr)^{\frac{1}{\alpha}} \wedge M_n \biggr\},
\end{align*}
where we define $\tilde{C}^{-1}_{a_1, a_2, \alpha} := \frac{1}{4 C_4^2} \bigl( \frac{a_1}{a_2} \bigr)^{\frac{3}{2-\alpha}} \bigl( \frac{1}{C_{1\alpha} \vee C_{2\alpha} \vee C_4} \frac{c_0^2 a_1^4}{a_2} \alpha^2 \bigr)^{\frac{1}{\alpha}}$. 

Now define $\mathcal{E}_3$ as the event that $| \hat{\theta}_{\text{mid}} - \theta_0 | \leq 2^{2 + \frac{2}{\alpha}} a_1^{- \frac{1}{\alpha}} \frac{1}{\alpha^{ \frac{1}{\alpha} + 1}} \frac{\log^{1 + \frac{1}{\alpha}} n}{n^{\frac{1}{\alpha}}}$. We have by Corollary~\ref{cor:midrange} that $\mathbb{P}( \mathcal{E}_3 ) \geq 1 - \frac{2}{n^{1/\alpha}}$. Therefore, on the event $\mathcal{E}_1 \cap \mathcal{E}_2 \cap \mathcal{E}_3$, we have by Lemma~\ref{lem:uniformcontrol} that
\begin{align*}
| \hat{\theta}_{\hat{\gamma}} - \theta_0 |
&\leq | \hat{\theta}_{\hat{\gamma}} - \hat{\theta}_{\text{mid}} | + | \hat{\theta}_{\text{mid}} - \theta_0 | \\
&\leq 4 \frac{\log n}{\hat{\gamma}} + 2^{2 + \frac{2}{\alpha}} a_1^{- \frac{1}{\alpha}} \frac{1}{\alpha^{\frac{1}{\alpha}+1}} \frac{\log^{1 + \frac{1}{\alpha}} n}{n^{\frac{1}{\alpha}}} \\
&\leq 
\tilde{C}_{a_1, a_2, \alpha} \biggl\{ \frac{\log^{1+\frac{1}{\alpha}} n}{ n^{\frac{1}{\alpha} } } \vee \frac{\log n}{M_n} \biggr\}  
+ 2^{2 + \frac{2}{\alpha}} a_1^{- \frac{1}{\alpha}} \frac{1}{\alpha^{\frac{1}{\alpha}+1}} \frac{\log^{1 + \frac{1}{\alpha}} n}{n^{\frac{1}{\alpha}}} \\
&\leq 
C_{a_1, a_2, \alpha} \biggl\{ \frac{\log^{1+\frac{1}{\alpha}} n}{ n^{\frac{1}{\alpha} } } \vee \frac{\log n}{M_n} \biggr\},
\end{align*}
where, in the final inequality, we define $C_{a_1, a_2, \alpha} := \tilde{C}_{a_1, a_2, \alpha} + 2^{2 + \frac{2}{\alpha}} a_1^{- \frac{1}{\alpha}} \frac{1}{\alpha^{\frac{1}{\alpha}+1}}$.

Since $\mathbb{P}( \mathcal{E}_1 \cap \mathcal{E}_2 \cap \mathcal{E}_3 ) \geq 1 - \frac{4}{n^{1/\alpha}} - \exp\bigl( - \frac{1}{\alpha} (\tau \wedge \sqrt{\log n}) \sqrt{\log n} \bigr)$, the desired conclusion~\eqref{eq:err_in_prob} follows. Hence, the Theorem follows as well. 

\end{proof}

\begin{theorem}
\label{thm:thetahat_gamma_bound}
Let $Z_1, \ldots, Z_n$ be independent and identically distributed random variables on $\mathbb{R}$ with a distribution $P$ symmetric around $0$ and write $\nu_{\gamma} := \mathbb{E} |Z|^{\gamma}$. Suppose there exists $\alpha \in (0, 2)$ and $a_1 \in (0, 1]$ and $a_2 \geq 1$ such that $\frac{a_1}{\gamma^\alpha} \leq \nu_{\gamma} \leq \frac{a_2}{\gamma^\alpha}$ for all $\gamma \geq 1$.

Let $C_1, C_2 > 0$ be universal constants and $C_{1\alpha} > 0$ be a constant depending only on $\alpha$, as defined in Proposition~\ref{prop:tailbound1}. Let $c_0 \in (0, 2^{-8})$, let $\gamma_u^{\alpha} = \frac{1}{C_{1\alpha} \vee C_{2,\alpha}} \bigl( \frac{c^2_0 a_1^6}{ a^3_2} \alpha^2 \bigr) \frac{n}{\log n}$, and let $\tau' \geq \frac{C_1 \sqrt{a_2}}{a_1} \sqrt{\log \log n}$.

Suppose $n$ is large enough so that $\gamma_u \geq 2$. Then, we have that
\begin{align}
\mathbb{P}\biggl\{ \sup_{\gamma \in [2, \gamma_u]} \frac{4\gamma}{a_1 c_0} | \hat{\theta}_{\gamma} - \theta_0 | \geq 1 \biggr\} \leq n^{-\frac{1}{\alpha}}.
\end{align}
Moreover, if $n$ is large enough such that $\gamma_u \geq e^{C_1 \frac{a_2}{a_1}} \geq 2$ and that $ \sqrt{\log n} \geq \frac{C_1 \sqrt{a_2}}{a_1} \sqrt{\log \log n}$. Then, we also have
\begin{align}
\mathbb{P}\biggl\{ \sup_{\gamma \in [2, \gamma_u]} \frac{  | \hat{\theta}_{\gamma} - \theta_0 |}{ \tau' \sqrt{ \frac{\gamma^{\alpha-2}}{n}}} \geq 1 \biggr\} 
\leq \exp \biggl\{ - \frac{a_1^2}{C_2 a_2} \biggl(\tau' \wedge \sqrt{\log n} \biggr) \sqrt{ \frac{n}{\gamma_u^\alpha}} \biggr\}
\end{align}

\end{theorem}

\begin{proof}

Since $\hat{\theta}_{\gamma}$ for any $\gamma \geq 2$ is location equivariant, we assume without loss of generality that $\theta_0 = 0$ so that $Y_i = Z_i$. 

Define $\tilde{\tau} = \tau' \wedge \sqrt{\log n}$ and note that $\frac{C_1 \sqrt{a_2}}{a_1} \sqrt{\log \log n} \leq \tilde{\tau} \leq \sqrt{\log n} \leq \frac{1}{4} \sqrt{ \frac{n}{\gamma_u^\alpha}}$ since $a_1 \leq 1$, $a_2 \geq 1$, and $c_0 \leq 2^{-8}$. We further note that
with our definition of and assumptions, the conditions in Proposition~\ref{prop:tailbound1} (i) and (ii) are all satisfied. 

Let $\{ \Delta_{\gamma} \}_{\gamma \geq 2}$ be a collection of positive numbers. For any $\gamma \geq 2$, we have by the second claim of Lemma~\ref{lem:tildetheta_moment} that, for $t \in \{-\Delta_{\gamma}, \Delta_{\gamma}\}$, 
\begin{align}
\bigl| \mathbb{E} \bigl[ - \text{sgn}(Z - t )|Z - t |^{\gamma - 1} \bigr] \bigr| \geq  \frac{a_1}{2} \Delta_{\gamma} \gamma^{1- \alpha}. \label{eq:denom_bound}
\end{align}

To prove the first claim of the theorem, we let $\Delta_{\gamma} = \frac{a_1 c_0}{4}$. We use Proposition~\ref{prop:tailbound1} (noting that the probability bound in~\eqref{eq:psi_tailbound1} is less than $\exp\{ - \frac{1}{\alpha} \log n\}$ under our definition of $\gamma_u$) and~\eqref{eq:denom_bound} to obtain that, with probability at least $1 - n^{- \frac{1}{\alpha}}$, the following holds simultaneously for all $\gamma \in [2, \gamma_u]$:
\begin{align*}
\frac{1}{n} \sum_{i=1}^n \bigl\{ - \text{sgn}(Y_i - \Delta_{\gamma}) |Y_i - \Delta_{\gamma} |^\gamma \bigr\} \geq \frac{1}{2} \mathbb{E} \bigl[ - \text{sgn}(Y - \Delta_{\gamma} )|Y - \Delta_{\gamma} |^{\gamma - 1} \bigr] > 0,
\end{align*}
where, in the last inequality, we use the fact that the function $\theta \mapsto \mathbb{E} |Y - \theta|^{\gamma}$ is strongly convex for all $\gamma > 1$ and minimized at $\theta = \theta_0 = 0$. 

Likewise, we have that
\begin{align*}
\frac{1}{n} \sum_{i=1}^n \bigl\{ - \text{sgn}(Y_i + \Delta_{\gamma}) |Y_i + \Delta_{\gamma} |^\gamma \bigr\} \geq \frac{1}{2} \mathbb{E} \bigl[ - \text{sgn}(Y + \Delta_{\gamma} )|Y + \Delta_{\gamma} |^{\gamma - 1} \bigr] < 0.
\end{align*}

By the strong convexity of the function $\theta \mapsto \frac{1}{n}\sum_{i=1}^n |Y_i - \theta|^\gamma$ therefore, we have that $|\hat{\theta}_{\gamma} - \theta_0 | = |\hat{\theta}_{\gamma}| \leq \Delta_{\gamma}$. The first claim thus follows as desired.

To prove the second claim, we let $\Delta_{\gamma} = \tilde{\tau} \sqrt{ \frac{\gamma^\alpha}{n}}$ and follow exactly the same argument. The only difference is that the probability bound of Proposition~\ref{prop:tailbound1} in this case becomes, under our assumptions on $\tilde{\tau}$, 
\[
\exp \biggl\{ - \frac{a_1^2}{C_2 a_2}  \biggl( \frac{\tilde{\tau}^2}{ \sqrt{\frac{\gamma_u^\alpha}{n} \log\log \gamma_u} } \wedge \frac{\tilde{\tau}}{ \sqrt{ \frac{\gamma_u^\alpha}{n}} } \biggr\} \leq \exp \biggl\{ - \frac{a_1^2}{C_2 a_2} \tilde{\tau} \sqrt{ \frac{n}{\gamma_u^\alpha}}  \biggr\}.
\]

The entire theorem then follows. 
\end{proof}

\begin{proposition}
\label{prop:tailbound1}
Let $Z_1, \ldots, Z_n$ be independent and identically distributed random variables on $\mathbb{R}$ with a distribution $P$ symmetric around $0$ and write $\nu_{\gamma} := \mathbb{E} |Z|^{\gamma}$. Suppose there exists $\alpha \in (0, 2)$ and $a_1 \in (0, 1]$ and $a_2 \geq 1$ such that $\frac{a_1}{\gamma^\alpha} \leq \nu_{\gamma} \leq \frac{a_2}{\gamma^\alpha}$ for all $\gamma \geq 1$.

For $\gamma \geq 1$ and $x \in \mathbb{R}$, define $\psi_{\gamma}(x) := - \text{sgn}(x) |x|^{\gamma - 1}$. Let $\gamma_u > 2$ and let $\{\Delta_{\gamma}\}_{\gamma \in [2, \gamma_u]}$ be a collection of positive numbers; define the event
\begin{align*}
\mathcal{E} \equiv \mathcal{E}_{\gamma_u, \{\Delta_\gamma\}, \alpha}  := \left\{ \sup_{ \substack{\gamma \in [2, \gamma_u] \\ |t| = \Delta_{\gamma}} } 
\biggl| \frac{ \frac{1}{n} \sum_{i=1}^n \psi_\gamma(Z_i - t) - \mathbb{E} \psi_{\gamma}(Z-t) }
{ \frac{a_1}{2} \Delta_{\gamma} \gamma^{1 - \alpha}} \biggr| \geq \frac{1}{2}
\right\} 
\end{align*}
Let $C_1, C_2 > 0$ be universal constants and $C_{1\alpha}$ be a constant depending only on $\alpha$ (the values of these are specified in the proof). Then, the following holds:


\begin{enumerate}
\item[(i)] Suppose $\Delta_{\gamma} \gamma = \frac{a_1 c_0}{4}$ for some $c_0 \in (0, 1)$ and suppose that $\gamma_u \geq 2$ and $\sqrt{ \frac{\gamma_u^\alpha}{n}} \leq \frac{1}{4} \frac{c_0 a_1}{C_1 \sqrt{a_2}}$, then, we have that 
\begin{align}
\mathbb{P}(\mathcal{E}^c) \leq \exp \biggl\{ - \frac{c_0^2 a_1^4}{C_{1\alpha} a_2}  \biggl( \frac{n}{\gamma_u^\alpha} \biggr) \biggr\}, \label{eq:psi_tailbound1}
\end{align}
\item[(ii)] Let $\tilde{\tau} \geq \frac{C_1 \sqrt{a_2}}{a_1} \sqrt{\log\log n}$ and suppose $\Delta_{\gamma} \gamma = \tilde{\tau} \sqrt{ \frac{\gamma^{\alpha}}{n}}$. Suppose also $\gamma_u \geq e^{C_1 \frac{a_2}{a_1}} $ and $\sqrt{ \frac{\gamma_u^\alpha}{n}} \leq \frac{1}{4 \tilde{\tau}}$. Then, we have that
\[
\mathbb{P}(\mathcal{E}^c) \leq \exp \biggl\{ - \frac{a_1^2}{C_2 a_2}  \biggl( \frac{\tilde{\tau}^2}{ \sqrt{\frac{\gamma_u^\alpha}{n} \log\log \gamma_u} } \wedge \frac{\tilde{\tau}}{ \sqrt{ \frac{\gamma_u^\alpha}{n}} } \biggr\}.
\]
\end{enumerate}

\end{proposition}

\begin{proof}

We define the function class 
\begin{align}
\mathcal{F} \equiv \mathcal{F}_{\gamma_u, \{\Delta_{\gamma}\}, \alpha} := \biggl\{ \frac{\psi_{\gamma}(z - t) }{\frac{a_1}{2} \Delta_{\gamma} \gamma^{1-\alpha}} \,:\, t \in \{-\Delta_{\gamma}, \Delta_{\gamma}\}, \gamma \in [2, \gamma_u] \biggr\}. \label{eq:F_class1}
\end{align}

We now use Talagrand's inequality (Theorem~\ref{thm:talagrand}) to prove the Proposition. To this end, we derive upper bounds on various quantities involved in Talagrand's inequality.

\noindent \textbf{Step 1:} bounding $ \sup_{f \in \mathcal{F}} \| f(Z) \|_{\text{ess-sup}}$ and $\tilde{\sigma}^2 := \sup_{f \in \mathcal{F}} \mathbb{E} f(Z)^2$.

Using the fact $\Delta_{\gamma} \leq \frac{1}{4\gamma}$ in both cases, we observe that for any $\gamma \geq 2$, if $|t| = \Delta_{\gamma}$ and $|z| \leq 1$, then $| \psi_{\gamma}(z - t) | \leq (1 + \Delta_{\gamma})^{\gamma - 1} \leq e$. Therefore, we have that, 
\begin{align*}
U := \sup_{\gamma \in [2, \gamma_u]} \biggl| \frac{ \frac{1}{n} \sum_{i=1}^n \psi_\gamma(Z_i - t)  }
{ \frac{a_1}{2} \Delta_{\gamma} \gamma^{1 - \alpha}} \biggr|
&\leq \frac{2e}{a_1} \sup_{\gamma \in [2, \gamma_u]}  \frac{ \gamma^{\alpha} }{ \Delta_{\gamma} \gamma}.
\end{align*}

Thus, it follows that
\begin{align}
U \leq \begin{cases}
  \frac{C}{c_0 a^2_1} \gamma_u^{\alpha} & \text{ if $\Delta_{\gamma} \gamma = \frac{a_1 c_0}{4}$} \\
  \frac{C}{a_1} \frac{ \sqrt{\gamma_u^{\alpha} n } }{\tilde{\tau}} & \text{ if $\Delta_{\gamma} \gamma = \tilde{\tau} \sqrt{ \frac{\gamma^{\alpha}}{n} } $}
  \end{cases}. 
\end{align}

Next, we have that, writing $\tilde{\sigma}^2 := \sup_{f \in \mathcal{F}} \mathbb{E} f(Z)^2$, 
\begin{align*}
\tilde{\sigma}^2 &\leq \frac{1}{4 a_1^2} \sup_{ \substack{ \gamma \in [2, \gamma_u] \\ |t| = \Delta_\gamma} } \frac{ \gamma^{2\alpha} \mathbb{E} |Z - t|^{2(\gamma-1)} }{ \Delta_{\gamma}^2 \gamma^2} \\
&= \frac{1}{4 a_1^2} \sup_{  \gamma \in [2, \gamma_u] } \frac{ \gamma^{2\alpha} \mathbb{E} |Z - \Delta_{\gamma}|^{2(\gamma-1)} }{ \Delta_{\gamma}^2 \gamma^2} \\
&\leq \frac{C a_2 }{a_1^2} \sup_{  \gamma \in [2, \gamma_u] } \frac{\gamma^{\alpha}}{ \Delta_{\gamma}^2 \gamma^2},
\end{align*}
where the last inequality follows from Lemma~\ref{lem:tildetheta_moment}.

Therefore, we have that
\begin{align}
\tilde{\sigma}^2 &\leq
\begin{cases}
\frac{C a_2}{c_0^2 a_1^4} \gamma_u^{\alpha} & \text{ if $\Delta_{\gamma} \gamma = \frac{a_1 c_0}{4}$} \\
\frac{C a_2}{ a_1^2} \frac{n}{ \tilde{\tau}^2} & \text{ if $\Delta_{\gamma} \gamma = \tilde{\tau} \sqrt{ \frac{\gamma^{\alpha}}{n} }$ }.
\end{cases}
\end{align}

When $\Delta_{\gamma} \gamma = \tilde{\tau} \sqrt{ \frac{\gamma^{\alpha}}{n} }$, we also see that
\begin{align}
\tilde{\sigma}^2 \geq \frac{C}{a_1} \frac{1}{4 \Delta_2^2 } \geq \frac{C}{a_1} \frac{n}{\tilde{\tau}^2}. \label{eq:sigma_lower_bound}
\end{align}

\noindent \textbf{Step 2:} bounding the envelope function.\\

Define $F(z) := \sup_{f \in \mathcal{F}} |f(z)|$. Since, for any $z \in \mathbb{R}$,
\[
\sup_{\gamma \in [2, \gamma_u], |t|=\Delta_{\gamma}} |\psi_\gamma(z - t)| =  |(|z| + \Delta_{\gamma})^{\gamma-1}|,
\]
we have that 
\begin{align}
F(z) = \frac{4}{a_1} \sup_{\gamma \in [2, \gamma_u]} \frac{ \gamma^{\alpha} |(|z| + \Delta_{\gamma})^{\gamma-1} |}{\Delta_{\gamma} \gamma}. \label{eq:envelope}
\end{align}

Using the fact that the distribution of $Z$ is symmetric around 0, and defining $K := \lceil \log_2 \gamma_u \rceil$, 
\begin{align}
\mathbb{E} F^2(Z) &= \frac{16}{a^2_1} \int_0^1 \sup_{\gamma \in [2, \gamma_u]} \frac{ \gamma^{2\alpha} (z + \Delta_{\gamma})^{2(\gamma-1)} }{\Delta_{\gamma}^2 \gamma^{2} } \, dP(z) \nonumber \\
&= \frac{16}{a^2_1} \sum_{k=1}^K \int_0^1 \sup_{\gamma \in [2^k, 2^{k+1}]}
\frac{ \gamma^{2\alpha} (z + \Delta_{\gamma})^{2(\gamma-1)} }{\Delta_{\gamma}^2 \gamma^{2} } \, dP(z) 
\end{align}

\noindent Case 1: suppose $\Delta_{\gamma} \gamma = \frac{a_1 c_0}{4}$. In this case, we have that
\begin{align*}
\mathbb{E} F^2(Z) &= \frac{16}{c_0^2 a^4_1 } \sum_{k=1}^K \int_0^1 \sup_{\gamma \in [2^k, 2^{k+1}]} \gamma^{2\alpha} (z + \Delta_{\gamma})^{2(\gamma-1)} \, dP(z) \\
&\leq \frac{C}{c_0^2 a^4_1 } \sum_{k=1}^K 2^{2k\alpha} \int_0^1 \sup_{\gamma \in [2^k, 2^{k+1}]} (z + \Delta_{\gamma})^{2(\gamma-1)} \, dP(z) \\
&\leq \frac{C}{c_0^2 a^4_1 } \sum_{k=1}^K 2^{2k\alpha} \cdot a_2 2^{-k\alpha} \\
&\leq \frac{C_\alpha a_2}{c_0^2 a^4_1 }  \gamma_u^{\alpha},
\end{align*}
where the second inequality follows from the third claim of Lemma~\ref{lem:tildetheta_moment}.

\noindent Case 2: suppose $\Delta_{\gamma} \gamma = \tilde{\tau} \sqrt{ \frac{\gamma^\alpha}{n} }$. In this case, 

\begin{align}
\mathbb{E} F^2(Z) &= \frac{C}{a_1^2} \frac{n}{\tilde{\tau}^2} \sum_{k=1}^K \int_0^1 \sup_{\gamma \in [2^k, 2^{k+1}]} \gamma^{\alpha} (z + \Delta_{\gamma})^{2(\gamma-1)} \, dP(z) \nonumber \\
&\leq \frac{C}{a_1^2} \frac{n}{\tilde{\tau}^2} \sum_{k=1}^K 2^{k\alpha} \int_0^1 \sup_{\gamma \in [2^k, 2^{k+1}]} (z + \Delta_{\gamma})^{2(\gamma-1)} \, dP(z) \nonumber \\
&\leq \frac{C}{a_1^2} \frac{n}{\tilde{\tau}^2} \sum_{k=1}^K 2^{k\alpha} \cdot C a_2 2^{-k\alpha} \leq \frac{C a_2}{a_1^2} \frac{n}{\tilde{\tau}^2}\log \gamma_u , \label{eq:F_case2}
\end{align}
where, in the second inequality, we use Lemma~\ref{lem:tildetheta_moment} again. \\

\noindent \textbf{Step 3:} bounding the VC-dimension of $\mathcal{F}$.\\

We first note that the class of univariate functions $\mathcal{G} := \bigl\{ \frac{| \, \cdot \, |^{\gamma - 1}}{\Delta_{\gamma} \gamma} \,:\, \gamma \geq 2 \bigr\}$ has VC dimension at most 4. This holds because $\log \mathcal{G}$ consists of functions of the form
\[
(\gamma-1) \log | \cdot | + \log (\Delta_\gamma \gamma)
\]
and thus lies in a subspace of dimension 2. It then follows from Lemma 2.6.15 and 2.6.18 (viii) of \cite{van1996weak} that $\mathcal{G}$ has VC-dimension at most 4. 

It then follows from Lemma 2.6.18 (vi) that $\mathcal{F}$ has VC-dimension at most 8. \\

\noindent \textbf{Step 4:} bounding the expected supremum.\\

Let us define
\begin{align}
\tilde{S}_n := \sup_{ \substack{\gamma \in [2, \gamma_u] \\ |t| = \Delta_{\gamma}} } \biggl| \frac{1}{n} \sum_{i=1}^n \frac{ \psi_{\gamma}(Z_i - t) - \mathbb{E} \psi_{\gamma}(Z - t)}{\frac{a_1}{2} \Delta_{\gamma} \gamma^{1-\alpha}} \biggr|.
\end{align}

\noindent Case 1: suppose $\Delta_{\gamma} \gamma = \frac{a_1 c_0}{4}$. Then, using the second claim of Theorem~\ref{thm:exp_sup_local}, we have that
\begin{align}
\mathbb{E} \tilde{S}_n \leq \frac{C_\alpha \sqrt{a_2}}{c_0 a^2_1} \sqrt{\frac{\gamma_u^{\alpha} }{n}}.
\end{align}

\noindent Case 2: suppose now that $\Delta_{\gamma} \gamma = \tilde{\tau} \sqrt{ \frac{\gamma^{\alpha}}{n} }$. 

We first note that, by~\eqref{eq:sigma_lower_bound} and~\eqref{eq:F_case2}, 
\begin{align}
\frac{\tilde{\sigma}}{ \| F \|_{L_2(P)}} \geq \frac{ C a_1}{a_2} \frac{1}{\sqrt{\log \gamma_u}}.
\end{align}

Define the entropy integral $J(\delta)$ as~\eqref{eq:entropy_integral} and note that $\frac{1}{\delta} J(\delta)$ is decreasing for $\delta \in (0, 1]$. By Corollary~\ref{cor:VC} and our bound on the VC-dimension of $\mathcal{F}$, we have that
\begin{align*}
\frac{\| F \|_{L_2(P)}}{\tilde{\sigma}} J\biggl( \frac{\tilde{\sigma}}{\| F \|_{L_2(P)}} \biggr) \leq
\sqrt{ 1 \vee \log \biggl( \frac{a_2}{C a_1} \sqrt{\log \gamma_u} \biggr)} \leq 
\sqrt{\log \log \gamma_u + \log \bigl( \frac{a_2}{a_1} \bigr) + C}.
\end{align*}

Therefore, using our upper and lower bounds on $\tilde{\sigma}$, upper bound on $U$ and upper bound on $\| F \|_{L_2(P)}$, we have, by the first claim of Theorem~\ref{thm:exp_sup_local}, that
\begin{align*}
\mathbb{E} \tilde{S}_n &\leq C \frac{\tilde{\sigma}}{\sqrt{n}} \bigl( \sqrt{\log \log \gamma_u + \log \bigl( \frac{a_2}{a_1} \bigr) + C} \bigr)
\biggl( 1 + \frac{ U }{\sqrt{n} \tilde{\sigma}} \sqrt{\log \log \gamma_u + \log \bigl( \frac{a_2}{a_1} \bigr) + C} \biggr) \\
&\leq \frac{C \sqrt{a_2}}{a_1} \frac{1}{\tilde{\tau}} \sqrt{\log \log \gamma_u + \log \bigl( \frac{a_2}{a_1} \bigr) + C} \leq \frac{C \sqrt{a_2}}{a_1} \frac{1}{\tilde{\tau}} \sqrt{\log \log \gamma_u}.
\end{align*}
where, in the second inequality, we used the fact that $\frac{U}{\sqrt{n} \tilde{\sigma}} \leq C \sqrt{ \frac{\gamma_u^\alpha}{n} } \leq C$, and in the last inequality, we used the hypothesis that $\gamma_u \geq e^{C_1 \frac{a_2}{a_1}}$ (with $C_1$ as a sufficiently large universal constant). \\

\noindent \textbf{Step 5:} bounding the tail probability.\\

Using our assumption that $\frac{C_1 \sqrt{a_2}}{ c_0 a^2_1} \sqrt{ \frac{\gamma_u^{\alpha}}{n}} \leq \frac{1}{4}$ and $\tilde{\tau} \geq \frac{C_1 \sqrt{a_2}}{ a_1} \sqrt{ \log \log n}$ (with $C_1$ as a sufficiently large universal constant), we have that $\mathbb{E} \tilde{S}_n \leq \frac{1}{4}$ in both the case where $\Delta_{\gamma} \gamma = \frac{a_1 c_0}{4}$ and the case where $\Delta_{\gamma} \gamma = \tilde{\tau} \sqrt{ \frac{\gamma^{\alpha}}{n}}$. 

\noindent Case 1: when $\Delta_\gamma \gamma = \frac{a_1 c_0}{4}$, we have that, writing $t = 3/4$,
\begin{align*}
\mathbb{P} (\mathcal{E}^c) &\leq \mathbb{P}( \tilde{S}_n - \mathbb{E} \tilde{S}_n \geq \frac{3}{4} ) \\
&\leq \exp\biggl\{ - \frac{ n t^2}{ U \mathbb{E} \tilde{S}_n + \tilde{\sigma}^2} \wedge \frac{ n t}{ \frac{2}{3} U} \biggr\} \\
&\leq \exp \biggl\{ - \frac{ c_0^2 a_1^4 }{ C_{\alpha} a_2} \biggl( \biggl(\frac{n}{\gamma_u^\alpha} \biggr)^{3/2} \wedge \frac{n}{\gamma_u^\alpha} \biggr)  \biggr\}.
\end{align*}

Case 2: when $\Delta_{\gamma} \gamma = \tilde{\tau} \sqrt{ \frac{\gamma^{\alpha}}{n}}$, we have

\begin{align*}
\mathbb{P} (\mathcal{E}^c) &\leq \mathbb{P}( \tilde{S}_n - \mathbb{E} \tilde{S}_n \geq \frac{3}{4} ) \\
&\leq \exp\biggl\{ - \frac{ n t^2}{ U \mathbb{E} \tilde{S}_n + \tilde{\sigma}^2} \wedge \frac{ n t}{ \frac{2}{3} U} \biggr\} \\
&\leq \exp \biggl\{ - \frac{a_1^2}{C a_2} \biggl( \frac{ \tilde{\tau}^2 }{ \sqrt{ \frac{\gamma_u^{\alpha}}{n}} \sqrt{ \log \log \gamma_u} } \wedge \frac{\tilde{\tau}}{ \sqrt{ \frac{\gamma_u^{\alpha}}{n}} } \biggr) \biggr\}.
\end{align*}

\end{proof}

\begin{theorem}
\label{thm:vhat_bound}
Let $Z_1, \ldots, Z_n$ be independent and identically distributed random variables on $\mathbb{R}$ with a distribution $P$ symmetric around $0$ and write $\nu_{\gamma} := \mathbb{E} |Z|^{\gamma}$. Suppose there exists $\alpha \in (0, 2)$ and $a_1 \in (0, 1]$ and $a_2 \geq 1$ such that $\frac{a_1}{\gamma^\alpha} \leq \nu_{\gamma} \leq \frac{a_2}{\gamma^\alpha}$ for all $\gamma \geq 1$.

Let $c_0 \in (0, 2^{-8}]$ and let $C_{1\alpha}, C_{2\alpha} > 0$ be constants depending only on $\alpha$ defined in Theorem~\ref{thm:thetahat_gamma_bound} and Proposition~\ref{prop:tailbound2}. Define $\gamma_u^\alpha = \frac{1}{C_{1\alpha} \vee C_{2\alpha}} \bigl( \frac{c_0^2 a_1^6}{a^3_2} \alpha^2 \bigr) \frac{n}{\log n}$ and suppose $n$ is large enough so that $\gamma_u \geq 2$. Then, with probability at least $1 - 2 n^{-\frac{1}{\alpha}}$, there exists a constant $C_3 \leq 2$ such that
\[
 C_3 V(\gamma) \geq \hat{V}(\gamma) \geq \frac{1}{C_3} V(\gamma), \quad \text{ for all $\gamma \in [2, (\gamma_u + 1)/2 ]$.}
\]
Moreover, on the same event, there exists a universal constant $C_4 \geq 1$ such that
\[
C_4 \frac{a_2}{a_1^2} \gamma^{\alpha - 2} \geq \hat{V}(\gamma) \geq \frac{1}{C_4} \frac{a_1}{a_2^2} \gamma^{\alpha - 2}, \quad \text{ for all $\gamma \in [2, (\gamma_u + 1)/2 ]$.}
\]
\end{theorem}

We note that, in Theorem~\ref{thm:vhat_bound}, by choosing $c_0$ arbitrarily close to 0, we can have $C_3$ be arbitrarily close to 1. 

\begin{proof}

By Theorem~\ref{thm:thetahat_gamma_bound}, with probability at least $1 - n^{-\frac{1}{\alpha}}$, we have that, simultaneously for all $\gamma \in [2, \gamma_u]$,
\[
|\hat{\theta}_{\gamma} - \theta_0 | \leq \frac{a_1 c_0}{4 \gamma}. 
\]

On this event, we have that
\[
\hat{\nu}_{\gamma} := \inf_{\theta \in \mathbb{R}} \frac{1}{n} \sum_{i=1}^n |Y_i - \theta |^\gamma = \inf_{|\theta -\theta_0| \leq \frac{a_1 c_0}{4}} \frac{1}{n} \sum_{i=1}^n |Y_i - \theta |^\gamma.
\]

Then, by Proposition~\ref{prop:tailbound2}, with probability at least $1 - n^{- \frac{1}{\alpha}}$, simultaneously for all $\gamma \in [2, \gamma_u]$, 
\[
1 - 3 \sqrt{c_0} \leq \frac{\hat{\nu}_\gamma}{\nu_{\gamma}} \leq 1 + 3 \sqrt{c_0}.
\]

Therefore, 
\begin{align*}
\hat{V}(\gamma) = \frac{ \hat{\nu}_{2(\gamma-1)} }{(\gamma-1)^2 \hat{\nu}_{\gamma-2}^2} 
\geq 
\frac{1 - 3 \sqrt{c_0}}{(1 + 3 \sqrt{c_0})^2} \frac{ \nu_{2(\gamma-1)} }{(\gamma-1)^2 \nu_{\gamma-2}^2} = \frac{1 - 3 \sqrt{c_0}}{(1 + 3 \sqrt{c_0})^2} V(\gamma).
\end{align*}

Likewise, we have that $\hat{V}(\gamma) \leq \frac{1 + 3 \sqrt{c_0}}{(1 - 3 \sqrt{c_0})^2} V(\gamma)$. Using our assumption that $c_0 \leq 2^{-8}$, the first claim of the theorem directly follows.

The second claim of the theorem follows then from Lemma~\ref{lem:Vgamma_bound}.

\end{proof}

\begin{proposition}
\label{prop:tailbound2}
Let $Z_1, \ldots, Z_n$ be independent and identically distributed random variables on $\mathbb{R}$ with a distribution $P$ symmetric around $0$ and write $\nu_{\gamma} := \mathbb{E} |Z|^{\gamma}$. Suppose there exists $\alpha \in (0, 2)$ and $a_1 \in (0, 1]$ and $a_2 \geq 1$ such that $\frac{a_1}{\gamma^\alpha} \leq \nu_{\gamma} \leq \frac{a_2}{\gamma^\alpha}$ for all $\gamma \geq 1$.

Let $\gamma_u \geq 2$ and $c_0 \in (0, 1)$. Define the event 
\begin{align}
\mathcal{A}_{\gamma_u, c_0} := \biggl\{ \sup_{\gamma \in [2, \gamma_u]} \biggl| \frac{ \inf_{|\theta - \theta_0| \leq \frac{a_1 c_0}{4 \gamma}} \frac{1}{n}\sum_{i=1}^n |Y_i - \theta|^{\gamma} - \nu_{\gamma}}{ \nu_\gamma } \biggr| \leq 3 \sqrt{c_0} 
\biggr\}.
\end{align}

Let $C_{2\alpha} > 0$ be a constant depending only on $\alpha$ (its value is specified in the proof). Suppose $\frac{\gamma_u^\alpha}{n} \leq c_0^2 \frac{a_1^2}{C_{2\alpha} a_2}$. Then, we have that
\[
\mathbb{P}( \mathcal{A}_{\gamma_u, c_0}^c) \leq \exp \biggl\{ - \frac{a_1^2}{C_{2\alpha} a_2} c_0^2 \biggl( \frac{ n }{\gamma_u^{\alpha} } \biggr) \biggr\}.
\]

\end{proposition}

\begin{proof}

First, we claim that, for all $\gamma \geq 1$, $z \in \mathbb{R}$, $t \geq 0$, and $\kappa > 0$, it holds that

\begin{align}
| z - t |^{\gamma} &\geq (|z| - t)_+^{\gamma} 
 = |z|^\gamma \biggl( 1 - \frac{t}{|z|} \biggr)_+^{\gamma} \nonumber \\
 &= |z|^\gamma \biggl( 1 - \frac{t}{\kappa} \frac{\kappa}{|z|} \biggr)_+^\gamma \geq |z|^\gamma \biggl( 1 - \frac{t}{\kappa} \biggr)_+^\gamma - \kappa^\gamma \nonumber \\
 &\geq |z|^\gamma \biggl( 1 - \gamma \frac{t}{\kappa} \biggr) - \kappa^\gamma.
 \label{eq:cz_ineq}
\end{align}

Now define $\Delta_\gamma = \frac{a_1 c_0}{4\gamma}$ and $\tilde{\theta} := \argmin_{ |\theta - \theta_0| \leq \Delta_{\gamma}} \frac{1}{n} \sum_{i=1}^n | Y_i - \theta |^\gamma$ and observe that $ \tilde{t} := \tilde{\theta} - \theta_0 = \argmin_{ |t| \leq \Delta_{\gamma}} \frac{1}{n} \sum_{i=1}^n | Z_i - t |^\gamma$. Suppose without loss of generality that $\tilde{t} \geq 0$. Then, using~\eqref{eq:cz_ineq}, we have that, for any $\kappa > 0$,
\begin{align*}
\inf_{ |\theta - \theta_0| \leq \Delta_\gamma} \frac{1}{n} \sum_{i=1}^n |Y_i - \theta |^\gamma 
&= \frac{1}{n} \sum_{i=1}^n |Z_i - \tilde{t}|^\gamma  \\
&\geq  \biggl( 1 - \frac{\gamma \Delta_{\gamma}}{\kappa} \biggr) \biggl( \frac{1}{n} \sum_{i=1}^n |Z_i|^\gamma \biggr) - \kappa^\gamma .
\end{align*}

We also trivially have that $\inf_{ |\theta - \theta_0| \leq \Delta_\gamma} \frac{1}{n} \sum_{i=1}^n |Y_i - \theta |^\gamma \leq \frac{1}{n} \sum_{i=1}^n |Z_i|^\gamma$. Therefore, writing $\mathbb{E}_n |z|^\gamma := \frac{1}{n} \sum_{i=1}^n |Z_i|^\gamma$ and $\mathbb{E}_n |y-\theta|^\gamma := \frac{1}{n} \sum_{i=1}^n |Y_i - \theta|^\gamma$, we have that, for any $\kappa > 0$,

\begin{align}
\frac{ \mathbb{E}_n |z|^\gamma - \nu_{\gamma}}{  \nu_\gamma} 
\geq \frac{ \min_{|\theta -\theta_0| \leq \Delta_{\gamma}} \mathbb{E}_n |y - \theta|^\gamma - \nu_\gamma}{  \nu_\gamma}
\geq \frac{ \bigl( 1 - \frac{\gamma \Delta_\gamma}{\kappa} \bigr) \mathbb{E}_n |z|^\gamma - \nu_\gamma - \kappa^\gamma}{ \nu_\gamma}.
\end{align}

Therefore, we have that
\begin{align}
\biggl| \frac{ \min_{|\theta -\theta_0| \leq \Delta_{\gamma}} \mathbb{E}_n |y - \theta|^\gamma - \nu_\gamma}{ \nu_\gamma} \biggr| 
&\leq \underbrace{ \biggl| \frac{\mathbb{E}_n |z|^{\gamma} - \nu_\gamma}{ \nu_\gamma} \biggr| }_{\text{Term 1}}
+ \underbrace{ \inf_{\kappa > 0} \biggr( \frac{\gamma \Delta_\gamma}{ \kappa} + \frac{\kappa^\gamma}{ \nu_\gamma} \biggr) }_{\text{Term 2}}. \label{eq:variation_decomp}
\end{align}

\noindent \textbf{Bounding Term 2:} \\

Since $\Delta_{\gamma} \gamma = \frac{a_1 c_0}{4}$, 
by setting $\kappa = \bigl( \frac{a^2_1 c_0}{4 \gamma^{\alpha}} \bigr)^{\frac{1}{\gamma+1}}$, we have
\begin{align*}
\text{Term 2} &\leq \inf_{\kappa > 0} \biggl( \frac{a_1 c_0}{4 \kappa} + \frac{\kappa^\gamma}{a_1 \gamma^\alpha} \biggr) \\
&\leq \frac{a_1 c_0}{4} \biggl( \frac{4}{a^2_1 c_0} \biggr)^{\frac{1}{\gamma+1}} \gamma^{ \frac{\alpha}{\gamma+1}} 
\leq 2 \sqrt{c_0}.
\end{align*}

\noindent \textbf{Bounding Term 1:} \\

To bound Term 1, we define the function class
\[
\mathcal{F}_{\gamma_u} := \biggl\{ z \mapsto \frac{|z|^\gamma}{\nu_\gamma} \,:\, \gamma \in [2, \gamma_u] \biggr\},
\]
so that we have $\sup_{\gamma \in [2, \gamma_u]} \bigl| \frac{ \mathbb{E}_n |z|^\gamma - \nu_\gamma}{\nu_\gamma} \bigr| = \sup_{f \in \mathcal{F}_{\gamma_u}} \bigl|  \mathbb{E}_n f(z) - \mathbb{E}f(Z) \bigr|$.

We observe that 
\begin{align*}
\tilde{\sigma}^2 &:= \sup_{\gamma \in [2, \gamma_u]} \mathbb{E} \frac{|Z|^\gamma}{ \nu_\gamma} \leq \frac{a_2}{a_1} \\
U &:= \sup_{\gamma \in [2, \gamma_u]} \frac{\| Z \|_{\text{ess-inf}}^\gamma}{\nu_\gamma} \leq \frac{1}{a_1} \gamma_u^\alpha. 
\end{align*}

Moreover, defining $F(z) := \sup_{\gamma \in [2, \gamma_u]} \frac{|z|^\gamma}{\nu_\gamma}$ and $K = \lceil \log \gamma_u \rceil$, we have that
\begin{align}
\mathbb{E} F^2(Z) &\leq \mathbb{E} \bigg( \sup_{\gamma \in [2, \gamma_u]} \frac{|Z|^{2\gamma}}{\nu^2_\gamma} \biggr) \\
&\leq \frac{1}{a^2_1} \sum_{k=1}^K \int_0^1 \sup_{\gamma \in [2^k, 2^{k+1}]} \gamma^{2\alpha} |z|^\gamma \, dP(z) \\
&\leq \frac{1}{a^2_1} \sum_{k=1}^K 2^{2k\alpha + 1} \nu_{2^k} 
\leq \frac{a_2}{a^2_1} \sum_{k=1}^K 2^{k\alpha + 1} \leq 
C_{\alpha} \frac{a_2}{a_1^2} \gamma_u^\alpha.
\end{align}

We note that $\log \mathcal{F}_{\gamma_u}$ is a subset of a linear subspace of dimension 2 (see Step 3 in the proof of Proposition~\ref{prop:tailbound1}). By Lemma 2.6.15 and 2.6.18 (viii) of \cite{van1996weak}, we know that the VC dimension of $\mathcal{F}_{\gamma_u}$ is at most 4. 

Write $\tilde{S}_n = \sup_{\gamma \in [2, \gamma_u]}  \bigl| \frac{1}{n} \sum_{i=1}^n |Z_i|^{\gamma} - \nu_{\gamma} \bigr\}$. Then, by Corollary~\ref{cor:VC} and the second claim of Theorem~\ref{thm:exp_sup_local}, we have that
\begin{align*}
\mathbb{E} \tilde{S}_n \leq  \frac{C_{\alpha} \sqrt{a_2}}{a_1} \sqrt{ \frac{\gamma_u^{\alpha}}{n} }.
\end{align*}

Therefore, using our hypothesis that $\frac{\gamma_u^\alpha}{n} \leq c_0^2 \frac{a_1^2}{C_{2\alpha} a_2}$ where $C_{2\alpha}$ is chosen to be sufficiently large, then $\mathbb{E} \tilde{S}_n \leq \frac{1}{2} a_0$. Then,

\begin{align*}
\mathbb{P} \biggl\{ \sup_{\gamma \in [2, \gamma_u]} \biggl| \frac{ \frac{1}{n} \sum_{i=1}^n |Z_i|^\gamma - \nu_\gamma }{\nu_\gamma} \biggr| \geq c_0 \biggr\} 
& \leq \mathbb{P}\biggl( \tilde{S}_n - \mathbb{E} \tilde{S}_n \geq \frac{c_0}{2} \biggr) \\
&\leq \exp \biggl\{ - \frac{a_1^2}{C_\alpha a_2} \biggl( c_0^2 \biggl( \frac{n}{ \gamma_u^\alpha} \biggr)^{3/2} \wedge \frac{c_0 n }{\gamma_u^{\alpha} } \biggr) \biggr\}.
\end{align*}

Therefore, by~\eqref{eq:variation_decomp}, it holds that
\begin{align*}
&\mathbb{P}\biggl\{ \sup_{\gamma \in [2, \gamma_u]} \biggl| \frac{ \inf_{|\theta - \theta_0| \leq \Delta_{\gamma}} \frac{1}{n} \sum_{i=1}^n |Y_i - \theta|^\gamma - \nu_\gamma }{\nu_\gamma} \biggr| \geq 3 \sqrt{c_0} \biggr\} \\
&\leq \exp \biggl\{ - \frac{a_1^2}{C_\alpha a_2} \biggl( c_0^2 \biggl( \frac{n}{ \gamma_u^\alpha} \biggr)^{3/2} \wedge c_0 \frac{ n }{\gamma_u^{\alpha} } \biggr) \biggr\}.
\end{align*}

By inflating the value of $C_{2\alpha}$ if necessary, the Proposition follows as desired. 
\end{proof}

\begin{lemma}
\label{lem:moment_tail}
Let $X$ be a random variable on $[-1, 1]$ with a distribution $P$ symmetric around $0$. If there exists $a_1 > 0$ and $\alpha \geq 0$ such that $\mathbb{E}|X|^\gamma \geq \frac{a_1}{\gamma^\alpha}$ for all $\gamma \geq 1$, then we have that
\[
\mathbb{P}\biggl( X \geq 1 - 2^{1 + \frac{2}{\alpha}} a_1^{- \frac{1}{\alpha}} \frac{1}{\alpha^{\frac{1}{\alpha}+1}} \frac{ \log^{1+ \frac{1}{\alpha}} n}{n^{\frac{1}{\alpha}}} \biggr) \geq \alpha^{-1} \frac{\log n}{n}.
\]
\end{lemma}

\begin{proof}

As a short-hand, write $\delta := 2^{1 + 2/\alpha} a_1^{-1/\alpha} \frac{1}{\alpha^{\frac{1}{\alpha}+1}} \frac{ \log^{1+1/\alpha} n}{n^{1/\alpha}}$ and $p =  \bigl( \alpha \frac{a_1 n }{4 \log n} \bigr)^{1/\alpha}$; note that $\delta p = 2 \alpha^{-1} \log n$. Then,
\begin{align*}
\mathbb{P}(X \geq 1 - \delta) 
&= \int_{1-\delta}^1 dP(x) 
\geq \int_{1-\delta}^1 x^p dP(x)  \\
&= \frac{1}{2} \mathbb{E}|X|^p - \int_0^{1-\delta} x^p dP(x) \\
&\geq \frac{a_1}{2 p^{\alpha}} - (1 - \delta)^p 
\geq \frac{a_1}{2 p^{\alpha}} - e^{- \delta p}  \\
&= 2 \frac{1}{\alpha} \frac{\log n}{n} - \frac{1}{\alpha} \frac{1}{n^2} \geq  \frac{1}{\alpha} \frac{\log n}{n}.
\end{align*}

\end{proof}

\begin{corollary}
\label{cor:midrange}
Let $X_1, \ldots, X_n$ be independent and identically distributed random variables on $[-1, 1]$ with a distribution $P$ symmetric around $0$ and let $X_{\text{mid}} = \frac{X_{(n)} + X_{(1)}}{2}$. If there exists $a_1 > 0$ and $\alpha \geq 0$ such that $\mathbb{E}|X|^\gamma \geq \frac{a_1}{\gamma^\alpha}$ for all $\gamma \in \mathbb{N}$, then we have that
\begin{align*}
\mathbb{P}\biggl( |X_{\text{mid}}| \leq 
2^{2 + \frac{2}{\alpha}} a_1^{- \frac{1}{\alpha}} \frac{1}{\alpha^{\frac{1}{\alpha}+1}} \frac{ \log^{1+ \frac{1}{\alpha}} n}{n^{\frac{1}{\alpha}}} \biggr) \geq 1 - \frac{2}{n^{1/\alpha}}.
\end{align*}

\end{corollary}

\begin{proof}

As a short-hand, write $\delta = 2^{1 + \frac{2}{\alpha}} a_1^{- \frac{1}{\alpha}} \frac{1}{\alpha^{\frac{1}{\alpha}+1}} \frac{ \log^{1+ \frac{1}{\alpha}} n}{n^{\frac{1}{\alpha}}}$. 
By the fact that $P$ is symmetric around $0$ and Lemma~\ref{lem:moment_tail}, we have
\begin{align*}
\mathbb{P}( |X_{\text{mid}}| \geq 2\delta) &\leq \mathbb{P}( X_{(n)} \leq 1 - \delta \text{ or } X_{(1)} \geq -1 + \delta ) \\
&\leq 2 \mathbb{P}( X_{(n)} \leq 1 - \delta) \\
&\leq 2 \bigl\{ \mathbb{P}( X_1 \leq 1 - \delta) \bigr\}^n \\
&\leq 2 \biggl(1 - \frac{1}{\alpha} \frac{\log n}{n} \biggr)^n  \leq 2 e^{- \frac{1}{\alpha} \log n} \leq \frac{2}{n^{1/\alpha}}. 
\end{align*}
The desired conclusion thus follows. 
\end{proof}

\subsection{Proof of Examples}
\label{sec:proof_of_examples}

\begin{proof} (of Proposition~\ref{prop:moment_example}) \\
It suffices to show that there exists constants $C''_{\alpha,1}, C''_{\alpha,2} > 0$ such that
\[
\frac{C''_{\alpha,1}}{\gamma^\alpha} \leq \int_{-1}^1 |x|^\gamma (1 - |x|)^{\alpha-1} dx \leq \frac{C''_{\alpha,2}}{\gamma^\alpha}.
\]
Indeed, we have by Stirling's approximation that
\begin{align*}
\int_{-1}^1 |x|^\gamma (1 - |x|)^{\alpha-1} dx 
&= 2 \int_0^1 x^\gamma (1 - x)^{\alpha-1} dx \\
&= 2 \frac{\Gamma(\alpha)\Gamma(\gamma+1)}{\Gamma(\gamma+\alpha+1)} \\
&\asymp 2 \Gamma(\alpha) \biggl( \frac{\gamma}{e} \biggr)^{\gamma} \biggl( \frac{\gamma+\alpha}{e} \biggr)^{-(\gamma+\alpha)} \biggl( \frac{\gamma}{\gamma+\alpha}\biggr)^{1/2} \\
&= 2 \Gamma(\alpha) \biggl( \frac{\gamma}{\gamma+\alpha} \biggr)^{\gamma + 1/2} e^{\alpha} \biggl( \frac{\gamma}{\gamma+\alpha} \biggr)^{\alpha} \gamma^{-\alpha}. 
\end{align*}

The conclusion of the Proposition then directly follows from the fact that $1 \geq ( \frac{\gamma}{\gamma+\alpha})^{\gamma+1/2} \geq e^{-3}$ for all $\gamma \geq 2$. 

\end{proof}

For a given density $p(\cdot)$, we define
\[
H^2(\theta_1, \theta_2) := \frac{1}{2} \int_{\mathbb{R}} \bigl( p(x - \theta_1)^{1/2} - p(x - \theta_2)^{1/2} \bigr)^2 \, dx
\]
for any $\theta_1, \theta_2 \in \mathbb{R}$. 
\begin{proposition}
\label{prop:hellinger_bound}
Let $\alpha \in (0, 2)$ and suppose $X$ is a random variable with density $p(\cdot)$ satisfying
\[
C_{\alpha,1} (1 - |x|)_+^{\alpha-1} \leq p(x) \leq C_{\alpha,2} (1 - |x|)_+^{\alpha-1} 
\]
for $C_{\alpha,1}, C_{\alpha,2} > 0$ dependent only on $\alpha$. Suppose also that $\bigl| \frac{p'(x)}{p(x)} \bigr| \leq \frac{C}{1 - |x|}$ for some $C > 0$. 

Suppose $p(\cdot)$ is symmetric around 0. Then, there exist $C'_{\alpha, 1}, C'_{\alpha,2}$ dependent only on $\alpha$ and $C$ such that
\[
C'_{\alpha,1} |\theta_1 - \theta_2|^\alpha \leq H^2(\theta_1, \theta_2) \leq C'_{\alpha,2} | \theta_1 - \theta_2 |^\alpha 
\]
for all $\theta_1, \theta_2 \in \mathbb{R}$.
\end{proposition}

\begin{proof}
Since $H^2(\theta_1, \theta_2) = H^2(0, \theta_1 - \theta_2)$, it suffices to bound $H^2(0, \theta)$ for $\theta \geq 0$. 

For the lower bound, we observe that
\begin{align*}
H^2(0, \theta) &= \int_{-1}^{1+\theta} 
  \bigl\{ p(x)^{1/2} - p(x - \theta)^{1/2} \bigr\}^2 dx \\
  &\geq \int_1^{1+\theta} p(x - \theta) dx 
  = \int_{1-\theta}^1 f(t) dt \\
  &\geq C_{\alpha,1} \int_{1-\theta}^1 (1 - t)^{\alpha-1} dt \\
  &= C_{\alpha,1} \biggl[ - \frac{ (1-t)^\alpha}{\alpha} \biggr]_{1-\theta}^1 
  = \frac{C_{\alpha,1}}{\alpha} \theta^\alpha.
\end{align*}

To establish the upper bound, observe that, by symmetry of $p(\cdot)$, 
\begin{align}
H^2(0, \theta) &= 2\int_{\theta/2}^{1+\theta} \bigl\{ p(x)^{1/2} - p(x - \theta)^{1/2} \bigr\}^2 dx \nonumber \\
&= 2 \int_{1-\theta}^{1+\theta} \bigl\{ p(x)^{1/2} - p(x - \theta)^{1/2} \bigr\}^2 dx + 2 \int_{\theta/2}^{1-\theta} \bigl\{ p(x)^{1/2} - p(x - \theta)^{1/2} \bigr\}^2 dx. \label{eq:H_decomp}
\end{align}

We upper bound the two terms of~\eqref{eq:H_decomp} separately. To bound the first term,
\begin{align*}
& \int_{1-\theta}^{1+\theta} \bigl\{ p(x)^{1/2} - p(x - \theta)^{1/2} \bigr\}^2 dx \\
&\leq \int_1^{1+\theta} p(x-\theta) dx + \int_{1-\theta}^{1} p(x) \vee p(x - \theta) dx \\
&\leq \frac{C_{\alpha,2}}{\alpha} \theta^\alpha +  C_{\alpha,2} \int_{1-\theta}^{1}
    \bigl\{ (1 - x)^{\alpha-1} \vee (1 - (x - \theta))^{\alpha-1} \bigr\} dx 
\end{align*}

If $\alpha \geq 1$, then $(1 - x)^{\alpha-1} \vee (1 - (x - \theta))^{\alpha-1} = (1 - (x - \theta))^{\alpha-1}$ and 
\[
\int_{1-\theta}^1 (1 - (x - \theta))^{\alpha-1} dx = \int_{1-2\theta}^{1-\theta} (1-x)^{\alpha-1} dx = (2^\alpha - 1) \frac{\theta^\alpha}{\alpha} \geq \frac{\theta^\alpha}{\alpha}.
\]
On the other hand, if $\alpha < 1$, then $(1 - x)^{\alpha-1} \vee (1 - (x - \theta))^{\alpha-1} = (1 - x)^{\alpha-1}$ and $\int_{1-\theta}^1 (1 - x)^{\alpha-1} dx = \frac{\theta^\alpha}{\alpha}$. Hence, we have that
\begin{align*}
\int_{1-\theta}^{1+\theta} \bigl\{ p(x)^{1/2} - p(x - \theta)^{1/2} \bigr\}^2 dx &\leq \frac{2 C_{\alpha,2}}{\alpha} \theta^\alpha. 
\end{align*}

We now turn to the second term of~\eqref{eq:H_decomp}. Write $\phi(x) = \log p(x)$ and note that $\phi'(x) = \frac{p'(x)}{p(x)}$. Then, by mean value theorem, there exists $\theta_x \in (0, \theta)$ depending on $x$ such that
\begin{align*}
2 \int_{\theta/2}^{1-\theta} \bigl\{ p(x)^{1/2} - p(x - \theta)^{1/2} \bigr\}^2 dx 
&= 2 \int_{\theta/2}^{1-\theta} \frac{\theta^2}{4} \phi'(x - \theta_x)^2 e^{\phi(x - \theta_x)} dx \\
&\leq \frac{C C_{\alpha,2}}{2} \theta^2  \int_{\theta/2}^{1-\theta} \biggl(\frac{1}{1 - |x - \theta_x|} \biggr)^2 (1 - |x- \theta_x|)^{\alpha-1} dx\\
&= \frac{C C_{\alpha,2}}{2} \theta^2  \int_{\theta/2}^{1-\theta} (1 - |x- \theta_x|)^{\alpha-3} dx \\
&\leq \frac{C C_{\alpha,2}}{2} \theta^2  \int_{0}^{1-\theta} (1 - x)^{\alpha-3} dx  \\
&= \frac{C C_{\alpha,2}}{2} \theta^2  \frac{\theta^{\alpha-2}}{2-\alpha} = \frac{C C_{\alpha,2}}{2(2-\alpha)} \theta^\alpha,
\end{align*}
where the second inequality follows because $\alpha - 3 < 0$. The desired conclusion immediately follows. 
\end{proof}

\begin{remark}
\label{rem:hellinger_bound}
We observe that if a density $p$ is of the form 
\[
p(x) = C_\alpha (1 - |x|)^{\alpha - 1} \mathbbm{1}\{ |x| \leq 1\},
\]
for a normalization constant $C_{\alpha} > 0$, 
then $\frac{p'(x)}{p(x)} \lesssim \frac{1}{1 + |x|}$ as required in Proposition~\ref{prop:hellinger_bound}. Therefore, we immediately see that for such a density, it holds that $H^2(\theta_1, \theta_2) \propto_{\alpha} |\theta_1 - \theta_2 |^{\alpha}$. 
\end{remark}

\subsection{Proof of Proposition~\ref{prop:pop_mle}} 
\label{sec:pop_mle_proof}

\begin{proof}

We first note that if $Y = \theta_0 + Z$ where $Z$ has a density $p(\cdot)$ symmetric around $0$, then, for $\gamma > 0$, 
\[
L(\gamma) = \frac{1}{\gamma} \log \bigl( \mathbb{E}|Z|^\gamma \bigr) + \frac{1 + \log \gamma}{\gamma} + \log \Gamma\bigl( 1 + \frac{1}{\gamma} \bigr).
\]

To prove the first claim, suppose that $Z$ is supported on all of $\mathbb{R}$. We observe that
\[
\lim_{\gamma \rightarrow \infty} L(\gamma) = \log \biggl( \lim_{\gamma \rightarrow \infty} \bigl\{ \mathbb{E}|Z|^{\gamma} \bigr\}^{\frac{1}{\gamma}} \biggr).
\]
We thus need only show that $\lim_{\gamma \rightarrow \infty} \bigl\{ \mathbb{E}|Z|^{\gamma} \bigr\}^{\frac{1}{\gamma}} = \infty$. Let $M > 0$ be arbitrary, then, for any $\gamma > 0$,
\begin{align*}
\{ \mathbb{E}|Z|^\gamma \}^{\frac{1}{\gamma}} 
&\geq
\{ \mathbb{E}\bigl[ |Z|^\gamma \mathbbm{1}\{ |Z| \geq M \} \bigr] \}^{\frac{1}{\gamma}} \\
&\geq M \cdot \mathbb{P}( |Z| \geq M)^{\frac{1}{\gamma}}.
\end{align*}
Since $\mathbb{P}(|Z| \geq M) > 0$ for all $M > 0$ by assumption, we see that $\lim_{\gamma \rightarrow \infty} \{ \mathbb{E}|Z|^\gamma \}^{\frac{1}{\gamma}} \geq M$. Since $M$ is arbitrary, the claim follows.

Now consider the second claim of the Proposition and assume that $\| Z \|_\infty = 1$; write $g(\cdot)$ as the density of $|Z|$. Writing $\eta = \frac{1}{\gamma}$, we have that
\begin{align*}
L(1/\eta) = \eta \log\bigl( \mathbb{E}|Z|^{\frac{1}{\eta}} \bigr) + \eta(1 - \log \eta) + \log \Gamma 
 (1 + \eta).
\end{align*}
Differentiating with respect to $\eta$, we have
\begin{align*}
\frac{ d L(1/\eta)}{d \eta} 
&= \log \frac{\mathbb{E}|Z|^{\frac{1}{\eta}}}{\eta} - \frac{ \mathbb{E}\{ |Z|^{\frac{1}{\eta}} \log|Z| \} }{\eta \mathbb{E}|Z|^{\frac{1}{\eta}} } + \frac{\Gamma'(1+\eta)}{\Gamma(1+\eta)} \\
&= \log \frac{\int_0^1 u^{\frac{1}{\eta}} g(u) d u }{\eta} - \frac{ \int_0^1 u^{\frac{1}{\eta}} \log(u) g(u) d u }{ \eta \int_0^1 u^{\frac{1}{\eta}} g(u) d u } + \frac{\Gamma'(1+\eta)}{\Gamma(1+\eta)}.
\end{align*}
We make a change of variable by letting $t = - \frac{1}{\eta} \log u$ to obtain
\begin{align*}
\frac{d L(1/\eta)}{d \eta} &= \log \frac{ \int_0^\infty e^{-t} g(e^{-\eta t}) e^{-\eta t} \eta d t }{\eta} - 
\frac{ \int_0^\infty e^{-t} (- \eta t) g(e^{-\eta t}) e^{- \eta t} \eta d t }{ \eta \int_0^\infty e^{-t} g(e^{ - \eta t} ) e^{-\eta t} \eta d t } + \frac{\Gamma'(1+\eta)}{\Gamma(1+\eta)} \\
&= \log \biggl\{ \int_0^\infty e^{-t} g(e^{-\eta t}) e^{-\eta t} d t \biggr\} + \frac{ \int_0^\infty t e^{-t} g(e^{- \eta t}) e^{- \eta t} d t }{ \int_0^\infty e^{-t} g(e^{-\eta t}) e^{- \eta t} dt } + \frac{\Gamma'(1+\eta)}{\Gamma(1+\eta)}.
\end{align*}

Therefore, using the fact that $\lim_{\eta \rightarrow 0} \frac{\Gamma'(1+\eta)}{\Gamma(1+\eta)} = - \gamma_{\text{E}}$, we have that
\begin{align*}
\lim_{\eta \rightarrow 0} \frac{d L(1/\eta)}{d \eta} &= \log \biggl\{ g(1) \int_0^\infty e^{-t} dt \biggr\} + \frac{\int_0^\infty t e^{-t} dt }{\int_0^\infty e^{-t} dt } - \gamma_{\text{E}} \\
&= \log g(1) + 1 - \gamma_{\text{E}}.
\end{align*}

Therefore, if $g(1) > e^{\gamma_{\text{E}} - 1}$, then $\lim_{\eta \rightarrow 0} \frac{d L(1/\eta)}{d \eta} > 0$ and hence, $\eta = 0$ is a local minimum of $L(1/\eta)$. On the other hand, if $g(1) < e^{\gamma_{\text{E}} - 1}$, then $\lim_{\eta \rightarrow 0} \frac{d L(1/\eta)}{d \eta} < 0$ and $\eta = 0$ is not a local minimum. The Proposition follows as desired. 

\end{proof}

\section{Other material}

\subsection{Technical Lemmas}

\begin{lemma}
\label{lem:Vgamma_bound}
Let $Z$ be a random variable supported on $[-1, 1]$. For $\gamma \geq 1$, define $\nu_{\gamma} := \mathbb{E} |Z|^{\gamma}$ and suppose there exists $\alpha \in (0, 2]$, $a_1 \in (0, 1]$, and $a_2 \in [1, \infty)$ such that $a_1 \gamma^{-\alpha} \leq \nu_{\gamma} \leq a_2 \gamma^{-\alpha}$ for all $\gamma \geq 1$. 

Define $V(\gamma) := \frac{ \mathbb{E} |Z|^{2(\gamma-1)}}{ (\gamma - 1)^2 \{ \mathbb{E}|Z|^{\gamma - 2} \}^2 }$. Then, for some universal constant $C \geq 1$, for all $\gamma \geq 2$,
\[
C \frac{a_2}{a_1^2} \gamma^{\alpha - 2} \geq V(\gamma) \geq \frac{1}{C} \frac{a_1}{a_2^2} \gamma^{\alpha - 2} .
\]
\end{lemma}

\begin{proof}
First suppose $\gamma \in [2, 3]$. Then we have that
\begin{align*}
a_1 \leq \mathbb{E}|Z| \leq \mathbb{E} |Z|^{\gamma - 2} \leq \{\mathbb{E} |Z| \}^{\gamma - 2} \leq a_2,
\end{align*}
where the second inequality follows because $|Z| \leq 1$, the third inequality follows from Jensen's inequality. Therefore, we have that
\begin{align*}
V(\gamma) \geq \frac{ a_1 \{ 2( \gamma - 1) \}^{-\alpha} }{(\gamma - 1)^2 a_2^2 } 
\geq \frac{1}{4^{\alpha+1} 3^{\alpha - 2} } \frac{a_1}{a_2^2} \gamma^{\alpha - 2}.
\end{align*}
The upper bound on $V(\gamma)$ follows similarly.

Now suppose $\gamma \geq 3$, then,
\begin{align*}
V(\gamma) &= \frac{ \nu_{2(\gamma-1)} }{(\gamma-1)^2 \nu_{\gamma - 2}^2 } \geq \frac{a_1 \{2(\gamma-1)\}^{-\alpha} }{(\gamma-1)^2 a_2^2 (\gamma - 2)^{-2\alpha}} \\
&= \frac{a_1}{a_2^2} 2^{-\alpha} \biggl( \frac{\gamma - 2}{\gamma-1} \biggr)^\alpha \biggl( \frac{\gamma - 2}{\gamma} \biggr)^\alpha \biggl( \frac{\gamma}{\gamma-1} \biggr)^2 \gamma^{\alpha - 2} \geq \frac{1}{C} \frac{a_1}{a_2^2} \gamma^{\alpha - 2}.
\end{align*}
The upper bound on $V(\gamma)$ follows in an identical manner. The conclusion of the Lemma then follows as desired.

\end{proof}

\begin{lemma}
\label{lem:tildetheta_moment}
Let $Z$ be a random variable on $[-1, 1]$ with a distribution symmetric around $0$ and write $\nu_\gamma := \mathbb{E}|Z|^{\gamma}$ for $\gamma \geq 1$. Suppose $a_1 \gamma^{-\alpha} \leq \nu_\gamma \leq a_2 \gamma^{-\alpha}$ for all $\gamma \geq 1$ and for some $\alpha \in (0, 2]$, $a_1 \in [0, 1]$ and $a_2 \in [1, \infty)$. Then, for any $\gamma \geq 1$ and any $0 \leq \Delta \leq \frac{1}{4\gamma}$, we have 
\begin{align*}
\mathbb{E}|Z - \Delta|^\gamma \leq C a_2 \gamma^{-\alpha}.
\end{align*}

Moreover, we have that for any $\gamma \geq 2$ and any $\Delta \in \mathbb{R}$, 
\begin{align*}
\mathbb{E}\bigl[ - |Z - \Delta|^{\gamma - 1} \text{sgn}(Z - \Delta) \bigr] \geq  \frac{a_1}{2}  |\Delta| \gamma^{1-\alpha}.
\end{align*}

Lastly, for any $k \in \mathbb{N}$ and any $\Delta_{\gamma}$ (allowed to depend on $\gamma$) such that $0 \leq \Delta_{\gamma} \leq \frac{1}{4 \gamma}$, we have
\begin{align*}
\mathbb{E} \biggl[ \sup_{\gamma \in [2^k, 2^{k+1}]} (|Z| + \Delta_{\gamma} )^{2 (\gamma - 1)} \biggr] \leq 
C a_2 2^{ - k \alpha}. 
\end{align*}

\end{lemma}

\begin{proof}
Consider the first claim. Observe that
\begin{align*}
\mathbb{E}|Z - \Delta|^\gamma = 
\underbrace{\mathbb{E}\biggl[ |Z - \Delta|^\gamma \mathbbm{1}\{ |Z| \leq 1/4 \} \biggr]}_{\text{Term 1}} + 
\underbrace{\mathbb{E}\biggl[ |Z - \Delta|^\gamma \mathbbm{1}\{ |Z| > 1/4 \} \biggr]}_{\text{Term 2}}.
\end{align*}

To bound Term 1, we have that
\[
|Z - \Delta|^\gamma \mathbbm{1}\{ |Z| \leq 1/4\} \leq 2^{-\gamma} \leq 2 \gamma^{-\alpha},
\]
where, in the last inequality, we use the fact that $\alpha \in (0, 2]$ and that $2^{-x} \leq 2 x^{-2}$ for all $x \geq 1$. It is clear then that Term 1 is bounded by $2 \gamma^{-\alpha}$. To bound Term 2, we have that
\begin{align*}
|Z - \Delta|^\gamma \mathbbm{1}\{ |Z| > 1/4\} 
&= |Z|^\gamma \biggl| 1 - \frac{\Delta}{Z} \biggr|^\gamma \mathbbm{1} \{ |Z| > 1/4 \} \\
&\leq |Z|^\gamma \biggl| 1 + \frac{1}{\gamma} \biggr|^\gamma \mathbbm{1} \{ |Z| > 1/4 \} \leq e | Z |^\gamma,
\end{align*}
where in the second inequality, we use the fact that $\Delta \leq \frac{1}{4 \gamma}$. Therefore, we have that
\[
\mathbb{E}\bigl[ |Z - \Delta|^\gamma \mathbbm{1}\{ |Z| < 1/4 \} \bigr] \leq e \mathbb{E}|Z|^{\gamma} \leq C a_2 \gamma^{-\alpha}.
\]
Combining the bounds on the two terms, we have that $\mathbb{E}|Z - \Delta|^\gamma \leq C a_2 \gamma^{-\alpha}$ as desired.

We now turn to the second claim. Without loss of generality, assume that $\Delta \geq 0$ so that, by symmetry of the distribution of $Z$, we have $\mathbb{E}\bigl[ - |Z - \Delta|^{\gamma - 1} \text{sgn}(Z - \Delta) \bigr] \geq 0$.

Since $\mathbb{E} \bigl[ -  |Z|^{\gamma - 1} \text{sgn}(Z) \bigr] = 0$, 
\begin{align*}
\mathbb{E}\bigl[ - |Z - \Delta|^{\gamma - 1} \text{sgn}(Z - \Delta) \bigr] 
&= \int_0^\Delta (\gamma - 1) \mathbb{E}\bigl[ |Z - t|^{\gamma - 2} \bigr] \, dt \\
&\geq |\Delta| (\gamma - 1)  \mathbb{E}\bigl[ |Z|^{\gamma - 2} \bigr] 
\end{align*}

For $\gamma \in [2, 3)$, it holds that $\mathbb{E} \bigl[ |Z|^{\gamma - 2} \bigr] \geq \mathbb{E} |Z| \geq a_1$ since $Z$ is supported on $[-1, 1]$. For $\gamma \geq 3$, it holds that $\mathbb{E} \bigl[ |Z|^{\gamma - 2} \bigr] = \nu_{\gamma - 2} \geq a_1 (\gamma - 2)^{-\alpha}$. Therefore, we have that
\begin{align*}
\mathbb{E}\bigl[ - |Z - \Delta|^{\gamma - 1} \text{sgn}(Z - \Delta) \bigr] 
&\geq  \begin{cases} a_1 |\Delta| (\gamma - 1) & \text{ if $\gamma \in [2, 3)$,} \\
a_1 | \Delta | (\gamma - 1) (\gamma - 2)^{-\alpha} & \text{ else.}
\end{cases} 
\end{align*}
Thus, for all $\gamma \geq 2$, we have that
\begin{align*}
\mathbb{E}\bigl[ - |Z - \Delta|^{\gamma - 1} \text{sgn}(Z - \Delta) \bigr] \geq \frac{a_1}{2} |\Delta| \gamma^{1 - \alpha}.
\end{align*}


Finally, we consider the third claim. The argument is similar to that of the first claim. We observe that
\begin{align}
\mathbb{E} \biggl[ \sup_{\gamma \in [2^k, 2^{k+1}]} (|Z| + \Delta_{\gamma})^{2(\gamma-1)} \biggr] 
&= \int_0^{\frac{1}{4}} \sup_{\gamma \in [2^k, 2^{k+1}]} (z+ \Delta_{\gamma})^{2 (\gamma-1)}\, dP(z) \nonumber \\
&\qquad \qquad + \int_{\frac{1}{4}}^1 \sup_{\gamma \in [2^k, 2^{k+1}]} (z + \Delta_{\gamma})^{2 (\gamma-1)}\, dP(z). \label{eq:zsup_decomp}
\end{align}

To bound the first term of~\eqref{eq:zsup_decomp}, we use the fact that $\Delta_{\gamma} \leq \frac{1}{4\gamma} \leq \frac{1}{4}$ and that $\alpha \in (0, 2]$ to obtain
\begin{align*}
\int_0^{\frac{1}{4}} \sup_{\gamma \in [2^k, 2^{k+1}]} (z + \Delta_{\gamma})^{2 (\gamma-1)} dP(z) 
\leq 2^{ - 2( 2^k - 1) } \leq 2^{ - k \alpha}. 
\end{align*}

To bound the second term of~\eqref{eq:zsup_decomp}, we have
\begin{align*}
\int_{\frac{1}{4}}^1 \sup_{\gamma \in [2^k, 2^{k+1}]} (z + \Delta_{\gamma})^{2 (\gamma-1)} dP(z) 
&\leq \int_{\frac{1}{4}}^1 \sup_{\gamma \in [2^k, 2^{k+1}]} z^{2(\gamma-1)} (1 + 4\Delta_{\gamma})^{2 (\gamma-1)} dP(z) \\
&\leq e^2 \mathbb{E} |Z|^{2(2^k - 1)} \leq C a_2 2^{- k \alpha}. 
\end{align*}
The third claim of the lemma thus follows as desired. 
\end{proof}

\begin{lemma}
\label{lem:conti_argmin_L}
Define $L(\gamma,\mathbf{P}) := \frac{1}{\gamma} \min_{\theta}\log \bigl(\int |y - \theta|^\gamma \mathbf{P}(dy) \bigr) + \frac{1+\log \gamma}{\gamma} + \log \Gamma\biggl( 1 + \frac{1}{\gamma} \biggr)$ for every $\gamma \geq 2$. Given $\lim_{\gamma\to\infty}L(\gamma,\mathbf{P}_1)=\lim_{\gamma\to\infty}L(\gamma,\mathbf{P}_2)=\infty$, $\gamma_1^*$ being the unique minimizer of $L(\gamma,\mathbf{P}_1)$, and $L(\gamma^*_1,\mathbf{P}_2)<\infty$, we have that $\gamma_1^*$ is the unique minimizer of $L(\gamma,(1-\delta)\mathbf{P}_1+\delta\mathbf{P}_2)$ for all small positive $\delta$.
\end{lemma}
\begin{proof}
    We first show that $\lim_{\gamma\to\infty}\inf_{0\leq\delta\leq1}L(\gamma,(1-\delta)\mathbf{P}_1+\delta\mathbf{P}_2)=\infty$. Given $M>0$, there exists a $N\in \mathbb{N}$ such that $\frac{1}{\gamma} \min_{\theta}\log \left(\int |y - \theta|^\gamma \mathbf{P}_1(dy) \right)\vee \frac{1}{\gamma} \min_{\theta}\log \left(\int |y - \theta|^\gamma \mathbf{P}_2(dy) \right)>M$ for every $\gamma>N$, and thus
    \begin{align*}
        L(\gamma,(1-\delta)\mathbf{P}_1+\delta\mathbf{P}_2)&\geq\frac{1}{\gamma} \min_{\theta}\log \left[\int |y - \theta|^\gamma ((1-\delta)\mathbf{P}_1+\delta\mathbf{P}_2)(dy) \right]\\
        &\geq\frac{1}{\gamma} \min_{\theta}\left[(1-\delta)\log \int |y - \theta|^\gamma \mathbf{P}_1(dy)+\delta\log \int |y - \theta|^\gamma \mathbf{P}_2(dy)\right]\\
        & \geq(1-\delta)\frac{1}{\gamma} \min_{\theta}\log \left(\int |y - \theta|^\gamma \mathbf{P}_1(dy) \right)+\delta\frac{1}{\gamma} \min_{\theta}\log \left(\int |y - \theta|^\gamma \mathbf{P}_2(dy) \right)\\
        &\geq M, \text{for every }\gamma>N.
    \end{align*}
For a fixed $\gamma \geq 2$, we have
\begin{align*}
    \lim_{\delta\to0^+}L(\gamma,(1-\delta)\mathbf{P}_1+\delta\mathbf{P}_2)=\begin{cases}
    L(\gamma,\mathbf{P}_1),& \text{if } L(\gamma,\mathbf{P}_2)<\infty\\
    \infty,              & \text{otherwise.}
\end{cases}
\end{align*}
\end{proof}

\subsection{Reference results}

We use the following statement of Talagrand's inequality:

\begin{theorem} 
\label{thm:talagrand} (Talagrand's Inequality; see e.g.~\citet[][Theorem 3.3.9]{gine2021mathematical})
Let $Z_1, \ldots, Z_n$ be independent and identically distributed random objects taking value on some measurable space $\mathcal{Z}$. Let $\mathcal{F}$ be a class of real-valued Borel measurable functions on $\mathcal{Z}$. 

Define $S_n = \sup_{f \in \mathcal{F}} \sum_{i=1}^n \bigl\{ f(Z_i) - \mathbb{E} f(Z) \bigr\} $. Let $U > 0$ be a scalar that $ \sup_{f \in \mathcal{F}} | f(Z)| \leq U$ almost surely; let $\sigma^2 := \sup_{f \in \mathcal{F}} \mathbb{E} f^2(Z)$. Then, for any $t > 0$, 
\begin{align*}
\mathbb{P}( S_n - \mathbb{E} S_n \geq t) \leq \exp \biggl\{ - \frac{t^2}{ 2 U \cdot \mathbb{E} S_n + n \sigma^2 } \wedge \frac{t}{ \frac{2}{3} U} \biggr\}. 
\end{align*}
\end{theorem}

We use the following bound on the expected supremum of the empirical process. For a class of real-valued functions $\mathcal{F}$ on some measurable domain $\mathcal{Z}$, we write $F(z) := \sup_{f \in \mathcal{F}} |f(z)|$ as its envelope function. For $\delta \in [0, 1)$, define the entropy integral
\begin{align}
    J(\delta) \equiv J(\delta, \mathcal{F}):=\int_0^\delta \sup_{Q}\sqrt{\log \mathcal{N}(\epsilon\|F\|_{L_2(Q)},\mathcal{F},L_2(Q))}\,d\epsilon,
    \label{eq:entropy_integral}
\end{align}
where the supremum is taken over all finitely discrete probability measures.

\begin{lemma} \cite[][Theorem 2.6.7]{van1996weak}
If $\mathcal{F}$ has finite VC dimension $V(\mathcal{F}) \geq 2$, then, for any $\epsilon \in (0, 1)$,
\begin{align*}
N( \epsilon \| F \|_{L_2(Q)}, \mathcal{F}, L_2(Q)) \leq C V(\mathcal{F}) (16e)^{V(\mathcal{F})} \biggl( \frac{1}{\epsilon} \biggr)^{2 (V(\mathcal{F}) - 1) }.
\end{align*}
\label{lem:VC}
\end{lemma}

\begin{corollary}
\label{cor:VC}
If $\mathcal{F}$ has finite VC dimension $V(\mathcal{F})$, then, for any $\delta \in (0, 1]$,
\[
J(\delta) \leq C \sqrt{ V(\mathcal{F}) }  \delta \sqrt{ \log \frac{1}{\delta} \wedge 1 } .
\]
\end{corollary}

\begin{proof}
Using Lemma~\ref{lem:VC}, we have that
\begin{align*}
J(\delta) &\leq \int_0^{\delta} \sqrt{ C V(\mathcal{F}) + 2( V(\mathcal{F}) - 1) \log \frac{1}{\epsilon} } \, d \epsilon \\
&\leq C \sqrt{V(\mathcal{F})} \biggl\{  \delta  + \int_0^{\delta} \sqrt{ \log \frac{1}{\epsilon}} \, d \epsilon \biggr\} \leq C \sqrt{ V(\mathcal{F}) } \biggl( \delta \sqrt{ \log \frac{1}{\delta} \wedge 1 } \biggr).
\end{align*}
\end{proof}


\begin{theorem}[\cite{van2011local}]
\label{thm:exp_sup_local}
Let $F(x):= \sup_{f\in\mathcal{F}}|f(x)|$, $M:=\max_{1\leq i\leq n}F(Z_i)$, and $\sigma^2:=\sup_{f\in\mathcal{F}}\mathbb{E} f(Z)^2$. Then the following two bounds hold:
\begin{align*}
    \mathbb{E} \sup_{f\in\mathcal{F}}\biggl|\frac{1}{n}\sum_{i=1}^n f(Z_i) - \mathbb{E}f(Z) \biggr|
    \lesssim\frac{\|F\|_{L_2(P)}}{\sqrt{n}}
    J\left(\frac{\sigma}{\|F\|_{L_2(P)}}\right)\left[1+\frac{\|M\|_{L_2(P)}\|F\|_{L_2(P)}J\left(\frac{\sigma}{\|F\|_{L_2(P)}}\right)}{\sqrt{n}\sigma^2}\right]
\end{align*}
as well as
\begin{align*}
\mathbb{E} \sup_{f\in\mathcal{F}}\biggl|\frac{1}{n}\sum_{i=1}^n f(Z_i) - \mathbb{E}f(Z) \biggr| \lesssim 
\frac{\| F \|_{L_2(P)} J(1)}{\sqrt{n}}.
\end{align*}

\end{theorem}

\end{spacing}

\end{document}